\pdfoutput=1
\documentclass[11pt]{article}
\usepackage{xcolor}
\usepackage{graphicx}
\usepackage{amsmath,amsfonts,amssymb,graphics,amsthm}
\usepackage{hyperref}
\usepackage{comment}
\usepackage{tabularx}
\usepackage[protrusion=true,expansion=true]{microtype}
\usepackage{enumerate}
\usepackage{bbm}
\usepackage{mathrsfs}
\usepackage[margin=1in]{geometry}
\usepackage{csquotes}
 \usepackage[normalem]{ulem}

\hypersetup{
    colorlinks=false,
    linktocpage,
    }

\numberwithin{equation}{section}

\newtheorem{theorem}{Theorem}[section]
\newtheorem{corollary}[theorem]{Corollary}
\newtheorem{lemma}[theorem]{Lemma}
\newtheorem{proposition}[theorem]{Proposition}

\newtheorem{lem}[theorem]{Lemma} 

\newtheorem{definition}[theorem]{Definition}

\theoremstyle{remark}

\theoremstyle{definition}
\newtheorem{remark}[theorem]{Remark}

\def\@rst #1 #2other{#1}
\newcommand\MR[1]{\relax\ifhmode\unskip\spacefactor3000 \space\fi
  \MRhref{\expandafter\@rst #1 other}{#1}}
\newcommand{\MRhref}[2]{\href{http://www.ams.org/mathscinet-getitem?mr=#1}{MR#2}}

\def\MR#1{\href{http://www.ams.org/mathscinet-getitem?mr=#1}{MR#1}}

\newcommand{\C}{\mathbbm{C}}
\newcommand{\D}{\mathbbm{D}}
\newcommand{\E}{\mathbbm{E}}
\newcommand{\N}{\mathbbm{N}}

\newcommand{\R}{\mathbbm{R}}
\renewcommand{\P}{\mathbbm{P}}
\newcommand{\bbH}{\mathbbm{H}}
\newcommand{\IG}{\mathrm{IG}}

\newcommand{\eps}{\varepsilon}
\newcommand{\1}{\mathbf{1}}

\newcommand{\sph}{\mathrm{sph}}
\newcommand{\disk}{\mathrm{disk}}
\newcommand{\wed}{\mathrm{wedge}}
\newcommand{\cone}{\mathrm{cone}}

\newcommand{\Leb}{\mathrm{Leb}}

\newcommand{\cc}{\textbf{c}}
\newcommand{\cMtwo}{{\cM^{\textrm{disk}}_2}}

\let\Re\undefined
\DeclareMathOperator{\Re}{Re}

\DeclareMathOperator{\Cov}{Cov}
\DeclareMathOperator{\Var}{Var}
\DeclareMathOperator{\SLE}{SLE}

\newcommand{\QD}{\mathrm{QD}}

\def\cX{\mathcal{X}}
\def\cW{\mathcal{W}}

\def\cS{\mathcal{S}}

\def\cP{\mathcal{P}}

\def\cM{\mathcal{M}}
\def\cL{\mathcal{L}}
\def\cK{\mathcal{K}}

\def\cH{\mathcal{H}}

\def\cF{\mathcal{F}}

\def\cD{\mathcal{D}}
\def\cC{\mathcal{C}}
\def\cB{\mathcal{B}}
\def\cA{\mathcal{A}}

\def\alb#1\ale{\begin{align*}#1\end{align*}}
\def\allb#1\alle{\begin{align}#1\end{align}}

\newcommand{\aryb}{\begin{eqnarray*}}
\newcommand{\arye}{\end{eqnarray*}}
\def\alb#1\ale{\begin{align*}#1\end{align*}}
\newcommand{\eqb}{\begin{equation}}
\newcommand{\eqe}{\end{equation}}
\newcommand{\eqbn}{\begin{equation*}}
\newcommand{\eqen}{\end{equation*}}

\newcommand{\BB}{\mathbbm}
\newcommand{\ol}{\overline}

\newcommand{\op}{\operatorname}

\newcommand{\rta}{\rightarrow}

\newcommand{\wt}{\widetilde}
\newcommand{\wh}{\widehat} 
\newcommand{\mcl}{\mathcal}

\let\originalleft\left
\let\originalright\right
\renewcommand{\left}{\mathopen{}\mathclose\bgroup\originalleft}
\renewcommand{\right}{\aftergroup\egroup\originalright}

\DeclareMathAlphabet{\mathpzc}{OT1}{pzc}{m}{it}

\begin{document}

\title{Conformal welding of quantum disks}
\author{
\begin{tabular}{c}Morris Ang\\[-5pt]\small MIT\end{tabular}\; 
\begin{tabular}{c}Nina Holden\\[-5pt]\small ETH Z\"urich\end{tabular}\; 
\begin{tabular}{c}Xin Sun\\[-5pt]\small University of Pennsylvania\end{tabular}
} 
\date{  }

\maketitle

\begin{abstract}
Two-pointed quantum disks with a weight parameter $W>0$ are a family of finite-area random surfaces that arise naturally in Liouville quantum gravity. In this paper we show that conformally welding two quantum disks according to their boundary lengths gives another quantum disk decorated with an independent chordal $\SLE_\kappa(\rho_-; \rho_+)$ curve. This is the finite-volume counterpart of 
the classical result of Sheffield (2010) and Duplantier-Miller-Sheffield (2014) on the welding of infinite-area two-pointed quantum surfaces called quantum wedges, which is fundamental to the mating-of-trees theory. Our results can be used to give unified proofs of the mating-of-trees theorems for the quantum disk and the quantum sphere, in addition to a mating-of-trees description of the  weight $W=\frac{\gamma^2}{2}$ quantum disk.
Moreover, it serves as a key ingredient in our companion work~\cite{AHS-SLE-integrability}, which proves an exact formula
for $\SLE_\kappa(\rho_-; \rho_+)$  using conformal welding of  random surfaces and a conformal welding result giving the so-called SLE loop.
\end{abstract}

\section{Introduction}
Liouville quantum gravity (LQG) is a theory of random surfaces  with close connections to conformal field theory and random planar maps \cite{polyakov-qg1, david-conformal-gauge, dk-qg}.
For $\gamma\in(0,2)$, the random area measure of a  $\gamma$-LQG surface is of the form $e^{\gamma h}d^2z$ where $h$ is a variant of Gaussian free field and $d^2z$ is the Euclidean area measure.
Although $h$ is only a Schwartz distribution which is not pointwise defined, the  area measure $e^{\gamma h}d^2z$ can be understood by regularizing $h$ and taking a renormalized limit~\cite{shef-kpz}.  This construction falls into the general framework of Gaussian multiplicative chaos; see \cite{kahane,rhodes-vargas-review}. 
Recently the metric associated with LQG surfaces was also rigorously constructed  by regularizing the field~\cite{dddf-lfpp,gm-uniqueness}.

Quantum wedges are a natural family of  infinite-area $\gamma$-LQG surfaces with two marked points on the boundary.
Neighborhoods of one point have finite $\gamma$-LQG area, whereas neighborhoods of the other one have infinite $\gamma$-LQG area. A quantum wedge is associated with a \emph{weight} parameter $W>0$ which describes the singularity at the two marked points.

\begin{figure}[ht!]
	\begin{center}
		\includegraphics[scale=0.6]{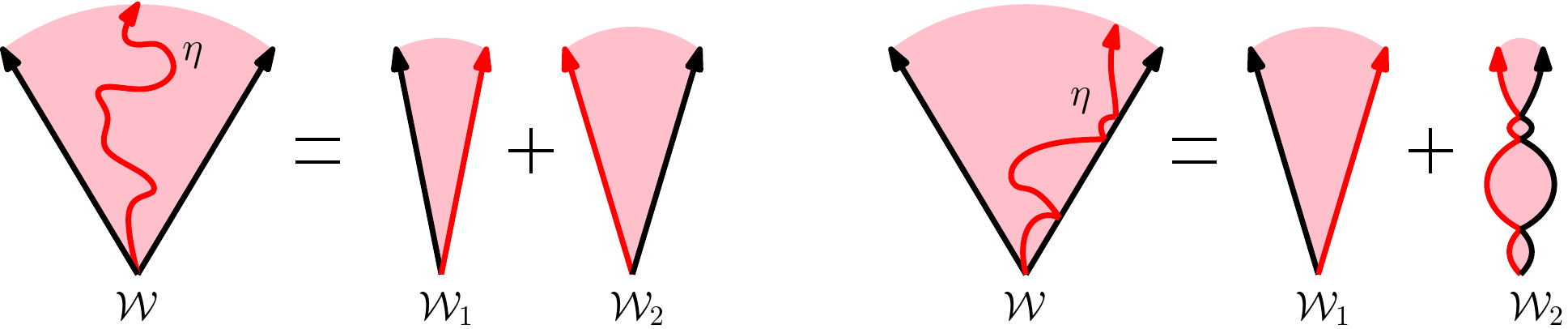}%
	\end{center}
	\caption{\label{fig-front} The conformal welding of quantum wedges \cite{wedges}. \textbf{Left.} The case $W_1, W_2 \geq \frac{\gamma^2}2$, so each quantum wedge has the disk topology. \textbf{Right.} The case $W_1 \geq \frac{\gamma^2}2>W_2$, so the second quantum wedge is a chain of disks. Not illustrated is 
		the case where $W_1, W_2 < \frac{\gamma^2}2$  and $W_1+W_2 \geq \frac{\gamma^2}2$and the case where $W_1+W_2 < \frac{\gamma^2}2$.}
\end{figure}

A particularly fruitful approach to studying LQG is through its coupling with Schramm-Loewner evolutions (SLE), which are an important family of conformally invariant random  planar curves associated with a parameter $\kappa>0$ \cite{schramm0}.
A key LQG/SLE coupling result is the \emph{conformal welding of quantum wedges}.
For $\rho_1,\rho_2>-2$, $\SLE_\kappa(\rho_1; \rho_2)$ is a variant of $\SLE_\kappa$; see Section~\ref{subsec-SLE}. The following result was proved in \cite{wedges}, see Figure \ref{fig-front} for an illustration.
\begin{displayquote}
Set $\kappa = \gamma^2$. For $W_1, W_2 >0$, a weight $(W_1 + W_2)$ $\gamma$-LQG quantum wedge $\cW$ cut by an independent $\SLE_\kappa(W_1 - 2; W_2 - 2)$ curve $\eta$ yields two \emph{independent} $\gamma$-LQG quantum wedges $\cW_1$ and $\cW_2$ of weights $W_1$ and $W_2$, respectively. Moreover, $(\cW, \eta)$ is measurable with respect to the quantum surfaces $(\cW_1, \cW_2)$. (See Theorem~\ref{thm-wedges} for the full statement.)
\end{displayquote}
 The conformal welding result for quantum wedges is arguably one of the deepest facts in random planar geometry. It was proved by Sheffield in~\cite{shef-zipper} when $W_1=W_2=2$ and generalized in \cite{wedges}. %
It  is a key input to the \emph{mating-of-trees} theory of Duplantier, Miller, and Sheffield~\cite{wedges},
which is a powerful framework to study SLE and LQG via Brownian motion, and is fundamental to the link between LQG and the scaling limits of random planar maps. See~\cite{ghs-mating-survey} for a survey. 

\begin{figure}[ht!]
	\begin{center}
		\includegraphics[scale=0.65]{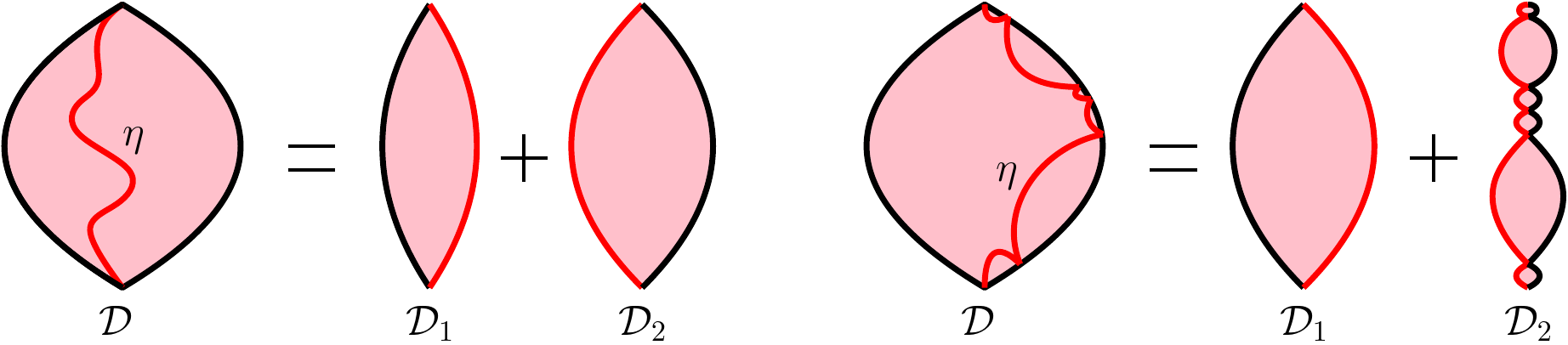}%
	\end{center}
	\caption{ \label{fig-front-disk} Theorem~\ref{thm-disk-cutting-2} describes the conformal welding of quantum disks. The cases $W_1, W_2 \geq \frac{\gamma^2}2$ (left) and $W_1 \geq \frac{\gamma^2}2 > W_2$ (right) are illustrated here. }
\end{figure}

For each weight parameter $W \geq \frac{\gamma^2}2$ there is also an infinite measure on $\gamma$-LQG surfaces with finite $\gamma$-LQG area called the \emph{(two-pointed) quantum disk of weight $W$}~\cite{wedges}. Quantum disks can be considered the finite-area analog of quantum wedges, and they also have the topology of a disk with two boundary marked points. 
We extend the definition of the quantum disk to $W \in (0, \frac{\gamma^2}2)$ in Section~\ref{subsec-thin-disks}, and we view them
as the finite-area analog of quantum wedges of weight $W$.
In this regime, the topology  of quantum wedges and disks is given by  a chain of  countably many disks; see Figures~\ref{fig-front} (right) and~\ref{fig-front-disk} (right) for an illustration.

The main result of this paper, Theorem~\ref{thm-disk-cutting-2}, is the conformal welding of quantum disks, which can be informally stated as follows; see Figure \ref{fig-front-disk}.
  \begin{displayquote}
  	Set $\kappa = \gamma^2$. 
  For $W_1, W_2 >0$, a weight $(W_1 + W_2)$ $\gamma$-LQG quantum disk $\cD$ cut by an independent $\SLE_\kappa(W_1 - 2; W_2 - 2)$ curve $\eta$ yields two quantum disks $\cD_1, \cD_2$ which are conditionally independent  given the $\gamma$-LQG length $\ell$ of $\eta$; the conditional law of $\cD_1$ (resp.\ $\cD_2$) is a weight $W_1$ (resp.\ $W_2$) quantum disk conditioned on having right (resp.\ left) boundary arc of length $\ell$.
  	Moreover, $(\cD, \eta)$ is measurable with respect to the quantum surfaces $(\cD_1, \cD_2)$. 
  \end{displayquote}
 We similarly show in Theorem~\ref{thm-sphere-cutting} that cutting a \emph{quantum sphere} by a certain SLE-type curve yields a quantum disk. Quantum spheres are quantum surfaces with the topology of the two-pointed sphere. This result is the finite-area analog of \cite[Theorem 1.4]{wedges}, which states that a \emph{quantum cone} cut by a certain SLE-type curve results in a quantum wedge. 
Using~\cite{mmq-welding}, the conformal welding of weight 2 quantum wedges was extended  to the critical case $\gamma = 2$ and $\kappa=4$ in \cite{hp-welding}.  We believe our results extend to $\gamma=2$ via similar considerations; see Remark \ref{rmk:critical}.
 
Our proof relies on the intuition that the quantum disk can be obtained from a quantum wedge by creating and pinching a suitable bottleneck. Using our approach,
the mating-of-trees theorems for the quantum sphere and disk can be easily deduced, as we sketch in Section~\ref{subsec-MOT-rederive}.
These results were originally proved in  \cite{sphere-constructions, ag-disk}.
Although the original proofs  are also based on pinching infinite area LQG surfaces, our treatment is unified and conceptually more straightforward.
Moreover, in Section~\ref{subsec-MOT} we derive the area and boundary length distribution of a weight-$\frac{\gamma^2}2$ quantum disk using mating-of-trees, properties of Brownian motion, and the main result of this paper.

Our paper is a key ingredient of several concurrent works.
 In our companion paper \cite{AHS-SLE-integrability}, we prove an exact result for  $\SLE_\kappa(\rho_-; \rho_+)$  
 and establish a quantum zipper result for the SLE loop~\cite{zhan-loop-measures} on the quantum sphere. 
 Both of the two results crucially rely on the conformal welding result proved in this paper.
 In the joint work of the first and third authors with Remy~\cite{ARS-FZZ}, another  conformal welding result is proved based on  our result to prove the so-called FZZ formula in Liouville conformal field theory (LCFT). More generally, conformal welding of finite-area quantum surfaces is a cornerstone of the ongoing program of the first and third authors 
 proving exact results for SLE, LCFT and mating-of-trees by exploring their connections; see~\cite{AS-CLE} for another example.

It has been shown in various senses that random planar maps weighted by certain statistical physics models converge to certain $\gamma$-LQG random surfaces, where $\gamma$ depends on the choice of model. 
For instance, random planar maps with the disk topology decorated by an FK cluster model, Potts model, or $O(n)$-loop model with monochromatic boundary conditions should converge to $\gamma$-LQG quantum disks with weight 2 for some $\gamma \in (\sqrt2,2)$ in the scaling limit. This was first demonstrated by the pioneering work of \cite{shef-burger} for the FK cluster model. See \cite{ghs-mating-survey} for a comprehensive review on the relation between LQG and random planar maps.
Applying our results to $W_1 = W_2 = 2$ implies that if the boundary condition is Dobrushin (rather than monochromatic), then the limiting surface should be a weight 4 quantum disk. Moreover, the natural chordal interface associated with the Dobrushin boundary condition should converge to $\SLE_{\kappa}$ with $\kappa = \gamma^2$. 
In the sense of metric geometry for $\gamma = \sqrt{\frac83}$, this follows from the work of Gwynne and Miller~\cite{gm-disk}, where the decorating model is the self avoiding walk; see Remark~\ref{rem-GM}. For Ising-weighted maps ($\gamma =\sqrt3$), \cite{ct-ising-length} proves interesting results consistent with this picture; see Remark~\ref{rem-interface-length}.

\bigskip
\noindent{\bf Proof ideas.} 
We derive our result from its counterpart for quantum wedges.
The crucial step is to  define a proper ``bottleneck'' around the origin of a weight $W > \frac{\gamma^2}2$ quantum wedge and, roughly speaking, condition on the bottleneck being small and the pinched region being large; under such conditioning, the pinched region becomes close in some sense to a weight $W$ quantum disk. If one then cuts the weight $W$ quantum wedge (using an $\SLE_\kappa(W_1-2;W_2-2)$ curve) into two quantum wedges of weights $W_1, W_2$, one expects that each of these is pinched to get quantum disks of weights $W_1, W_2$, as desired.  

While the high level picture is clear, a direct implementation of this argument seems exceedingly difficult when $W_1, W_2 \geq \frac{\gamma^2}2$ because it is hard to define a tractable bottleneck event which pinches all three quantum wedges $\cW_1, \cW_2, \cW$ to yield quantum disks. For instance, if one defines a bottleneck for each quantum wedge $\cW_1, \cW_2$, then these bottlenecks together should be a bottleneck for $\cW$, but the analysis of this bottleneck on $\cW$ must consider the conformal welding of $\cW_1, \cW_2$. 

To resolve this, the key insight is our new definition of quantum disks with weight $W' < \frac{\gamma^2}2$. Weight $W'$ quantum wedges are defined as an infinite Poissonian chain of weight $(\gamma^2-W')$ quantum disks, and we define a weight $W'$ quantum disk as a finite truncation of this chain. Consequently, it is easy to ``pinch'' a weight $W'$ quantum wedge to obtain a weight $W'$ quantum disk, and this enables us to define a tractable bottleneck for the above proof sketch when $W_1, W_2 < \frac{\gamma^2}2$ and $W_1 + W_2 > \frac{\gamma^2}2$. 

The same argument shows that for $W_1, \dots, W_n \in(0, \frac{\gamma^2}2)$ with $W = \sum W_j > \frac{\gamma^2}2$, cutting a weight $W$ quantum disk by a certain collection of SLE-type curves yields a collection of quantum disks with weights $W_1, \dots, W_n$.  Soft arguments then allow us to remove the weight restrictions, yielding the full theorem.
\medskip

\noindent\textbf{Paper outline. }  We explain preliminaries and state our main results in Section~\ref{sec-main-results}, and prove our main results in Sections~\ref{sec-extrinsic}--\ref{sec-general}. In Section~\ref{sec-MOT} we give alternative proofs of finite area mating-of-trees theorems, in addition to giving a novel mating-of-trees representation of the weight $W=\frac{\gamma^2}{2}$ quantum disk. %
A more detailed overview of Sections~\ref{sec-extrinsic}--\ref{sec-general} can be found at the end of Section~\ref{sec-main-results}.
\medskip 

\noindent\textbf{Acknowledgements. } We thank Yilin Wang for inspiring discussions on the conformal welding of quantum disks, and
Guillaume Baverez, Ewain Gwynne and Jason Miller for helpful comments on the first version of the paper. 
 M.~A. was supported by NSF grant DMS-1712862. N.~H. was supported by Dr. Max Rossler, the Walter Haefner Foundation,
and the ETH Zurich Foundation. X.~S. was supported by a Junior Fellow award from the Simons Foundation, and NSF Grant DMS-1811092 and DMS-2027986.

\section{Definition of LQG surfaces and statement of the main results}
\label{sec-main-results}

The main goal of this section is to precisely state our main results and give sufficient background to make these statements. 
In Section~\ref{subsec-prelim} we give some preliminaries, and then we state the main results in Section~\ref{subsec-main-results}.
	Some quantum surfaces and curves are only discussed at high level in Section~\ref{subsec-main-results}, and the rest of this section is devoted to introducing these random objects. In Section~\ref{subsec-wedges} we define quantum wedges, sphere and cones, in Section~\ref{subsec-thin-disks} we introduce the thin quantum disk, and in Section~\ref{subsec-basic} we explain some basic properties of quantum disks. In Section~\ref{subsec-thick-disint} we carry out the disintegrations of quantum disks with respect to boundary arc lengths. In Section~\ref{subsec-SLE} we explain some SLE preliminaries. Finally we give an outline for the rest of the paper in Section~\ref{subsec-outline}.

\subsection{Preliminaries}\label{subsec-prelim}

We define the Neumann Gaussian free field (GFF) on the strip $\cS := \R \times (0,\pi)$; the definition extends to other domains by conformal invariance. 
With slight abuse of notation we sometimes consider $\cS$ as a subset of $\R^2$ and other times of $\C$, so for instance $\{0\} \times [0,\pi]$ is also written as $[0,i\pi]$. We also write $\R_+ := (0, \infty)$, $\R_- := (-\infty, 0)$ and $\cS_\pm := \R_\pm \times (0,\pi)$.

Consider the space of smooth functions on $\cS$ with bounded support and mean zero on $[0,i\pi]$, and define the Dirichlet inner product 
\[(f,g)_\nabla = \frac1{2\pi}\int_{\cS}\nabla f(z) \cdot \nabla g(z) \, d^2z. \]
Let $H(\cS)$ be the Hilbert space closure of this space with respect to $(\cdot, \cdot)_\nabla$. Then the \emph{Neumann GFF on $\cS$ normalized to have mean zero on $[0,i\pi]$} is the random distribution 
\[h = \sum_{i = 1}^\infty \alpha_i f_i, \]
where $(\alpha_i)_{i=1}^\infty$ are i.i.d.\ standard Gaussians and $(f_i)_{i=1}^\infty$ is an orthonormal basis for $H(\cS)$; one can show the law of $h$ does not depend on the choice of $(f_i)_{i=1}^\infty$. The above summation does not converge in $H(\cS)$, but a.s.\ converges in the space of distributions \cite[Section 4.1.4]{wedges}.

Define $\cH_1(\cS) \subset H(\cS)$ (resp. $\cH_2(\cS) \subset H(\cS)$) to be the space of functions which are constant (resp. have mean zero) on every vertical segment $[t, t+i\pi]$ for $t \in \R$. Then $H(\cS) = \cH_1(\cS) \oplus \cH_2(\cS)$ is an orthogonal decomposition, and hence the projections of $h$ to $\cH_1(\cS)$ and $\cH_2(\cS)$ are independent. See \cite[Section 4.1.6]{wedges} for more details. In this paper, we mainly consider generalized functions which are GFFs plus (possibly random) continuous functions; we call these fields. 
For a field $\psi$ on $\cS$, we write $\psi_t$ for the average of $\psi$ on $[t,t+i\pi]$, and identify the projection of $\psi$ to $\cH_1(\cS)$ with the function  $(\psi_t)_{t \in \R}$.

In this paper, we will always consider LQG with parameter $\gamma \in (0, 2)$, and write $Q = \frac\gamma2 + \frac2\gamma$. We will often keep the dependences on $\gamma$ implicit for notational simplicity. 
Let 
\[\cD \cH := \{ (D,h) \: : \: D \subset \C \text{ is open},\,h \text{ is a distribution on }D\}.\]
We will typically take $h$ to be a variant of the GFF.
For $(D,h), (\wt D, \wt h) \in \cD \cH$, we say that $(D,h) \sim_\gamma (\wt D, \wt h)$ if there exists a conformal map $\varphi: \wt D \to D$ such that 
\eqb \label{eq-quantum-surface}
\wt h = h \circ \varphi + Q \log |\varphi'|.
\eqe
For $\gamma \in (0,2)$, a \emph{$\gamma$-LQG surface} (or quantum surface) is an equivalence class of pairs $(D,h) \in \cD\cH$ under the equivalence relation $\sim_\gamma$, and an \emph{embedding} of a quantum surface is a choice of representative $(D,h)$ from the equivalence class. We sometimes abuse notation and let $(D,h)$ denote a $\gamma$-LQG surface (i.e., an equivalence class) rather than an embedding of this $\gamma$-LQG surface; the meaning will be clear from the context. We often want to decorate a quantum surface by one or more marked points or curves. In this case we define equivalence classes via~\eqref{eq-quantum-surface}, and further require that the conformal map $\varphi$ maps decorations on the first surface to corresponding decorations on the second surface.

We frequently consider non-probability measures in this paper, and extend the usual language of probability theory to this setting. Precisely, consider a triple $(\Omega, \cF, M)$ with $\Omega$ a sample space,  $\cF$ a $\sigma$-algebra on $\Omega$, and $M: \cF \to [0, \infty]$ a measure (not necessarily with $M(\Omega) = 1$). If $X$ is a $\cF$-measurable function (``random variable''), its \emph{law} is the pushforward measure $M_X = X_*M$. We write $X \sim M_X$ and say that $X$ is \emph{sampled} from $M_X$. \emph{Weighting} the law of $X$ by $f(X) \in \R_+$ corresponds to defining the measure $\wt M_X$ via the Radon-Nikodym derivative $\frac{d \wt M_X}{dM_X} = f$. %
For an event $E \in \cF$ with $M[E] \in (0, \infty)$, \emph{conditioning} on $E$ yields the probability measure $\frac{M[\,\cdot\, \cap E]}{M[E]}$ on the measurable space $(E, \cF_E)$ with $\cF_E = \{ S \cap E \: :\: S \in \cF\}$.

For $W \geq \frac{\gamma^2}2$, the \emph{weight $W$ quantum disk} was introduced in \cite[Section 4.5]{wedges} in terms of Bessel processes (see also \cite[Section 3.5]{ghs-mating-survey}). Since this quantum surface has the topology of the disk, we will call it a \emph{thick quantum disk}. 
\begin{definition}[Thick quantum disk]\label{def-thick-disk}
For $W \geq \frac{\gamma^2}2$, we define an infinite measure $\cM_2^\disk(W)$ on two-pointed quantum surfaces $(\cS, \psi, +\infty, -\infty)$ with field $\psi$ as follows. Write $\beta := \frac\gamma2 + Q - \frac W\gamma$. Sample the field $\wh h$ on $\cS$ having independent projections to $\cH_1(\cS)$ and $\cH_2(\cS)$ given by: 
\begin{itemize}
\item 
\[ \wh h_t =
\left\{
	\begin{array}{ll}
		B_{2t} - (Q -\beta)t  & \mbox{if } t \geq 0 \\
		\wt B_{-2t} +(Q-\beta) t & \mbox{if } t < 0
	\end{array}
\right. , \]
where $(B_s)_{s \geq 0}, (\wt B_s)_{s \geq 0}$ are independent standard Brownian motions conditioned on $B_{2s} - (Q-\beta)s<0$ and $\wt B_{2s} - (Q-\beta)s < 0$ for all $s>0$.
\item The projection of an independent Neumann GFF on $\cS$ to $\cH_2(\cS)$. 
\end{itemize}
Independently take the real number $\textbf{c} \sim \frac\gamma2 e^{(\beta-Q)c}dc$, write $\psi := \wh h + \textbf{c}$, and output the quantum surface $(\cS, \psi, +\infty, -\infty)$. 
\end{definition}
Although the Brownian motions are conditioned on a probability zero event, they can be understood by  limiting procedures. Alternatively, with $\delta = 2+\frac2\gamma (Q-\beta)$, the process $(B_{2s} - (Q-\beta)s)_{s \geq 0}$ conditioned on $B_{2s} - (Q-\beta)s <0$ for all $s > 0$  can be sampled by running a dimension $(4-\delta)$ Bessel process  $(Z_s)_{[0,\tau]}$ started from $Z_0 = 0$ until the first time $\tau$ that $Z_\tau = 1$. Then $(B_{2t} - (Q-\beta)t)_{t \geq 0}$ is the time-reversal of $(\frac2\gamma \log Z_t)_{t \in [0,\tau]}$ with time reparametrized in $[0,\infty)$ so the process has quadratic variation $2dt$. 

Definition~\ref{def-thick-disk} is a rephrasing of \cite[Definition 4.21]{wedges} using the Bessel process description \cite[Remark 3.7]{wedges} (see also \cite{py-bessel-decomp}).
The law of $\cc$ corresponds to the  fact that the maximum value of a dimension $\delta = 2+\frac2\gamma (Q-\beta)$ Bessel excursion has the power law $1_{m>0}m^{\delta - 3}dm$. %

In Section~\ref{subsec-thin-disks} we will extend the definition of $\cM_2^\disk(W)$ to $W \in (0, \frac{\gamma^2}2)$, and call these quantum surfaces \emph{thin quantum disks}. The adjective ``thin'' here is inherited from thin quantum wedges defined in~\cite{wedges}.
\subsection{Main results}\label{subsec-main-results}
In this section we state our main results. There are several definitions and details which we only describe at high level; we discuss these more comprehensively in later subsections.

We will want to conformally weld quantum disks according to the natural $\gamma$-LQG boundary measure called \emph{quantum length}: If $h$ on $\cS$ is locally absolutely continuous with respect to a Neumann GFF, then the quantum boundary length measure $\nu_h(dx)$ can be defined as ``$e^{\frac\gamma2 h(x)}dx$'' (this is done rigorously by mollifying and renormalizing \cite{shef-kpz}), and satisfies for continuous $g$ the scaling relation $\nu_{h+g}(dx) = e^{\frac\gamma2 g(x)} \nu_h(dx)$. 

For $\kappa \in (0,4)$ and $\rho_1, \rho_2 > -2$, in a simply connected domain with two marked boundary points $\SLE_\kappa(\rho_1; \rho_2)$ is a conformally invariant chordal random curve between these points, with $\SLE_\kappa = \SLE_\kappa(0;0)$ \cite{lsw-restriction,dubedat-rho}. This is defined in \cite[Section 2.2]{ig1} using Loewner evolutions; the details are not needed for our work so we omit them.
 While ordinary $\SLE_\kappa$ does not hit the domain boundary, when $\rho_1$ (resp. $\rho_2$) is strictly less than $\frac{\kappa}{2}-2$ the $\SLE_\kappa(\rho_1; \rho_2)$ curve a.s.\ hits the left (resp.\ right) boundary arc. 
 $\SLE_\kappa$ curves arise as the welding interface of $\gamma$-LQG surfaces by quantum length when $\kappa = \gamma^2$, and hence one can define the $\gamma$-LQG length of $\SLE_\kappa$-type curves for $\kappa = \gamma^2$ by pushing forward the length measure along the boundary \cite{shef-zipper}. Here and in the rest of the paper, we will take $\kappa = \gamma^2$. 

We will define for $W > 0$ and $\ell, \ell' > 0$ the family of measures $\{\cM_2^\disk(W; \ell, \ell')\}_{\ell, \ell'>0}$ such that $\cM_2^\disk(W; \ell, \ell')$ is supported on quantum surfaces with left and right boundary arcs having quantum lengths $\ell$ and $\ell'$, respectively. This family satisfies
\eqb \label{eq-intro-disint}
\cM_2^\disk(W) = \iint_0^\infty \cM_2^\disk(W; \ell, \ell') \,d \ell\, d \ell'.
\eqe
The relation~\eqref{eq-intro-disint} in fact characterizes $\cM_2^\disk(W, \ell, \ell')$ modulo a Lebesgue measure zero set of values of $(\ell, \ell')$. We will remove this ambiguity by introducing a suitable topology in Sections~\ref{subsec-thick-disint} and~\ref{subsec-thin-decomp} for which $\cM_2^\disk(W, \ell, \ell')$ is continuous in $\ell,\ell'$.
 
When $W \geq \frac{\gamma^2}2$, we let $\cM_2^\disk(W; \ell, \ell') \otimes \SLE_\kappa(\rho_1; \rho_2)$ denote the measure on curve-decorated quantum surfaces obtained by taking a quantum disk $(\cS, \psi, +\infty, -\infty) \sim  \cM_2^\disk(W;\ell,\ell')$ with an arbitrary embedding in $\cS$, independently sampling $\eta \sim \SLE_\kappa(\rho_1;\rho_2)$ in $(\cS, +\infty,-\infty)$, and outputting $(\cS,\psi, +\infty, -\infty, \eta)$. When $W \in (0, \frac{\gamma^2}2)$, the measure $\cM_2^\disk(W; \ell, \ell') \otimes \SLE_\kappa(\rho_1; \rho_2)$ corresponds to sampling independent $\SLE_\kappa(\rho_1;\rho_2)$-curves in each component of the thin quantum disk. We emphasize that for all $W > 0$ our definition of the measure $\cM_2^\disk(W; \ell, \ell') \otimes \SLE_\kappa(\rho_1; \rho_2)$ does not depend on the choice of embedding. Moreover, $\cM_2^\disk(W; \ell, \ell')$ is a measure on $\gamma$-LQG quantum surfaces with dependence on $\gamma$ implicit. 

For fixed $\ell, \ell', \ell_1$, a pair of quantum disks from $\cM_2^\disk(W_1; \ell, \ell_1) \times \cM_2^\disk(W_2; \ell_1, \ell')$ can almost surely be \emph{conformally welded}  along their length $\ell_1$ boundary arcs according to quantum length, to obtain a quantum surface with two marked points joined by an interface locally absolutely continuous with respect to $\SLE_\kappa(W_1 -2; W_2-2)$. 
See e.g.\ \cite{shef-zipper}, \cite[Section 3.5]{wedges} or \cite[Section 4.1]{ghs-mating-survey} for more details on the conformal welding of quantum surfaces.
In the following theorem, we  identify $\cM_2^\disk(W_1; \ell, \ell_1) \times \cM_2^\disk(W_2; \ell_1, \ell')$ with the law of the curve-decorated quantum surface obtained from conformal welding.

\begin{theorem}[Conformal welding of quantum disks]\label{thm-disk-cutting-2}
Suppose $W_1, W_2 > 0$. There exists a constant $c_{W_1,W_2}\in (0,\infty)$ such that for all $\ell, \ell' >0$, the following identity holds as measures on the space of curve-decorated quantum surfaces:
\eqb\label{eq-thm-cut-disk}
\cM_2^\disk(W_1+W_2; \ell, \ell') \otimes \SLE_\kappa(W_1-2; W_2-2)
= c_{W_1,W_2} \int_0^\infty \cM_2^\disk(W_1; \ell, \ell_1)\times \cM_2^\disk(W_2; \ell_1, \ell') \, d \ell_1.
\eqe
\end{theorem}

We now generalize to the multiple curve setting. For $n>1$ and $W_1, \dots, W_n >0$, we define a conformally invariant probability measure $\cP^\disk(W_1, \dots, W_n)$ on $(n-1)$-tuples of $\SLE_\kappa$-type chordal curves in a simply connected domain with two marked points.  This is a special case of \emph{multiple SLE}; see Section~\ref{subsec-SLE} for a precise definition. We note that if $(\eta_1, \dots, \eta_{n-1}) \sim \cP^\disk(W_1, \dots, W_n)$, then for $j = 0, \dots, n-1$, a.s.\ $\eta_j$ and $\eta_{j+1}$ intersect (other than at their endpoints) if and only if $W_j < \frac{\gamma^2}2$; here $\eta_0$ and $\eta_n$ denote the domain boundary arcs. 

We define $\cM_2^\disk(W; \ell, \ell') \otimes \cP^\disk(W_1, \dots, W_n)$ in the same way as $\cM_2^\disk(W; \ell, \ell') \otimes \SLE_\kappa(\rho_1; \rho_2)$, and emphasize that this is a measure on curve-decorated quantum surfaces (with no dependence on embedding). %

\begin{theorem}[Welding multiple disks]\label{thm-disk-cutting}
For $n\in\{2,3,\dots \}$, consider $W_1, \dots, W_n>0$ and $W = \sum_i W_i$. There exists a constant $c_{W_1, \dots, W_n}\in (0,\infty)$ such that for all $\ell, \ell' >0$, the following identity holds as measures on the space of curve-decorated quantum surfaces:
\begin{align}
&\cM_2^\disk(W; \ell, \ell') \otimes \cP^\disk(W_1, \dots, W_n)\nonumber
\\&= c_{W_1,\dots,W_n} \iiint_0^\infty \cM_2^\disk(W_1; \ell, \ell_1)\times \cM_2^\disk(W_2; \ell_1, \ell_2) \times \cdots \times \cM_2^\disk(W_n; \ell_{n-1}, \ell') \, d \ell_1 \dots  d\ell_{n-1}.
\label{eq-disk-cutting-2}
\end{align}
Here the integrand on the right hand side is understood as the law of the quantum surface decorated by $n-1$ curves obtained from conformal welding.
\end{theorem}

Now, we treat the case of quantum spheres. For each $W > 0$ there is a natural infinite measure $\cM_2^\sph(W)$ on sphere-homeomorphic doubly-marked quantum surfaces called \emph{weight $W$ quantum spheres} (see Section~\ref{subsec-wedges}). We can also define a conformally invariant probability measure $\cP^\sph(W_1, \dots, W_n)$ on $n$-tuples of curves in the Riemann sphere between two marked points (see Section~\ref{subsec-SLE}). 
\begin{theorem}[The quantum sphere case]\label{thm-sphere-cutting}
For $n \geq 1$,  consider $W_1, \dots, W_n \in (0, \infty)$ and $W = \sum_i W_i$. There is a constant $\wh c_{W_1, \dots, W_n} \in (0, \infty)$ such that the following identity holds as measures on curve-decorated quantum surfaces:
\begin{align*}
&\cM_2^\mathrm{sph}(W) \otimes \cP^\sph(W_1, \dots, W_n)
\\
&= \wh c_{W_1, \dots, W_n} \iiint_0^\infty \cM_2^\disk(W_1; \ell_0, \ell_1) \times \cM_2^\disk(W_2; \ell_1, \ell_2)\times \cdots \times \cM_2^\disk(W_n; \ell_{n-1}, \ell_0)\, d \ell_0 \dots d\ell_{n-1}.  
\end{align*}
\end{theorem}

Finally, we comment on  several works related to Theorem~\ref{thm-disk-cutting-2}. 

\begin{remark}[{Relation to \cite{cle-lqg}}]
\cite[Lemma 3.3]{cle-lqg} is a version of Theorem~\ref{thm-disk-cutting-2} for $\gamma > \sqrt2$ and weights $W_1 = 2$ and $W_2 \in (0, \gamma^2-2)$. They identify the law of the left quantum surface as $\cM_2^\disk(2)$, and note that conditioned on the left boundary lengths of the components of the right quantum surface, it is a conditionally independent collection of weight $\gamma^2 - W_2$ quantum disks.
\end{remark}

\begin{remark}[Relation to \cite{gm-disk}] \label{rem-GM}
The argument of \cite[Theorem 1.5]{gm-disk} can be adapted to show that the free Boltzmann chordal-self-avoiding-walk-decorated quadrangulation of the disk converges in the Gromov-Hausdorff-Prokhorov-uniform topology to the welding of $\sqrt{\frac83}$-LQG quantum disks of weight 2 along a boundary arc. Theorem~\ref{thm-disk-cutting-2} identifies this limit as a weight 4 quantum disk decorated by an independent $\SLE_{8/3}$ curve. This establishes a scaling limit result of self-avoiding walks to $\SLE_{8/3}$, and can be considered a finite-area analog of \cite[Theorem 1.1]{gwynne-miller-saw}.  
Moreover, based on \cite{gm-disk} and our Theorem~\ref{thm-sphere-cutting}, we prove in~\cite[Theorem 6.10]{AHS-SLE-integrability} that random quadrangulations decorated by a self-avoiding polygon converge to a quantum sphere decorated by an SLE loop.
\end{remark}

\begin{remark}[{Relation to \cite{ct-ising-length}}]\label{rem-interface-length}
Theorem~\ref{thm-disk-cutting-2} implies that when a weight $W_1 + W_2$ quantum disk conditioned to have boundary arc lengths $\ell, \ell'$ is decorated by an independent $\SLE_\kappa(W_1-2; W_2-2)$, then the law of the interface length $L$ is given by 
\eqb \label{eq-len-law}
\P[L \in dx] = Z^{-1} \left|\cM_2^\disk(W_1, \ell ,x)\right| \left| \cM_2^\disk(W_2, x, \ell')\right|,
\eqe
where $Z$ is the normalizing constant. In fact, by scaling and resampling properties of the weight 2 quantum disk we have $|\cM_2^\disk(2, \ell, x)| = C(\ell + x)^{-\frac{4}{\gamma^2} - 1}$ for all $\ell, x > 0$ for some constant $C>0$, and consequently, for $\gamma =\sqrt3$ and $W_1 = W_2 = 2$ the law~\eqref{eq-len-law} agrees with the scaling limit result of \cite[Equation (46)]{ct-ising-length} for the critical Boltzmann triangulation of the disk decorated by an Ising model with Dobrushin boundary conditions. The computation identifying~\eqref{eq-len-law} and this limit law is carried out (in different language) immediately after \cite[Equation (46)]{ct-ising-length}. 
\end{remark}

\begin{remark}
	We have only stated Theorems \ref{thm-disk-cutting} and \ref{thm-sphere-cutting} for the subcritical case $\kappa=\gamma^2\in(0,4)$ but we expect that they also hold in the critical case $\kappa=\gamma^2=4$ via similar techniques as in Appendix \ref{appendix-A} and \cite{hp-welding} (where the latter work builds on  \cite{aps-critical-lqg-lim,mmq-welding}). To carry out such an argument, we would take the subcritical theorems as input and  take the limit $\gamma\uparrow 2$ similarly as in Appendix \ref{appendix-A}, but using additional ingredients from \cite{hp-welding,aps-critical-lqg-lim,mmq-welding}.
	\label{rmk:critical} 
\end{remark}

\subsection{Quantum wedges, spheres and cones}\label{subsec-wedges}
In this section, we recall the definitions of various quantum surfaces from \cite[Section 4.2--4.5]{wedges}. See also \cite[Sections 3.4--3.5]{ghs-mating-survey}. We omit the weight $\frac{\gamma^2}2$ quantum wedge description as it is not needed.

\begin{definition}[Thick quantum wedge] \label{def-thick-wedge}
For $W > \frac{\gamma^2}2$, let $\beta := \frac\gamma2 + Q - \frac W \gamma$ (so $\beta < Q$). Then $\cM^\wed(W)$ is the probability measure on doubly-marked quantum surfaces $(\cS, \wt h, +\infty, -\infty)$, where the field $\wt h$ has independent projections to $\cH_1(\cS)$ and $\cH_2(\cS)$ given by:
\begin{itemize}
\item 
\[ \wt h_t :=
\left\{
	\begin{array}{ll}
		B_{2t} - (Q- \beta) t  & \mbox{if } t \geq 0 \\
		\wt B_{-2t} - (Q-\beta)t & \mbox{if } t < 0
	\end{array}
\right. , \]
where $(B_s)_{s \geq 0}, (\wt B_s)_{s \geq 0}$ are independent standard Brownian motions conditioned on $\wt B_{2s} + (Q-\beta)s > 0$ for all $s>0$.

\item The projection of an independent Neumann GFF on $\cS$ to $\cH_2(\cS)$. 
\end{itemize}
\end{definition}

The thin quantum wedge arises as a concatenation of thick quantum disks. 
\begin{definition}[Thin quantum wedge]\label{def-thin-wedge}
Fix $W \in (0, \frac{\gamma^2}2)$ and sample a Poisson point process $\{ (u, \cD_u)\}$ from the measure $\Leb_{\R_+} \otimes \cM_{2}^\disk(\gamma^2 - W)$. The weight $W$ quantum wedge is the infinite beaded surface obtained by concatenating the $\cD_u$ according to the ordering induced by $u$. We write $\cM^\wed(W)$ for the probability measure on weight $W$ quantum wedges. 
\end{definition}
This agrees with the definition in \cite[Section 4.4]{wedges}. Indeed, let $M$ be the excursion measure of a Bessel process with dimension less than $2$, then a Bessel process can be obtained by sampling a Poisson point process from $\mathrm{Leb}_{\R_+}\times M$ and concatenating the excursions according to the ordering induced by the first coordinate.

Now we define analogs of quantum disks and wedges with no boundaries. For notational simplicity we work in the cylinder $\cC := \R \times [0,2\pi] /\sim$, where we identify $\R \times \{0\}$ and $\R \times \{2\pi\}$ by the equivalence $x \sim x + 2\pi i$. We can define $H(\cC), \cH_1(\cC)$, and $\cH_2(\cC)$ as in the strip setting, and thus make sense of the Neumann GFF in $\cC$ with mean zero on $[0,2\pi i]$.

\begin{definition}[Quantum sphere]\label{def-sphere}
For $W >0$, we define a infinite measure $\cM_2^\sph(W)$ on two-pointed quantum surfaces $(\cC, \psi, +\infty, -\infty)$ with field $\psi$ defined as follows. Write $\alpha := Q- \frac W{2\gamma}$. Sample the field $\wh h$ on $\cC$ having independent projections to $\cH_1(\cC)$ and $\cH_2(\cC)$ given by: 
\begin{itemize}
\item 
\[ \wh h_t =
\left\{
	\begin{array}{ll}
		B_t - (Q -\alpha)t  & \mbox{if } t \geq 0 \\
		\wt B_{-t} +(Q-\alpha) t & \mbox{if } t < 0
	\end{array}
\right. , \]
where $(B_s)_{s \geq 0}, (\wt B_s)_{s \geq 0}$ are independent standard Brownian motions conditioned on $B_s - (Q-\beta)s<0$ and $\wt B_s - (Q-\beta)s < 0$ for all $s>0$.

\item The projection of an independent Neumann GFF on $\cC$ to $\cH_2(\cC)$. 
\end{itemize}
Independently take the real number $\mathbf{c} \sim \frac\gamma2 e^{2(\alpha-Q)c}\,dc$, write $\psi := \wh h + \mathbf{c}$, and output the quantum surface $(\cC, \psi, +\infty, -\infty)$. 
\end{definition}

\begin{definition}[Quantum cone] \label{def-cone}
For $W > 0$, let $\alpha := Q - \frac W {2\gamma}$ (so $\alpha < Q$). Then $\cM^\cone(W)$ is the probability measure on doubly-marked quantum surfaces $(\cC, \wt h, +\infty, -\infty)$, where the field $\wt h$ has independent projections to $\cH_1(\cC)$ and $\cH_2(\cC)$ given by:
\begin{itemize}
\item 
\[ \wt h_t :=
\left\{
	\begin{array}{ll}
		B_{t} - (Q- \alpha) t  & \mbox{if } t \geq 0 \\
		\wt B_{-t} - (Q-\alpha)t & \mbox{if } t < 0
	\end{array}
\right. , \]
where $(B_s)_{s \geq 0}, (\wt B_s)_{s \geq 0}$ are independent standard Brownian motions conditioned on $\wt B_{s} + (Q-\alpha)s > 0$ for all $s>0$.

\item The projection of an independent Neumann GFF on $\cC$ to $\cH_2(\cC)$. 
\end{itemize}
\end{definition}

\subsection{Thin quantum disks}\label{subsec-thin-disks}
In this section we define the thin quantum disk with weight $W \in (0, \frac{\gamma^2}2)$. There is a thin-thick duality $W \leftrightarrow \gamma^2-W$ in the sense that thin quantum disks are a Poissonian chain of thick quantum disks from $\cM_2^\disk(\gamma^2 - W)$. At first glance the nontrivial topology of thin quantum disks seems unnatural, but this topology enables our arguments in this paper and subsequent work, and we will see that the thin quantum disks are the natural analogue of thick quantum disks for $W\in (0, \frac{\gamma^2}2)$.

In Definition~\ref{def-thin-wedge}, although the thin quantum wedge $\cW$ only comes with the ordering on thick quantum disks and not the labels $u$, we will see in Corollary~\ref{cor-cut-measure-intrinsic} that $\{(u, \cD_u)\}$ is measurable with respect to $\cW$. Therefore it makes sense to define the \emph{quantum cut point measure} which assigns mass $x$ to the collection of cut points between the quantum disks $\{ \cD_u \: : \: u \leq x\}$. 

For each $\alpha \in (0,1)$, a L\'evy process $(L_t)_{t \geq 0}$ is called an \emph{$\alpha$-stable subordinator} if it is a.s.\ increasing and $Y_{a t} \stackrel d= a^{1/\alpha} Y_t$ \cite{bertoin-book}. The following result is proved (in different language) in \cite[Lemma 2.6]{ghm-kpz}. See also Remark~\ref{rem-stable-subordinators} for a succinct self-contained proof. 
\begin{lemma}\label{lem-stable-subordinators}
Consider a weight $W \in (0, \frac{\gamma^2}2)$ thin quantum wedge. 
Let $L_t$ denote the total quantum length of the left boundary arcs of the thick quantum disks within quantum cut point distance $t$ of the root of the quantum wedge. Then $(L_t)_{t\geq0}$ is a stable subordinator with exponent $1 - \frac{2W}{\gamma^2}$. The same holds if $L_t$ is instead defined as the sum of the right boundary arc lengths, or the sum of the perimeters.
\end{lemma}

\begin{corollary}[Intrinsic definition of cut point measure]\label{cor-cut-measure-intrinsic}
Parametrize the left (or right) boundary of a thin quantum wedge by quantum length, and let $\mathcal T \subset \R_+$ be the set corresponding to cut points along this boundary. Then the quantum cut point measure is given by the $(1 - \frac{2W}{\gamma^2})$-Minkowski content measure of $\mathcal T$ multiplied by a deterministic constant.
\end{corollary}
\begin{proof}
$\mathcal T$ is the range of the stable subordinator $L_t$ with exponent $(1-\frac{2W}{\gamma^2})$, so the $(1-\frac{2W}{\gamma^2})$-Minkowski content on $\mathcal T$ agrees (up to deterministic constant) with the pushforward of Lebesgue measure on $\R$ under the map $t \mapsto L_t$ \cite[Lemma 5.13]{hs-cardy-embedding}.
\end{proof}

We now introduce the infinite measure on thin quantum disks, and note that each thin quantum disk is a concatenation of quantum surfaces with finite total area and boundary length.

\begin{definition}[Thin quantum disk]\label{def-thin-disk}
For $W \in (0, \frac{\gamma^2}2)$, we can define the infinite measure $\cM_2^\disk(W)$ on two-pointed beaded surfaces as follows. %
Take $T \sim  (1-\frac2{\gamma^2} W)^{-2}\mathrm{Leb}_{\R_+}$, then sample a Poisson point process $\{(u, \cD_u)\}$ from the measure $\Leb_{[0,T]} \times \cM_2^\disk(\gamma^2 - W)$, and concatenate the $\cD_u$ according to the ordering induced by $u$.  
\end{definition}
The choice of constant $(1-\frac2{\gamma^2} W)^{-2}$ will be justified in \cite[Section 3.2]{AHS-SLE-integrability}, where, roughly speaking, we understand thin quantum disks as an ``analytic continuation'' of thick quantum disks. The quantum cut point measure on a thin quantum disk is well defined and measurable with respect to the thin quantum disk, see Corollary~\ref{cor-cut-measure-intrinsic}.

\subsection{Basic properties of quantum disks}\label{subsec-basic}
In this section, we discuss some basic properties of quantum disks. For a domain $(D,x,y)$ we define its left boundary arc to be the clockwise arc from $x$ to $y$. 

\begin{lemma}[Thick quantum disk boundary length law]\label{lem-length-exponent}
For a quantum disk of weight $W\in [\frac{\gamma^2}2, \gamma Q)$, the $\cM_2^\disk(W)$-law of the quantum length of its left boundary arc is $c_W \ell^{-\frac2{\gamma^2}W} d \ell$ for some $c_W \in (0,\infty)$. If $W \geq \gamma Q$, then the $\cM_2^\disk(W)$-measure of $\{\text{left boundary length} \in I\}$ is infinite for any open interval $I \subset\R_+$. 

The same is true for the quantum length of the right boundary arc or whole boundary. 
\end{lemma}
\begin{proof}
Write $\beta = \frac\gamma2 + Q - \frac W\gamma$.
For $0< L < L'$ we have, with $(\wh h, c)$ as in Definition~\ref{def-thick-disk},
\alb
&\cM_2^\disk [ \nu_{{\wh h}+c}(\R) \in (L, L')] = \E\left[ \int_{-\infty}^\infty \1_{e^{\frac\gamma2c} \nu_{\wh h}(\R) \in (L, L')} \frac\gamma2e^{(\beta - Q)c}\,dc \right]\\ 
&\qquad\qquad= \E \left[\int_L^{L'} \nu_{\wh h}(\R)^{\frac2\gamma(Q-\beta)} y^{\frac2\gamma (\beta - Q)} \cdot  y^{-1} \, dy\right]
= \E \left[ \nu_{\wh h}(\R)^{\frac2\gamma(Q-\beta)}  \right]\int_L^{L'}  y^{\frac2\gamma (\beta - Q)-1}  \, dy,
\ale
where we have used the change of variables $y = e^{\frac\gamma2c}\nu_{\wh h}(\R)$ (so $dc = \frac2\gamma y^{-1} dy$). When $W \in [\frac{\gamma^2}2, \gamma Q)$ (so $\frac2\gamma (Q-\beta) < \frac4{\gamma^2}$), the expectation is finite. This can be proved by H\"older's inequality as in the proof of \cite[Lemma A.4]{wedges}, or alternatively see \cite[Theorem 1.7 and (1.32)]{rz-boundary}. Since $\frac2\gamma(\beta - Q) - 1 = -\frac2{\gamma^2}W$, we conclude that the law of the left boundary arc quantum length is $c \ell^{-\frac2{\gamma^2}W} d\ell$. If $W \geq \gamma Q$ (so $\frac2\gamma (Q-\beta) \geq \frac4{\gamma^2}$)   then the expectation is infinite \cite[Proposition 3.5]{robert-vargas-revisited}, as desired. 
\end{proof}

\begin{remark}[Proof of Lemma~\ref{lem-stable-subordinators}]\label{rem-stable-subordinators}
By Lemma~\ref{lem-length-exponent} the left boundary length of a quantum disk from $\cM_2^\disk(\gamma^2-W)$ has law $\Pi(dL) = c L^{-\frac2{\gamma^2}(\gamma^2-W)} dL$ for some $c \in (0,\infty)$. Therefore, with $\Lambda$ a Poisson point process on $\Leb_{\R_+} \times \Pi(dL)$, $L_t$ is given by the sum of $L$ over all $(u, L) \in \Lambda$ with $u \leq t$, for $t \geq 0$. Since $\Pi$ is the L\'evy measure of a stable subordinator with exponent $\frac{2}{\gamma^2} (\gamma^2 - W) - 1 = 1-\frac{2W}{\gamma^2}$, by the L\'evy-It\^o decomposition we conclude that $(L_t)_{t \geq 0}$ is a stable subordinator with exponent $1-\frac{2W}{\gamma^2}$.  

\end{remark}

\begin{lemma}[Thin quantum disk boundary length law]\label{lem-thin-disk-perim}
For $W \in (0, \frac{\gamma^2}2)$, the left boundary length of a thin quantum disk from $\cM_2^\disk(W)$ is distributed as $c_W \ell^{-\frac2{\gamma^2}W}d\ell$ for some constant $c_W \in (0,\infty)$. The same is true for the right boundary length or the total boundary length.
\end{lemma}
\begin{proof}
Let $(L_t)_{t\geq 0}$ be the stable subordinator of exponent $\alpha := 1 - \frac{2W}{\gamma^2} \in (0,1)$ described in Lemma~\ref{lem-stable-subordinators} for the left boundary length. Writing $C =  (1-\frac2{\gamma^2} W)^{-2}$, the $\cM_2^\disk(W)$-measure of the event that the left boundary arc length lies in $[\ell, \ell + d\ell]$ is given by
\[ \E\left[ \int_0^\infty 1_{L_T \in [\ell, \ell + d\ell]} C\, dT\right] = C \E\left[\int_0^\infty 1_{T^{1/\alpha} L_1 \in [\ell, \ell + d\ell]} \, dT \right] = C \E\left[\int_{(\ell/L_1)^\alpha}^{((\ell + d\ell)/L_1)^\alpha} 1\, dT \right] = C \alpha \ell^{\alpha - 1} \E[L_1^{-\alpha}].\]
We are done once we check that $\E[L_1^{-\alpha}] < \infty$. Since $(L_t)_{t\geq 0}$ is a stable subordinator with index $\alpha$, we know $\E[e^{-\lambda L_1}] = e^{-c\lambda^\alpha}$ for some $c>0$ for all $\lambda > 0$; indeed we have $\lambda L_1 \stackrel d= L_{\lambda^\alpha}$ 
and, since $(L_t)_{t\geq 0}$ has nonnegative and stationary increments, $\E[e^{-L_{t}}] = e^{-ct}$.
Hence
\[\Gamma(\alpha) \E[L_1^{-\alpha}] = \E \left[ \int_0^\infty e^{-\lambda L_1} \lambda^{\alpha - 1} \,d\lambda \right]= \int_0^\infty e^{-c \lambda^\alpha} \lambda^{\alpha - 1} d\lambda <\infty. \qedhere \]
\end{proof}

The exponent $-\frac2{\gamma^2}W$ in Lemmas~\ref{lem-length-exponent} and~\ref{lem-thin-disk-perim} is natural in light of the following lemma, which explains how the quantum disk measure scales after adding a constant to the quantum disk field. We note that Lemmas~\ref{lem-length-exponent} and~\ref{lem-thin-disk-perim} are not immediate consequences of Lemma~\ref{lem-scaling} as it does not yield finiteness/infiniteness of the constant $c_W$.
\begin{lemma}[Scaling property of quantum disks]\label{lem-scaling}
For $W > 0$, the following procedures agree for all $\lambda > 0$:
\begin{itemize}
\item Take a quantum disk from $\cM_2^\disk(W)$.
\item Take a quantum disk from $\lambda^{-\frac2{\gamma^2}W +1} \cM_2^\disk(W)$ then add $\frac2\gamma \log \lambda$ to its field. 
\end{itemize}
\end{lemma}
\begin{proof}
\noindent When $W \geq \frac{\gamma^2}2$, this is immediate from Definition~\ref{def-thick-disk} because of the constant term $\mathbf{c} \sim e^{ (-\frac2{\gamma^2}W-1)\frac\gamma2c}\,dc$ (written here in terms of $W$ rather than $\beta$). 

When $W \in (0,\frac{\gamma^2}2)$, by the previous case the following procedures yield the same law on beaded quantum surfaces for fixed $T>0$:
\begin{itemize}
\item Sample a Poisson point process from $\Leb_{[0,T]} \times \cM_2^\disk(\gamma^2 - W)$.
\item Sample a Poisson point process from $\Leb_{[0,\lambda^{-\frac2{\gamma^2}(\gamma^2 - W)+1}T]} \times \cM_2^\disk(\gamma^2 - W)$ then add $\frac2\gamma \log \lambda$ to the field.
\end{itemize}
In Definition~\ref{def-thin-disk} we take $T$ from a multiple of Lebesgue measure, so the scaling by $\lambda^{-\frac2{\gamma^2}(\gamma^2 - W) + 1} = \lambda^{-\frac2{\gamma^2}W + 1}$ yields the claim for thin quantum disks. 
\end{proof}

\subsection{Definition of $\cM_2^\disk(W; \ell_1, \ell_2)$}\label{subsec-thick-disint}
In this section, we explain  the construction of the disintegration from~\eqref{eq-intro-disint}, which is
\[
\cM_2^\disk(W) = \iint_0^\infty \cM_2^\disk(W; \ell_1, \ell_2) \,d \ell_1 \, d \ell_2
\]
where $\cM_2^\disk (W; \ell_1, \ell_2)$ is supported on the set of doubly-marked quantum surfaces with left and right boundary arcs having quantum lengths $\ell_1$ and $\ell_2$, respectively.

We first discuss in more generality disintegrations of measures. Suppose $M$ is a $\sigma$-finite measure on a Radon space $(X, \cX)$ and $T: X \to \R^n$ is a measurable function such that the pushforward $T_*M$ is absolutely continuous with respect to  $\Leb_{\R^n}$ (Lebesgue measure on $\R^n$). Then there exists a collection of $\sigma$-finite measures $\{ M_{t_1, \dots, t_n}\}_{(t_1, \dots, t_n)\in\R^n}$ such that:
\begin{itemize}
\item $M_{t_1, \dots, t_n}$ is supported on $T^{-1}(t_1, \dots, t_n)$ for all $(t_1, \dots, t_n) \in \R^n$;
\item For any $A \in \cX$ the function $(t_1, \dots, t_n) \mapsto M_{t_1, \dots, t_n}[A]$ is measurable, and 
\[M[A] = \iiint M_{t_1, \dots, t_n}[A] \,dt_1 \dots dt_n.\]
\end{itemize}
We call this collection $\{M_{t_1, \dots, t_n}\}$ a \emph{disintegration}, and write $M = \iiint M_{t_1, \dots, t_n} \,dt_1\dots dt_n$. Disintegrations are unique in the sense that for any two disintegrations $\{M_{t_1, \dots, t_n}\}$, $\{\wt M_{t_1, \dots, t_n}\}$ we have $M_{t_1, \dots, t_n} = \wt M_{t_1, \dots, t_n}$ for $\Leb_{\R^n}$-a.e. $(t_1, \dots, t_n)$. 

We briefly justify the above claims on disintegrations. When $M$ is a probability measure we can take $M_{t_1, \dots, t_n}$ to be the regular conditional probability distribution multiplied by $f(t_1, \dots, t_n)$, where $T_*M = f(t_1, \dots, t_n) \Leb_{\R^n}(dt_1 \dots dt_n)$. Uniqueness follows from that of regular conditional probability distributions \cite[Chapter 6]{kallenberg}. The extension to $\sigma$-finite $M$ follows by exhaustion.

If one can specify a choice of disintegration $\{ M_{t_1, \dots, t_n}\}$ and a sufficiently strong topology for which the map $(t_1, \dots, t_n) \mapsto M_{t_1, \dots, t_n}$ is continuous, then the disintegration is canonically defined for \emph{all} (and not just a.e.) $(t_1, \dots, t_n)$. We will do this in detail for the disintegration
~\eqref{eq-intro-disint} for the case $W > \frac{\gamma^2}2$ and sketch the necessary modifications for $W = \frac{\gamma^2}2$.

Let $W > \frac{\gamma^2}2$ and define the event
\eqb\label{eq-E'-disk}
\wh E_\zeta := \{\sup_{t \in \R} \psi_t > -\zeta  \}.
\eqe
We now provide an alternative description of $\cM_2^\disk(W)$ restricted to the event $\wh E_\zeta$. 
\begin{lemma}\label{lem-disk-beta}
For $W > \frac{\gamma^2}2$ and $\beta = \frac\gamma2 + Q - \frac W\gamma$, with $\wh E_\zeta$ defined in~\eqref{eq-E'-disk}, we have $\cM_2^\disk (W)|_{\wh E_\zeta} = (\frac{2W}{\gamma^2}-1)^{-1} e^{(Q-\beta)\zeta} \cdot P_\zeta$. Here, $P_\zeta$ is a probability measure on quantum surfaces $(\cS, \psi, +\infty, -\infty)$ where the field $\psi$ has independent projections to $\cH_1(\cS)$ and $\cH_2(\cS)$ given by:
\begin{itemize}
\item 
\[ \psi_t :=
\left\{
	\begin{array}{ll}
		-\zeta + B_{2t} - (Q -\beta)t  & \mbox{if } t \geq 0 \\
		-\zeta + \wt B_{-2t} +(Q-\beta) t & \mbox{if } t < 0
	\end{array}
\right. , \]
where $(B_s)_{s \geq 0}, (\wt B_s)_{s \geq 0}$ are independent standard Brownian motions  conditioned on $\wt B_{2s} - (Q-\zeta)s < 0$ for all $s>0$.

\item The projection of an independent Neumann GFF on $\cS$ to $\cH_2(\cS)$. 
\end{itemize}
\end{lemma}
\begin{proof}
\cite[Proposition 2.14]{ag-disk} explains that when we condition on $\wh E_\zeta$ we get the probability measure $P_\zeta$. To conclude we observe that $\cM_2^\disk(W) [ \wh E_\zeta] = \int_{-\zeta}^\infty \frac\gamma2 e^{(\beta - Q)c}dc = \frac\gamma2 (Q-\beta)^{-1} e^{(Q-\beta)\zeta} = (\frac{2W}{\gamma^2}-1)^{-1} e^{(Q-\beta)\zeta}$. 
\end{proof}
The following lemma provides a  decomposition of $\psi$ into independent components. 
\begin{lemma}\label{lem-decomposition-psi}
Fix $\zeta \in \R$ and a pair of arbitrary functions $f_1,f_2 \in \cH_2(\cS)$ such that $f_1$ (resp. $f_2$) is supported on $[0,1] \times [0,\frac\pi2]$ (resp. $[0,1] \times [\frac\pi2, \pi]$), is positive on the interval $(0,1)$ (resp. $(i\pi, i\pi+1)$), and has Dirichlet energy 1. For $\psi\sim P_\zeta$ sampled as in Lemma~\ref{lem-disk-beta}, we have the following decomposition of $\psi$ into four independent components:
\eqb\label{eq-psi-decomp}
\psi = \psi_+ + \alpha_1 f_1 + \alpha_2 f_2 + \psi_-.
\eqe
Here $\alpha_1,\alpha_2\sim N(0,1)$, $\psi_+ + \alpha_1 f_1+\alpha_2f_2$ is a distribution which is harmonic in $\cS_-$, and $\psi_-$ is a distribution supported in $\cS_-$. The process $((\psi_-)_{-t})_{t \geq 0} \stackrel d= (B_{2t} - (Q-\beta)t)_{t\geq 0}$ where $(B_s)_{s \geq 0}$ is standard Brownian motion conditioned on $B_{2t} - (Q-\beta)t \leq 0$ for all $t > 0$, and independently the projection of $\psi_-$ to $\cH_2(\cS)$ agrees in distribution with the projection to $\cH_2(\cS)$ of a GFF on $\cS_-$ with Neumann boundary conditions on $\partial \cS_- \backslash [0,i\pi]$ and zero boundary conditions on $[0,i\pi]$. 
\end{lemma}
\begin{proof}
Let the projection of $\psi_+$ to $\cH_1(\cS)$ be given by $(\psi_+)_t \equiv -\zeta$ for $t \leq 0$ and $(\psi_+)_t = \psi_t$ for $t \geq 0$. Let the projection of $\psi_-$ to $\cH_1(\cS)$ be given by $(\psi_-)_t = \psi_t + \zeta$ for $t \leq 0$ and $(\psi_-)_t \equiv 0$ for $t\geq 0$. 

Now we describe the (independent) projections of $\psi_+$ and $\psi_-$ to $\cH_2(\cS)$.
Let $H_\mathrm{harm} \subset \cH_2(\cS)$ (resp. $H_\mathrm{supp}\subset \cH_2(\cS)$) be the subspace of functions harmonic (resp. supported) in $\cS_-$. Then $H_\mathrm{harm}$ and $H_\mathrm{supp}$ are orthogonal complements in $\cH_2(\cS)$. We may extend $\{f_1, f_2\}$ to an orthonormal basis $\{f_i\}_\N$ of $H_\mathrm{harm}$ and let $\{g_j\}_\N$ be an orthonormal basis of $H_\mathrm{supp}$, then the projection of a Neumann GFF on $\cS$ to $\cH_2(\cS)$ can be written as $\sum \alpha_i f_i + \sum \beta_j g_j$ where $\alpha_i, \beta_j \sim N(0,1)$ are independent. Writing the projections of $\psi_+$ and $\psi_-$ to $\cH_2(\cS)$ as $\sum_{i=3}^\infty \alpha_i f_i$ and $\sum_{j=1}^\infty \beta_j g_j$, respectively, gives the desired decomposition. See \cite[Section 3.2.4]{ghs-mating-survey} for further discussion.
\end{proof}

We note that conditioned on $\psi_+$, the conditional law of $(\nu_\psi([-1,0]), \nu_\psi([-1+i\pi,i\pi]))$ has a density $g_{\psi_+}(x,y)\, dx\, dy$, where $g_{\psi_+}$ is a nonnegative bounded continuous function on $\R^2$; indeed, $\nu_\psi([0,1]) = \int_0^1 e^{\frac\gamma2 \alpha_1 f_1(x)} \nu_{\psi_+}(dx)$ and $\alpha_1 \sim N(0,1)$, and likewise for $\nu_\psi([i\pi,i\pi+1])$.

\begin{definition}[Disintegration of thick quantum disk measure]
For $W > \frac{\gamma^2}2$ and $\ell_1,\ell_2>0$, define $\cM_2^\disk(W; \ell_1, \ell_2) = \lim_{\zeta \to \infty} \cM_2^\disk(W;\ell_1,\ell_2,\zeta)$, where $\cM_2^\disk(W; \ell_1, \ell_2, \zeta)$ is the measure on quantum surfaces $(\cS,\psi,+\infty,-\infty)$ defined as follows: 

Take $(\psi_+, \psi_-)$ from $(\frac{2W}{\gamma^2}-1)^{-1} e^{(Q-\beta)\zeta} \cdot P_\zeta$ (see Lemmas~\ref{lem-disk-beta} and \ref{lem-decomposition-psi}) and restrict to the event that $d_1 := \ell_1 - \nu_{\psi_+ + \psi_-}(\R \backslash [0,1]) >0$ and $d_2 := \ell_2 - \nu_{\psi_+ + \psi_-}((\R + i\pi) \backslash [i\pi,i\pi+1]) >0$. Weight the measure by $g_{\psi_+}(d_1, d_2)$, and let $\alpha_1, \alpha_2 \in \R$ be the values such that $\psi:= \psi_+ + \psi_- + \alpha_1f_1 + \alpha_2f_2$ satisfies $\nu_\psi(\R) = \ell_1$ and $\nu_\psi(\R + i\pi) = \ell_2$. Output $(\cS,\psi, +\infty, -\infty)$. 
\end{definition}
In the above definition, the $\zeta\to\infty$ limit makes sense because it is straightforward to check that
\[\cM_2^\disk(W;\ell_1,\ell_2,\zeta) = \cM_2^\disk(W;\ell_1,\ell_2,\zeta')|_{\wh E_\zeta} \quad \text{ for }\zeta' > \zeta.\] Therefore we have $\cM_2^\disk(W; \ell_1, \ell_2)|_{\wh E_\zeta} = \cM_2^\disk(W; \ell_1, \ell_2,\zeta)$.

It is clear that $\cM_2^\disk(W)|_{\wh E_\zeta} = \iint_0^\infty \cM_2^\disk(W; \ell_1, \ell_2, \zeta) \,d \ell_1\, d \ell_2$. Sending $\zeta\to \infty$, we obtain~\eqref{eq-intro-disint}, as desired. 
 We see in the next proposition that this disintegration is canonical in the sense that the measures $\cM_2^\disk(W; \ell_1, \ell_2)$ are continuous in $(\ell_1, \ell_2)$ in a suitable topology. In terms of notation, we will write $M^\#:=M/|M|$ to denote the normalized probability measure of a measure $M$.
\begin{proposition}\label{prop-thick-cts}
For $W > \frac{\gamma^2}2$, the family of measures $\{\cM_2^\disk(W;\ell_1, \ell_2)\}_{\ell_1, \ell_2}$ is continuous in $(\ell_1, \ell_2)$ with respect to the metric
\[d(M, \wt M) = \int_0^\infty e^{-\zeta} d_\zeta\left((M|_{\wh E_\zeta})^\#, (\wt M|_{\wh E_\zeta})^\# \right) \,d \zeta,\]
where $d_\zeta$ is the total variation distance between the laws of $\psi(\cdot - \tau_{-\zeta})|_{\cS_+}$ for $\psi \sim (M|_{\wh E_\zeta})^\#$ and $\psi \sim (\wt M|_{\wh E_\zeta})^\#$; here, $\tau_{-\zeta} := \inf  \{t \in \R \: : \: \psi_t = -\zeta\}$.
\end{proposition}
\begin{proof}
It suffices to show that $\{\cM_2^\disk(W;\ell_1, \ell_2)\}_{\ell_1, \ell_2}$ is continuous in $(\ell_1, \ell_2)$ with respect to each $d_\zeta$. This follows from the continuity of $g_{\psi_+}$ for each $\psi_+$. %
\end{proof}
 
Now we sketch the construction and continuity of $\{ \cM_2^\disk(\frac{\gamma^2}2; \ell_1, \ell_2)\}_{\ell_1, \ell_2}$. The previous construction is not applicable for $W = \frac{\gamma^2}2$ for two reasons: Firstly, $\{ \sup_{t \in \R} \psi_t > -\zeta\}$ has infinite $\cM_2^\disk(\frac{\gamma^2}2; \ell_1, \ell_2)$-mass, so we instead use $\wh E_N = \{ \sup_{t \in \R} \psi_t \in [-N,N]\}$ to define $\cM_2^\disk(\frac{\gamma^2}2;\ell_1, \ell_2, N)$. Secondly, the  description Lemma~\ref{lem-disk-beta} does not apply to $W = \frac{\gamma^2}2$, so we use the description of $\psi$ in Definition~\ref{def-thick-disk} (i.e. with embedding so $\psi_0 \geq \psi_t$ for all $t$) to establish a field decomposition like~\eqref{eq-psi-decomp}. Then proceeding as before, we can construct $\cM_2^\disk(\frac{\gamma^2}2; \ell_1, \ell_2) := \lim_{N\to \infty} \cM_2^\disk(\frac{\gamma^2}2; \ell_1, \ell_2, N)$. This family $\{\cM_2^\disk(\frac{\gamma^2}2; \ell_1, \ell_2)\}_{\ell_1, \ell_2}$ is continuous in $(\ell_1, \ell_2)$ with respect to the metric 
\[d'(M, \wt M) = \int_0^\infty e^{-N} d'_N\left((M|_{\wh E_N})^\#, (\wt M|_{\wh E_N})^\#\right) \, d N,\]
where $d'_N$ is the total variation distance between the laws of $\psi|_{\cS_+ - N}$ for $\psi \sim (M|_{\wh E_N})^\#$ and $\psi \sim (\wt M|_{\wh E_N})^\#$, where the fields $\psi$ are chosen with the embedding that $\psi_0 \geq \psi_t$ for all $t$. We note that this approach also works for $W > \frac{\gamma^2}2$, but the previous writeup is more convenient for our later proof. 

For $W \in (0, \frac{\gamma^2}2)$, we can likewise construct a family $\{\cM_2^\disk(W; \ell_1, \ell_2)\}_{\ell_1, \ell_2}$ satisfying
\[\cM_2^\disk(W) = \iint_0^\infty \cM_2^\disk(W; \ell_1, \ell_2)\, d \ell_1 \, d \ell_2 \]
via the earlier discussion on disintegrations. A priori, this family is only unique for a.e.\ $\ell, \ell'$, but we will extend this to a pointwise definition such that $(\ell_1, \ell_2) \mapsto \cM_2^\disk(W; \ell_1, \ell_2)$ is continuous in a similar topology as in the thick case. See Section~\ref{subsec-thin-decomp}.

We now explain how the measure $\cM_2^\disk(W; \ell, \ell')$ scales when adding a constant to the field.
\begin{lemma}\label{lem-scaling-disint}
For $W , \ell, \ell'>0$, the following procedures agree for all $\lambda > 0$:
\begin{itemize}
\item Take a quantum disk from $\cM_2^\disk(W; \lambda\ell, \lambda\ell')$;

\item Take a quantum disk from $\lambda^{-\frac2{\gamma^2}W -1} \cM_2^\disk(W; \ell,\ell')$ then add $\frac2\gamma \log \lambda$ to its field. 
\end{itemize}
\end{lemma}
\begin{proof}
Consider the case $W \geq \frac{\gamma^2}2$ (the other case is similar). Lemma~\ref{lem-scaling} tells us that the measure \\ $\eps^{-2}\iint_1^{1+\eps}\cM_2^\disk(W; a \lambda \ell, b \lambda \ell') \, da\, db$ agrees with $\lambda^{\frac2\gamma (\beta - Q)}\eps^{-2}\iint_1^{1+\eps} \cM_2^\disk(W; a\ell, b\ell') \, da \, db$ (when we add $\frac2\gamma \log \lambda$ to the field).

Send $\eps \to 0$ and note that the first measure converges to $\lambda^2\cM_2^\disk(W; \lambda \ell, \lambda \ell')$, while the second converges to $\lambda^{-\frac2{\gamma^2}W + 1} \cM_2^\disk(W; \ell, \ell')$ (when we add $\frac2\gamma \log \lambda$ to the resulting field), as desired. 
\end{proof}

\subsection{Schramm-Loewner evolution}\label{subsec-SLE}
Now that we have provided details on the quantum surfaces involved in our main theorems, we turn to the relevant SLE curves in these theorems.

$\SLE_\kappa$ is a one-parameter family of conformally invariant random curves introduced by \cite{schramm0}, %
 which arises in the scaling limit of many statistical physics models. It is conformally invariant in the sense that for any pair $(D,x,y), (\wt D, \wt x, \wt y)$ of simply connected domains with two marked boundary points and any conformal map $\varphi:  D \to \wt D$ with $\varphi(x) = \wt x$ and $\varphi(y) = \wt y$, the law of $\SLE_\kappa$ in $D$ from $x$ to $y$ agrees with the pullback of the law of $\SLE_\kappa$ in $\wt D$ from $\wt x$ to $\wt y$. 
For the regime $\kappa \in (0,4]$, $\SLE_\kappa$ is a.s.\ simple and does not hit the boundary of $D$ except at $x$ and $y$.  

As explained in Section~\ref{subsec-main-results}, for $\rho_L, \rho_R > -2$ we can define a variant called $\SLE_\kappa(\rho_L; \rho_R)$. These curves are still a.s.\ simple, but when $\rho_L$ or $\rho_R$ is less than $\frac\kappa2 - 2$ the random curve a.s.\ hits the corresponding boundary arc.

When $\kappa = \gamma^2$, the $\gamma$-LQG length of $\SLE_\kappa$-type curves can be defined via conformal welding \cite{shef-zipper}, or equivalently as a Gaussian multiplicative chaos measure on the measure defined by the Minkowski content \cite{benoist-lqg-chaos}. The quantum length of a curve is measurable with respect to the curve-decorated quantum surface.

 We now inductively define the measure on curves $\cP^\disk(W_1, \dots, W_n)$ featured in Theorem~\ref{thm-disk-cutting}.
\begin{definition}[Multiple SLE]\label{def-cP}
Consider a simply connected domain $D$ with boundary marked points $x,y$, and weights $W_1, \dots, W_n >0$ for some $n \geq 2$. 
For $n = 2$, we define $\cP^\disk(W_1, W_2)$ to be $\SLE_\kappa(W_1-2; W_2 -2)$ in $(D,x,y)$. To define the probability measure $\cP^\disk(W_1, \dots, W_n)$ on curves $(\eta_1, \dots, \eta_{n-1})$ for $n \geq 3$, we first sample $\eta_{n-1}\sim \SLE_\kappa(W_1+\dots + W_{n-1} -2;W_n-2)$ in $(D,x,y)$, then for each connected component $(D',x',y')$ of $D \backslash \eta_{n-1}$ lying to the left of $\eta_{n-1}$ (with marked points the first and last points visited by $\eta_{n-1}$), we independently sample $(n-2)$ curve segments from $\cP^\disk(W_1, \dots, W_{n-1})$, and concatenate them to obtain $n-2$ curves $\eta_1, \dots, \eta_{n-2}$. 
\end{definition}
By the conformal invariance of $\SLE_\kappa(\rho_1;\rho_2)$ curves, the measure $\cP^\disk(W_1, \dots, W_n)$ is also  conformally invariant. 

We now state the quantum wedge welding theorem of \cite{wedges}, which should be compared to Theorem~\ref{thm-disk-cutting}. Although \cite[Theorem 1.2]{wedges} is stated only for the $n=2$ case, the general case is not hard to derive from the $n=2$ case and is used in, e.g., \cite[Appendix B]{wedges}, and we explicitly describe it here for the reader's convenience. Here, $\cM^\wed(W) \otimes \cP^\disk(W_1, \dots, W_n)$ is a measure on curve-decorated surfaces.
\begin{theorem}[Conformal welding of quantum wedges \cite{wedges}]\label{thm-wedges}
Consider weights $W_1, \dots, W_n>0$ and $W = \sum_i W_i$.  Then 
\[
\cM^\wed(W) \otimes \cP^\disk(W_1, \dots, W_n)
= \cM^\wed(W_1) \times \cM^\wed(W_2)\times \cdots \times \cM^\wed(W_n).
\]
\end{theorem}

Finally we define the analogous probability measure $\cP^\sph$ for curves between two marked points in a sphere with conformal structure. We state the definition for $(\wh \C, 0, \infty)$, where $\wh \C = \C \cup \{\infty\}$ is the Poincar\'e sphere. 
\begin{definition}[Multiple SLE on sphere]
For $n \geq 1$ and $W_1, \dots, W_n >0$, the probability measure $\cP^\sph (W_1, \dots, W_n)$ on $n$-tuples of curves $(\eta_0, \dots, \eta_{n-1})$ in $\wh \C$ from $0$ to $\infty$ is defined as follows. First sample $\eta_0$ as a whole-plane $\SLE_\kappa((\sum_{j=1}^n W_j) - 2)$ process from $0$ to $\infty$, then sample an $(n-1)$-tuple of curves from $\cP^\disk(W_1, \dots, W_n)$ in each connected component of $\wh \C \backslash \eta_0$, and concatenate them to get $\eta_1, \dots, \eta_{n-1}$.
\end{definition}

\subsection{Outline of proofs}\label{subsec-outline}
We now outline the proof of Theorem~\ref{thm-disk-cutting}. The proof of Theorem~\ref{thm-sphere-cutting} is similar and discussed in Section~\ref{sec-general}. 

We start with a thick quantum wedge $(\cS, h, +\infty, -\infty)$ embedded so that neighborhoods of $+\infty$ have finite quantum boundary length,  decorated by independent curves $(\eta_1, \dots, \eta_n)$ which cut it into independent thin quantum wedges. 

In Section~\ref{sec-extrinsic} we define a ``field bottleneck'' event which, roughly speaking, says that when we explore the field from left to right and stop when the field average process $h_t$ first takes a large negative value, then the quantum lengths of the unexplored boundary segments and curves are macroscopic. We show that conditioned on the existence of the field bottleneck, the unexplored region resembles a thick quantum disk decorated by curves, conditioned on having macroscopic boundary arcs and interfaces. 

In Section~\ref{subsec-thin-decomp} we show that pinching a thin quantum wedge yields a thin quantum disk, which allows us to define a ``curve bottleneck'' event in  Section~\ref{sec-intrinsic}. This event  roughly says that, letting $z_0 \in \R$ be the point such that the quantum length $\nu_h([z_0,\infty))$ is 1, certain curve segments of $\eta_1, \dots, \eta_{n-1}$ near $z_0$ are short, and the curve lengths to the right of $z_0$ are macroscopic. When we condition on this event, the region to the right of the curve bottleneck resembles a welding of thin quantum disks. We prove that the field bottleneck and curve bottleneck are compatible in a certain sense, and hence conclude that a thick curve-decorated quantum disk with macroscopic interfaces, cut along its curves, yields a collection of thin quantum disks with macroscopic side lengths. 

The arguments of Sections~\ref{sec-extrinsic}--\ref{sec-intrinsic} yield a weaker version of Theorem~\ref{thm-disk-cutting}. In Section~\ref{sec-general} we bootstrap this to the full Theorem~\ref{thm-disk-cutting} using a short and relatively easy argument.
	
\section{Pinching a thick quantum wedge yields a quantum disk} \label{sec-extrinsic}
The goal of this section is to prove Proposition~\ref{prop-extrinsic}, which for $W > \frac{\gamma^2}2$ constructs a curve-decorated weight $W$ quantum disk from a curve-decorated weight $W$ quantum wedge. It does so by identifying a ``field bottleneck'' when a quantum wedge is explored from its infinite end to its finite end, then conditioning on the surface to the right of the bottleneck being large; in the limit this pinched surface converges to a quantum disk. More strongly, Proposition~\ref{prop-extrinsic} identifies the law of the triple (field at the bottleneck, boundary arc length of pinched surface, field and curves in the bulk of pinched surface); this information will be used in Section~\ref{sec-intrinsic} to show that the field bottleneck is compatible with the ``curve bottleneck'' introduced there. 

The limit surface will be $\cM_2^\disk(W;1)$ with some conditioning on curves, where $\cM_2^\disk(W;\ell)$ is defined as follows. 
\begin{definition}[Disintegration on one boundary length]\label{def-disint-one}
 We define 
 \[\cM_2^\disk(W; \ell) := \int_0^\infty \cM_2^\disk(W; \ell, \ell')\, d \ell',\] 
i.e. we only disintegrate on the left boundary arc length. 
\end{definition}
Lemmas~\ref{lem-length-exponent} and~\ref{lem-thin-disk-perim} tell us that when $W \in (0, \gamma Q)$ the measure $\cM_2^\disk(W; \ell)$ is finite, and hence so is $\cM_2^\disk(W; \ell, \ell')$ for any $\ell'>0$. Conversely when $W \geq \gamma Q$, the measure $\cM_2^\disk(W; \ell)$ is infinite but $\sigma$-finite.

Recall that doubly-marked quantum surfaces embedded in $(\cS, +\infty,-\infty)$ have a field that is determined modulo translation: $(\cS, h, +\infty,-\infty)$ and $(\cS, h(\cdot - t), +\infty,-\infty)$ are equivalent as quantum surfaces. We say the \emph{canonical embedding} of $(\cS, h, +\infty,-\infty)$ is the embedding where $\nu_{h}(\R_+) = \frac12$.
Recall also that for a field $h$ on $\cS$, we write $h_t$ for the average of $h$ on $[t, t+i\pi]$. 

Fix $W > \frac{\gamma^2}2$ and nonnegative $W_1, \dots, W_n$ with $\sum W_j = W$. 
Consider a canonically embedded weight $W$ quantum wedge $(\cS,h, +\infty, -\infty)$, and let $\tau_{-r} = \inf\{ u \: : \: h_u = -r\}$. Independently sample curves $(\eta_1, \dots, \eta_{n-1}) \sim \cP^\disk(W_1, \dots, W_n)$ in $(\cS, +\infty,-\infty)$ and write $\eta_n := \R + i\pi$.

For $r,\zeta,K,\eps>0$ define
\begin{align}
\label{eq-E}
E_{r, K} &= \{1 \leq \nu_h(\R_+ + \tau_{-r}) \leq 1 + e^{\frac\gamma2(-r+K^3)} \}, \\
\sigma_{-\zeta} &= \inf \{ t > \tau_{-r}\: : \: X_t = -\zeta\}
, \\
\label{eq-E'}
E'_{r,\zeta,\eps} &= \{ \sigma_{-\zeta} < \infty \text{ and } \nu_h(\eta_j \cap (\cS_+ + \sigma_{-\zeta}+1)) > \eps \text{ for } j = 1, \dots, n\},
\end{align}
and define the \emph{field bottleneck event}
\eqb\label{eq-BN}
\mathrm{BN}_{r, \zeta, K, \eps} := E_{r, K} \cap E'_{r, \zeta, \eps}.
\eqe
\begin{proposition}[Surface description given $\mathrm{BN}_{r,\zeta,K,\eps}$]\label{prop-extrinsic}
Define the rectangle $R = [0, S] \times [0,\pi]$ for some $S>0$.
Then for fixed $S,K,\eps$, and $U \subset \cS$ a neighborhood of $+\infty$ excluding $-\infty$, as $r \to \infty$ then $\zeta \to \infty$ the following two probability measures have total variation distance $o(1)$:
\begin{itemize}
\item 
The law of $(h(\cdot + \tau_{-r})|_R, \nu_h(\R_+ + \tau_{-r}), 
(h|_U, \eta_1\cap U, \dots, \eta_{n-1}\cap U))$ conditioned on $\mathrm{BN}_{r,\zeta,K,\eps}$;

\item The law of the mutually independent triple $(\phi - r, V, (\wh \psi|_U, \wh \eta_1 \cap U, \dots, \wh \eta_{n-1}\cap U))$, whose components we now define. The field $\phi$ is given by 
\eqb\label{eq-phi}
\phi = \left( \wh h + (Q - \beta) \Re (\cdot) \right) \Big|_R
\eqe
where $\wh h$ is a Neumann GFF on $\cS$ normalized to have mean zero on $[0,i\pi]$, and $\beta  = \frac\gamma2 + Q - \frac W\gamma$. The random variable $V$ is sampled from $\mathrm{Unif}([1, 1+e^{\frac\gamma2(-r+K^3)}])$.

For the last field-curves tuple, take $(\cS, \wh \psi, +\infty,-\infty, \wh\eta_1, \dots, \wh\eta_{n-1})$ a canonically embedded sample from $\cM_2^\disk(W;1) \otimes \cP^\disk(W_1, \dots, W_n)$ conditioned on $\nu_{\wh \psi}(\wh \eta_j)>\eps$ for $j=1,\dots, n$.
\end{itemize}
\end{proposition}

\begin{figure}[ht!]
	\begin{center}
\includegraphics[scale=0.75]{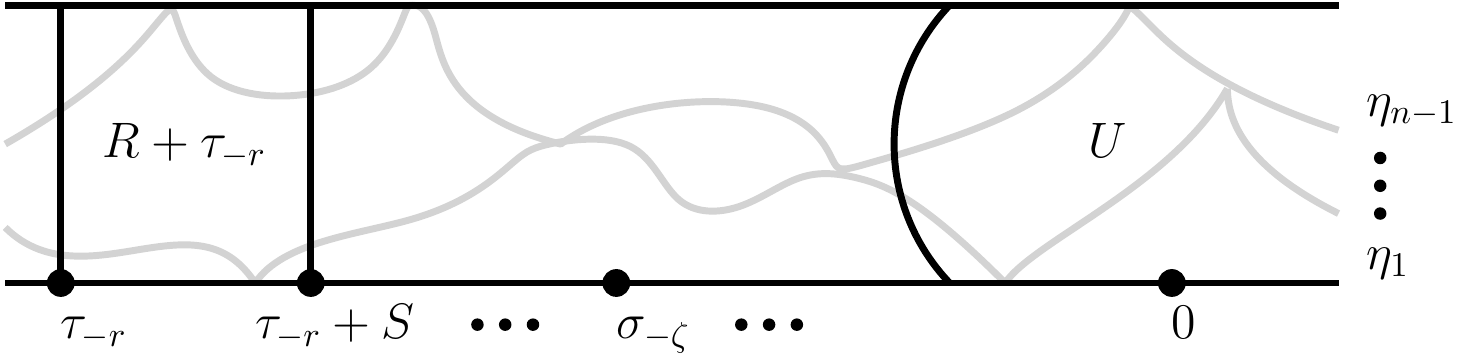}%
	\end{center}
	\caption{\label{fig-extrinsic} Setup for Proposition~\ref{prop-extrinsic}. Since we are sending $r \to \infty$ then $\zeta \to \infty$, the regions $R + \tau_{-r}$ and point $\sigma_{-\zeta}$ are sent to $-\infty$, and are well separated from each other and from $U$ with high probability.}
\end{figure}

We prove the near-independence of the field near $\tau_{-r}$ and the field in $U$ so that, conditioned on $\mathrm{BN}_{r, \zeta, K, \eps}$, the ``curve bottleneck'' event in Section~\ref{sec-intrinsic} is almost measurable with respect to the field and curves near $\tau_{-r}$, so further conditioning on the curve bottleneck yields the same limit law of the field and curves in $U$.

We now explain our choice of $\mathrm{BN}_{r,\zeta,K,\eps}$. Let $z_0\in \R$ be the point such that $\nu_h([z_0, +\infty)) = 1$. In Section~\ref{sec-intrinsic} we will define a ``curve bottleneck'' near $z_0$. The choice of upper bound $1+ e^{\frac\gamma2(-r+K^3)}$ in $E_{r, K}$ means that when we condition on $E_{r,K}$, with high probability $\tau_{-r}$ is close (in the Euclidean metric) to $z_0$, so the field bottleneck is close to the curve bottleneck. This is necessary for showing compatibility of the bottlenecks. The definition of $E'_{r, \zeta, \eps}$ comprises two events: 1) the growth of the field average process $(h_t)_{t \geq \tau_{-r}}$ to the value $-\zeta$, and 2) the curves in the ``bulk'' having macroscopic length. 1) allows us to compare the quantum wedge field to that of a quantum disk via Lemma~\ref{lem-disk-beta}, and 2) is a technical condition that later allows us to work with probability measures rather than infinite measures: although $\cM_2^\disk(W; 1)$ may be an infinite measure, if we sample $(\eta_1, \dots, \eta_n)\sim \cP^\disk(W_1, \dots, W_n)$ and restrict to the event that all curves have macroscopic quantum lengths, the resulting measure is finite. 

In order to prove Proposition~\ref{prop-extrinsic}, we switch to a more convenient field $\psi$ that resembles $h$ (Lemma~\ref{lem-equivalent-description}), identify the law of the field and curves in the bulk when we condition on side length (Lemma~\ref{lem-extrinsic-bulk}), and, after further conditioning on the bulk field and curves, identify the law of the unexplored boundary arc length and field near the bottleneck (Lemma~\ref{lem-extrinsic-tip}). Combining these yields Proposition~\ref{prop-extrinsic}.

\begin{lemma}\label{lem-equivalent-description}
For fixed $\zeta > 0$, consider $(\cS, \psi, +\infty,-\infty) \sim P_\zeta$ (defined in Lemma~\ref{lem-disk-beta}). Then $h(\cdot + \tau_{-r})|_{\cS_+}$ conditioned on $\{ \sigma_{-\zeta} < \infty\}$ agrees in distribution with $\psi(\cdot + \wt \tau_{-r})|_{\cS_+}$, where $\wt \tau_{-r} := \inf\{ t \in \R \: : \: \psi_t = -r\}$. 
\end{lemma}
\begin{proof}
Note that conditioned on $\{ \sigma_\zeta < \infty\}$, the process $(h_{t + \tau_{-r}})_{t \geq 0}$ has the law of Brownian motion started at $-r$ with variance 2 and downward drift of $-(Q-\beta)t$ (with $\beta = \frac\gamma2 + Q - \frac W\gamma$), conditioned to hit $-\zeta$. By \cite[Lemma 3.6]{wedges}, this is the same law as $(\psi_{t + \wt \tau_{-r}})_{t \geq 0}$. Finally, each field $h, \psi$ has independent projections to $\cH_1(\cS)$ and $\cH_2(\cS)$, and their projections to $\cH_2(\cS)$ agree in law with the projection of a Neumann GFF on $\cS$ to $\cH_2(\cS)$. This proves the lemma.
\end{proof}

Recall from Lemma~\ref{lem-decomposition-psi} the following decomposition of $\psi$:
\[\psi = \psi_+ + \psi_- + \alpha_1 f_1,\]
(we have absorbed the term $\alpha_2 f_2$ into $\psi_+$ to simplify notation), 
and that the conditional law of $\nu_\psi([0,1])$ given $\psi_+$ has a density $g_{\psi_+}(x) dx$, where $g_{\psi_+}$ is a nonnegative bounded continuous function. Also note that $\psi_+|_{\cS_++1}$ agrees with $\psi|_{\cS_++1}$ by definition. 

\begin{lemma}[Bulk field and curves given bottleneck]\label{lem-extrinsic-bulk}
Independently of $\psi$ sample $(\eta_1, \dots, \eta_{n-1}) \sim \cP^\disk(W_1, \dots, W_n)$ and write $\eta_n = \R + i\pi$. 
The conditional law of $(\psi_+, \eta_1\cap \cS_+, \dots, \eta_{n-1}\cap \cS_+)$ given $\{\nu_{\psi} ([\wt \tau_{-r}, \infty)) \in [1, 1+e^{\frac\gamma2(-r+K^3)}]\} \cap \{ \nu_{\psi_+}(\eta_j \cap (\cS_++1)) > \eps \text{ for }j = 1, \dots, n \}$ converges in total variation as $r \to \infty$ to that of $(\psi_+, \eta_1\cap \cS_+, \dots, \eta_{n-1} \cap \cS_+)$ conditioned on $\nu_{\psi}(\R) = 1$ and $\{ \nu_{\psi_+}(\eta_j \cap (\cS_++1)) > \eps \text{ for }j = 1, \dots, n \}$.
\end{lemma}
\begin{proof}
Define the event $A:=\{ \nu_{\psi_+}(\eta_j \cap (\cS_++1)) > \eps \text{ for all }j=1,\dots, n-1\}$ and the shorthand $X := (\psi_+, \eta_1 \cap \cS_+, \dots, \eta_{n-1} \cap \cS_+)$. 

Writing $I_r = [1, 1+e^{\frac\gamma2(-r+K^3)}]$, we claim that
\[
\lim_{r \to \infty} \frac{\P[\nu_{\psi} ([\wt \tau_{-r}, \infty)) \in I_r ]}{\P[\nu_\psi (\R) \in I_r]} = 1,
\]
and moreover, conditioned on any $X = (\psi_+, \eta_1 \cap \cS_+, \dots, \eta_{n-1} \cap \cS_+)$ for which $\nu_{\psi_+}(\R_+) < 1$ and $A$ holds, we have the almost sure limit
\eqb\label{eq-extrinsic-ratio}
\lim_{r \to \infty} \frac{\P[\nu_{\psi} ([\wt \tau_{-r}, \infty)) \in I_r \mid X]}{\P[\nu_\psi (\R) \in I_r \mid X]} = 1.
\eqe

The lemma follows from these two assertions. Indeed, write $\cL$ for the law of $X$ conditioned on $\{ \nu_\psi(\R) = 1\} \cap A$, $\cL_r$ for the law of $X$ given $\{\nu_\psi(\R) \in I_r\}\cap A$ and $\wt \cL_r$ for the law of $X$ given $\{\nu_\psi([\wt \tau_{-r}, \infty)) \in I_r\} \cap A$. By Bayes' theorem we get $X$-a.s. that when $\nu_{\psi_+}(\R_+) < 1$, we have
\[\lim_{r \to \infty} \frac{d\wt\cL_r}{d\cL_r} (X) = \lim_{r \to \infty} \frac{\P[\nu_\psi([\wt \tau_{-r}, \infty]) \in I_r \mid X] }{\P[\nu_\psi(\R) \in I_r \mid X]} \cdot \frac{\P[\nu_\psi(\R) \in I_r]}{\P[\nu_\psi([\wt\tau_{-r},\infty]) \in I_r]} = 1,\]
so $\lim_{r \to \infty} \wt \cL_r = \lim_{r \to \infty} \cL_r = \cL$ where the convergence is in total variation distance.

We justify~\eqref{eq-extrinsic-ratio}; the other limit is similar. By the continuity of $g_{\psi_+}$ we have
\[\frac{\P[\nu_\psi(\R) \in I_r \mid X]}{e^{\frac\gamma2(-r+K^3)}} = \frac{\int_1^{1+e^{\frac\gamma2(-r+K^3)}}\E\left[g_{\psi_+}(x - \nu_\psi(\R \backslash [0,1]))\mid X\right]dx}{e^{\frac\gamma2(-r+K^3)}} = \E\left[g_{\psi_+}(1- \nu_\psi(\R \backslash [0,1])) \mid X\right] + o_r(1),\]
\[\frac{\P[\nu_\psi([\wt \tau_{-r}, \infty)) \in I_r \mid X]}{e^{\frac\gamma2(-r+K^3)}} = \E\left[g_{\psi_+}(1- \nu_\psi(\R \backslash [0,1]) + \nu_\psi((-\infty, \wt \tau_{-r}]))\mid X\right] + o_r(1).\]
Clearly $\nu_\psi((-\infty, \wt \tau_{-r}])$ converges to 0 in probability as $r \to \infty$, so since $g_\psi$ is bounded and continuous, we may divide one of the above equations by the other to obtain~\eqref{eq-extrinsic-ratio}.
\end{proof}

The next result is similar to \cite[Proposition 5.5]{shef-zipper}, with additional details. 
\begin{lemma}[{Field near $\tau_{-r}$ and boundary length given bottleneck}]\label{lem-extrinsic-tip}
Condition on any $\psi_+$ for which $\nu_{\psi_+}(\R_+) < 1$. We have, $\psi_+$-a.s., that the total variation distance between the following two laws goes to zero as $r \to \infty$:
\begin{itemize}
\item The law of $(\psi(\cdot + \wt \tau_{-r})|_R, \nu_\psi((\wt \tau_{-r}, \infty)))$ when we further condition on $\{\nu_\psi((\wt \tau_{-r}, \infty)) \in [1, 1+e^{\frac\gamma2(-r+K)}]\}$.
\item The law of $(\phi - r, V)$ as described in Proposition~\ref{prop-extrinsic}.
\end{itemize}
\end{lemma}
\begin{proof}
For $N>0$, we will show that for sufficiently large $r > r_0(N)$, the two laws of Lemma~\ref{lem-extrinsic-tip} are within $o_N(1)$ in total variation. Sending $N \to \infty$ then implies the desired result. Elements of this argument are similar to those of Lemma~\ref{lem-extrinsic-bulk}, so we will be brief. 

Because $g_{\psi_+}$ is bounded and continuous and the length of the interval $[1, 1+e^{\frac\gamma2(-r+K^3)}]$ goes to zero as $r \to \infty$, when we condition on $\psi_+$ and $\{\nu_{\psi} ([\wt \tau_{-r}, \infty)) \in [1, 1+e^{\frac\gamma2(-r+K^3)}]\}$ the law of the pair $(\psi|_{\cS_- - N}, \nu_{\psi} ([\wt \tau_{-r}, \infty)))$ is within $o_r(1)$ in total variation to an independent pair, and the conditional marginal law of $\nu_{\psi} ([\wt \tau_{-r}, \infty)$ is close in total variation to that of $V$. 

Fix $N > 0$. By the Markov property of the GFF we may further decompose $\psi = \psi_+ + \alpha f + \psi_N + h_N$ as a sum of mutually independent distributions. Here, $\psi_N$ is harmonic in $\cS_- - N$, and $h_N$ is a distribution supported in $\cS_- - N$ with the following description: The field average process $(h_N)_{-N -t}$ agrees in law with $(B_{2t} - (Q - \beta)t)_{ t\geq 0}$ where $(B_t)_{t \geq 0}$ is standard Brownian motion conditioned on $B_{2t} - (Q - \beta)t < 0$ for all $t>0$. The (independent) projection of $h_N$ to $\cH_2(\cS_- - N)$ agrees in law with the projection to $\cH_2(\cS_- - N)$ of a GFF in $\cS_- - N$ with Dirichlet boundary conditions on $[-N, -N + i\pi]$ and Neumann boundary conditions elsewhere. Here, $\cH_2(\cS_- - N)\subset \cH_2(\cS)$ is the subspace of functions supported in $\cS_- - N$. 

If we sample $\psi$ given $\psi_+$ and $\{\nu_{\psi} ([\wt \tau_{-r}, \infty)) \in [1, 1+e^{\frac\gamma2(-r+K^3)}]\}$, then for large $r$ the conditional law of $h_N$ is within $o_N(1)$ in total variation from its unconditioned law --- essentially, conditioned only on $\psi_+$, the length $\nu_\psi((\wt \tau_{-r}, -N))$ converges to zero in probability as $N \to \infty$, so further conditioning on $\{\nu_{\psi} ([\wt \tau_{-r}, \infty)) \in [1, 1+e^{\frac\gamma2(-r+K^3)}]\}$ weights the law of $h_N$ by a $g_{\psi_+}$-dependent factor that is uniformly bounded above and converges to 1 in probability. 

We claim that for large $r$, the law of $\psi(\cdot + \wt \tau_{-r})|_R$ conditioned on  $\psi_+$ and $\{\nu_{\psi} ([\wt \tau_{-r}, \infty)) \in [1, 1+e^{\frac\gamma2(-r+K^3)}]\}$ is within $o_N(1)$ in total variation from the law of $\phi - r$. By our earlier discussion we may resample $h_N$ from its unconditioned law (incurring an $o_N(1)$ total variation error). We have $\wt \tau_{-r} \to -\infty$ in probability as $r \to \infty$, so since $\psi_+ + \psi_N$ converges in probability to a constant function in neighborhoods of $-\infty$, and $h_N$ is independent of $\psi_+ + \psi_N$, we conclude that the law of $(h_N + \psi_+ + \psi_N)(\cdot + \wt \tau_{-r})|_R + r$ converges\footnote{This follows immediately from the following description of the (independent) projections of $h_N + \psi_+ + \psi_N$ to $\cH_1(\cS_- - N)$ and $\cH_2(\cS_- - N)$. The projection to $\cH_1(\cS_- - N)$, viewed as a stochastic process from right to left ($-N$ to $-\infty$) is Brownian motion with variance 2 and downward drift, with random starting value and conditioned to stay below $-\zeta$. The projection to $\cH_2(\cS_- - N)$ in neighborhoods of $-\infty$ is close in total variation to that of a Neumann GFF on $\cS$ (\cite[Proposition 2.4]{ag-disk}).  } in total variation to that of $\phi$ as $r \to \infty$. This (with our earlier $o_N(1)$ error) yields the claim. 

To summarize, for large $r > r_0(N)$, we know that when we condition $\psi$ on $\psi_+$ and on $\{\nu_{\psi} ([\wt \tau_{-r}, \infty)) \in [1, 1+e^{\frac\gamma2(-r+K^3)}]\}$, the law of $\psi(\cdot + \wt \tau_{-r})|_R$ is $o_N(1)$ in total variation to that of $\phi - r$, and $\nu_\psi((\wt \tau_{-r}, \infty))$ is close to independent of $\psi(\cdot + \wt \tau_{-r})|_R$ and has law close to that of $V$. We are done.
\end{proof}

\begin{proof}[Proof of Proposition~\ref{prop-extrinsic}]
Let the field $\psi$ be sampled as in Lemma~\ref{lem-equivalent-description}, and independently sample $(\eta_1, \dots, \eta_{n-1})\sim \cP^\disk(W_1, \dots, W_n)$ and write $\eta_n = \R + i\pi$.  Let $(\cS, \wh \psi, +\infty,-\infty, \wh \eta_1, \dots, \wh \eta_{n-1})$ be the canonical embedding of $(\cS,  \psi, +\infty,-\infty,  \eta_1, \dots, \eta_{n-1})$.
Consider the triple 
\eqb \label{eq-triple}
(\psi(\cdot + \wt \tau_{-r})|_{R}, \nu_{\psi}((\wt \tau_{-r} , \infty)), (\wh \psi|_U, 
\wh \eta_1\cap U, \dots, \wh \eta_{n-1}\cap U))
\eqe
 conditioned on $\{\nu_{\psi}((\wt \tau_{-r} , \infty)) \in [1,1+e^{\frac\gamma2(-r+K^3)}]\}$ and on $\{ \nu_{\psi_+}(\eta_j \cap (\cS_++1)) > \eps \text{ for }j = 1, \dots, n \}$. Note that outside an event of probability $o(1)$ as $r \to \infty, \zeta \to \infty$, the tuple $(\wh \psi|_U, 
\wh \eta_1\cap U, \dots, \wh \eta_{n-1}\cap U)$ is a function of $(\psi_+, \eta_1 \cap \cS_+, \dots, \eta_{n-1} \cap \cS_+)$. 

Combining Lemmas~\ref{lem-extrinsic-bulk} and~\ref{lem-extrinsic-tip}, we see that as $r \to \infty$, the total variation distance between~\eqref{eq-triple} and the second triple of Proposition~\ref{prop-extrinsic} goes to zero. Also, Lemma~\ref{lem-equivalent-description} says that~\eqref{eq-triple} agrees in law with the first triple of Proposition~\ref{prop-extrinsic}. These two observations prove Proposition~\ref{prop-extrinsic}. 
\end{proof}

\section{Decompositions of thin quantum wedges and disks} \label{subsec-thin-decomp}
In this section, we show that thin quantum wedges and disks can be decomposed as a certain concatenation of beaded quantum surfaces. %

In this section, for $W \in (0, \frac{\gamma^2}2)$ we write $\cM_{2, \bullet}^\disk (\gamma^2 - W)$ for the infinite measure on simply connected three-pointed quantum surfaces, obtained by first taking a two-pointed surface $\cD \sim \cM_2^\disk(\gamma^2-W)$, then taking a boundary point from quantum measure on its left boundary arc (inducing a weighting by the left boundary arc length). 

\begin{lemma}[]\label{lem-fixed-cut-mass-decomp}
Fix $T>0$ and $W \in (0, \frac{\gamma^2}2)$. The following three procedures yield the same infinite measure on $\gamma$-LQG quantum surfaces.
\begin{itemize}
\item Sample $\cD'$ from $\cM_2^\disk(W)$ conditioned on having quantum cut point measure $T$ (i.e. concatenate the quantum surfaces of a Poisson point process on $\Leb_{[0,T]} \times \cM_2^\disk(\gamma^2 - W)$). Then take a point from the quantum length measure on the left boundary arc of $\cD'$ (this induces a weighting by the left arc length).

\item  Sample $\cD'$ from $\cM_2^\disk(W)$ conditioned on having quantum cut point measure $T$, then independently take $(u,\cD^\bullet) \sim \Leb_{[0,T]} \times \cM_{2, \bullet}^\disk (\gamma^2 - W)$. Insert $\cD^\bullet$ into $\cD'$ at cut point location $u$. 

\item Take $(u,\cD^\bullet) \sim \Leb_{[0,T]}\times \cM_{2, \bullet}^\disk (\gamma^2 - W)$, then independently sample two thin quantum disks from $\cM_2^\disk(W)$ conditioned on having cut point measures $u$ and $T- u$. Concatenate the three surfaces. 
\end{itemize}
\end{lemma}
\begin{proof}
The equivalence of the first two procedures above follows immediately from \cite[Lemma 4.1]{size-biased-PPP} applied to the Poisson point process on $\mathrm{Leb}_{[0,T]} \times \cM_2(\gamma^2-W)$.

The equivalence of the second and third procedures follows from the fact that a Poisson point process on $\Leb_{[0,T]} \times \cM_2(\gamma^2-W)$ can be obtained as the union of independent Poisson point processes on $\Leb_{[0,u]} \times \cM_2(\gamma^2-W)$ and $\Leb_{[u,T]} \times \cM_2(\gamma^2-W)$.

Finally, the measure on quantum surfaces is infinite because $\cM_{2,\bullet}^\disk(\gamma^2 - W)$ is infinite (indeed, the $\cM_{2,\bullet}^\disk(\gamma^2 - W)$-law of the left boundary arc length is a power law by Lemma~\ref{lem-length-exponent}).
\end{proof}

\begin{proposition}[Decomposition of marked thin quantum wedge]\label{prop-thin-wedge-decomp}
For $W \in (0,\frac{\gamma^2}2)$ we have for some constant $c \in (0, \infty)$ that
\[\cM^\wed (W) \times \mathrm{Leb}_{\R_+} = c\cdot  \cM^\disk_2(W) \times \cM_{2, \bullet}^\disk (\gamma^2 - W) \times \cM^\wed(W). \]
That is, taking a thin quantum wedge and taking a boundary point from quantum measure on its left boundary arc yields a finite beaded surface, a finite triply-marked surface, and an infinite beaded surface; these three surfaces are independent and have the laws described above. 
\end{proposition}
\begin{proof}
Since thin quantum wedges are uniquely characterized by their components up to cut point measure $T$ (for arbitrary $T$), it suffices to consider the thin quantum wedge up to this point. When we then restrict the two measures to the event that the marked point lies in this initial part of the thin quantum wedge, they agree by Lemma~\ref{lem-fixed-cut-mass-decomp} (first and third procedures). Sending $T \to \infty$ yields the result.
\end{proof}

\begin{figure}[ht!]
\begin{center}
  \includegraphics[scale=0.75]{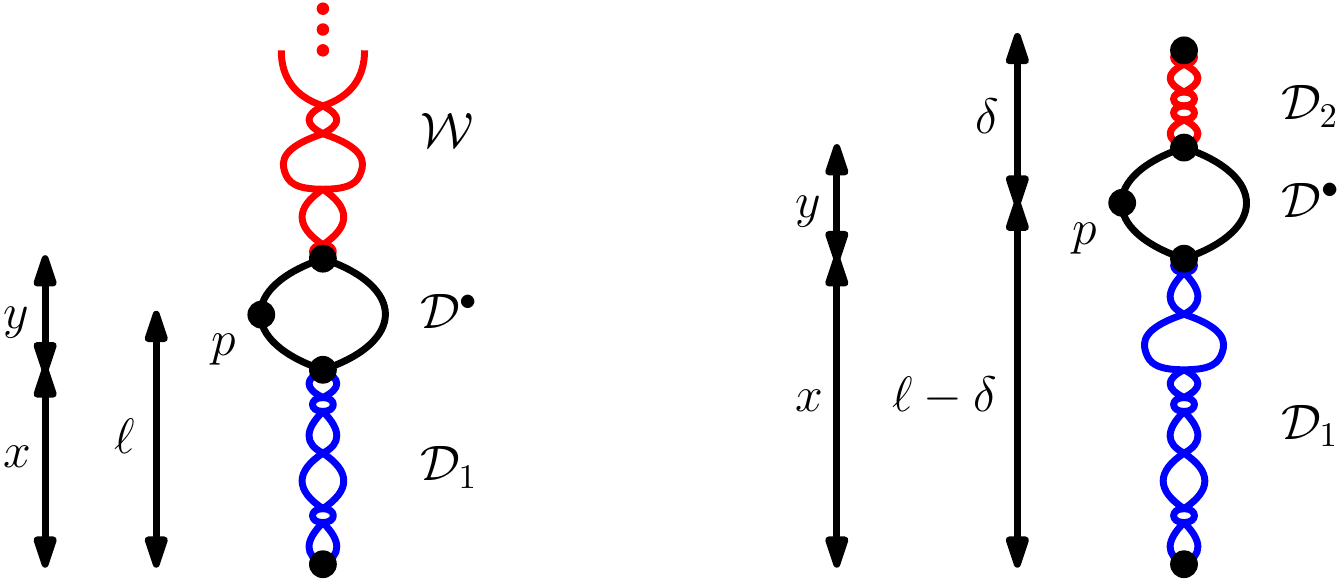}%
 \end{center}
\caption{\label{fig-thin-decomp} \textbf{Left.} In Corollary~\ref{cor-thin-wedge-decomp}, a thin quantum wedge with point $p$ on its left boundary at quantum length $\ell$ from the root decomposes as a concatenation of three quantum surfaces with some length conditioning. \textbf{Right.} In Corollary~\ref{cor-thin-disk-decomp}, a thin quantum disk conditioned to have left boundary length $\ell$ with point $p$ on its left boundary at quantum length $\delta$ from $v$ decomposes as a concatenation of three quantum surfaces with some length conditioning.}
\end{figure}

\begin{corollary}\label{cor-thin-wedge-decomp}
Fix $\ell > 0$. 
For $W \in (0,\frac{\gamma^2}2)$, the following procedures yield the same probability measure on pairs of quantum surfaces $(\cD_1,\cD^\bullet, \cW)$; see Figure~\ref{fig-thin-decomp}.
\begin{itemize}
\item Sample a thin quantum wedge $\wt \cW \sim \cM_2^\wed(W)$ and let $p$ be the point on the left boundary arc of $\wt \cW$ at distance $\ell$ from the root. Let $\cD^\bullet$ be the three-pointed disk containing $p$, and $\cD_1$ (resp. $\cW$) the finite (resp. infinite) beaded component of $\wt \cW \backslash \cD^\bullet$.
\item Take $(\cD_1, \cD, \cW) \sim   \cM_2^\disk(W)\times \cM_2^\disk(\gamma^2-W)\times \cM^\wed(W)$ and condition on the event of finite measure that the left boundary lengths $x,y$ of $\cD_1$ and $\cD$  satisfy $x < \ell < x+y$. Mark the point $p$ on $\cD$ such that the length from $p$ to the quantum disk tip is $x+y-\ell$, and call this surface $\cD^\bullet$.
\end{itemize}
\end{corollary}
\begin{proof}
As we explain, this follows from Proposition~\ref{prop-thin-wedge-decomp} by conditioning on the location of the marked point. Consider a thin quantum wedge decorated by a uniformly chosen point from its left boundary arc, and condition on the event that the left boundary interval $I_\eps$ at distance between $\ell - \eps$ and $\ell + \eps$ from the thin quantum wedge root lies on a single thick quantum disk and that the marked point lies in $I_\epsilon$. By Proposition~\ref{prop-thin-wedge-decomp} we may express this as a concatenation of quantum surfaces $(\cD_1, \cD^\bullet, \cW)\sim c\cdot  \cM^\disk_2(W) \times \cM_{2, \bullet}^\disk (\gamma^2 - W) \times \cM^\wed(W)$ conditioned on the left boundary lengths $x,y$ of $\cD_1, \cD$ satisfying $x < \ell - \eps < \ell + \eps < x+y$ and on the marked point of $\cD^\bullet$ lying in the corresponding length $2\eps$ interval. Although $\cD^\bullet$ is weighted by its left boundary length, restricting to the event that the marked point lies in an interval of length $2\eps$ removes this weighting. Sending $\eps \to 0$ thus yields the result. 
\end{proof}

\begin{proposition}[Decomposition of marked thin quantum disk]\label{prop-thin-disk-decomp}
The following two procedures yield the same measure on quantum surfaces.
\begin{itemize}
\item Take a thin quantum disk from $\cM_2^\disk (W)$, then take a boundary point from quantum measure on its left boundary arc (this induces a weighting by its left boundary length). 

\item Take a triple of quantum surfaces from $(1-\frac2{\gamma^2} W)^2\cM_2^\disk(W) \times \cM_{2, \bullet}^\disk (\gamma^2 - W) \times \cM_2^\disk(W)$ and concatenate them.  
\end{itemize}
\end{proposition}
\begin{proof}
In Lemma~\ref{lem-fixed-cut-mass-decomp} take $T \sim \Leb_{\R_+}$ then apply the first or third procedure. Here we use the equivalence of the following procedures. 
\begin{itemize}
\item First take $T \sim \Leb_{\R_+}$ then $t \sim \Leb_{[0,T]}$ (inducing a weighting by $T$), and output $(t,T)$;
\item Independently take $t, t' \sim \Leb_{\R_+}$, then output $(t, t+t')$. \qedhere
\end{itemize}
\end{proof}

As an immediate corollary, we have for any disintegration $\cM_2^\disk(W) = \int_0^\infty \wt \cM_2^\disk(W; \ell) d\ell$ with respect to the left boundary length (so $\wt \cM_2^\disk(W; \ell)$ is only uniquely defined $\ell$-a.e.) that the following holds. Recall that $M^\#$ is the normalization of $M$ to be a probability measure. 
\begin{corollary}\label{cor-thin-disk-decomp}
For Lebesgue a.e. $\ell > 0$  the following holds. 
Fix $\delta \in (0, \ell)$. 
For $W \in (0,\frac{\gamma^2}2)$, the following procedures yield the same measure on triples of quantum surfaces $(\cD_1,\cD^\bullet, \cD_2)$; see Figure~\ref{fig-thin-decomp}.
\begin{itemize}
\item Take a thin quantum disk $(\wt\cD,u,v) \sim \wt \cM_2^\disk(W, \ell)$ and let $p$ be the point on the left boundary arc of $\wt\cD$ at distance $\delta$ from $v$. Let $\cD^\bullet$ be the three-pointed disk containing $p$, and $\cD_1$, $\cD_2$ the two finite beaded components of $\cD \backslash \cD^\bullet$ (with $u \in \cD_1$). 
\item Sample $(\cD_1, \cD, \cD_2) \sim \cM_2^\disk(W) \times \cM_2^\disk(\gamma^2- W) \times \cM_2^\disk(W)$ and ``restrict to'' $\{x+y+z=\ell\}$ and $\{x < \ell -\delta < x+y\}$, where $x,y,x$ are the left boundary arc lengths of $\cD_1, \cD, \cD_2$. 
Mark the point $p$ on $\cD$ to get $\cD^\bullet$.
 
In other words, take $(x,y) \in \R_+^2$ from the measure $1_{x<\ell - \delta< x + y < \ell}\,C x^{-\frac2{\gamma^2}W} y^{-\frac2{\gamma^2}(\gamma^2 - W)} (\ell - x - y)^{-\frac2{\gamma^2}W}\, dx \, dy$, then sample $\cD_1 \sim \wt\cM_2^\disk(W; x)^\#, \cD \sim \wt\cM_2^\disk(\gamma^2 - W; y)^\#$ and $\cD_2 \sim \wt\cM_2^\disk(W; \ell - x - y)^\#$, and mark the point on the left boundary of $\cD$ to get $\cD^\bullet$.
Here, the constant $C$ is $(1-\frac2{\gamma^2} W)^2c_W^2 c_{\gamma^2-W}$, where $c_W$ and $c_{\gamma^2-W}$ are the constants of Lemmas~\ref{lem-length-exponent} and~\ref{lem-thin-disk-perim}.
\end{itemize}
\end{corollary}

The second procedure of Corollary~\ref{cor-thin-disk-decomp} does not depend on the choice of $\{\wt \cM_2^\disk(W; \ell)\}_\ell$. This gives us a way to bootstrap the disintegration $\{\wt \cM_2^\disk(W; \ell)\}_\ell$ (which is only uniquely defined $\ell$-a.e.) to a disintegration $\{\cM_2^\disk(W; \ell)\}_{\ell}$ which is well defined for \emph{all} $\ell > 0$.
\begin{definition}\label{def-thin-disint}
We define $\cM_2^\disk(W; \ell)$ to be the measure on quantum surfaces uniquely specified by the disintegration in the previous paragraph. 
\end{definition}

One can check that $\cM_2^\disk(W; \ell)$ defined in Definition~\ref{def-thin-disint} does not depend on $\delta$; this reduces to a computation on the joint law of quantum lengths arising when we mark two points at distances $\delta, \delta'$ from one of the quantum disk marked points. Moreover, as we will see in Lemma~\ref{lem-trim}, for each $\delta>0$ the measures $\{\cM_2^\disk(W;\ell)\}_{\ell>\delta}$ are continuous with respect to the total variation distance of the \emph{$\delta$-trimming} of the quantum disk. 

For a thin quantum disk $(\cD, u, v)$ with left side length greater than $\delta > 0$, we define the \emph{$\delta$-trimming} of $\cD$ to be the beaded surface obtained by marking the point on the left boundary of $\cD$ at distance $\delta$ from $v$, then discarding the beads from $v$ to this marked point (inclusive), to obtain a beaded quantum surface containing $u$. Note that this surface is a.s.\ nonempty because in the beaded quantum surface $(\cD, u, v)$, there are a.s.\ infinitely many small beads near $u$.

\begin{lemma}[Continuity properties of $\delta$-trimming]\label{lem-trim}
Fix $\ell> \delta > 0$ and $W \in (0, \frac{\gamma^2}2)$. Sample a quantum disk $(\cD,u,v) \sim \cM_2^\disk(W; \ell)^\#$, and let $\cD^\delta$ be its $\delta$-trimming. Repeat the above procedure replacing $\ell$ with $\wt\ell$ to get $\wt \cD^\delta$. Then the quantum surfaces $\cD^\delta$ and $\wt \cD^\delta$ can be coupled so that as $\wt\ell \to \ell$, we have $\cD^\delta = \wt \cD^\delta$ with probability approaching 1. Moreover, if we then take $\delta \to 0$, the left side length of $\cD^\delta$ converges in probability to $\ell$. 
\end{lemma}
\begin{proof}
By Corollary~\ref{cor-thin-disk-decomp}, the left side length of $\cD^\delta$ has probability density function
\alb
f_{\ell,\delta}(x) &:=  \frac{1_{x \in (0,\ell - \delta)}}{Z_{\ell, \delta}}\int_{\ell - \delta - x}^{\ell - x} x^{-p} y^{p-2} (\ell - x - y)^{-p} \, dy \\
&=\frac{1_{x \in (0,\ell - \delta)}}{Z_{\ell, \delta}} \frac{\delta^{1-p}}{1-p}  \frac{x^{-p}}{(\ell-x)(\ell - \delta - x)^{1-p}},
\ale
where $p = \frac2{\gamma^2}W$, and $Z_{\ell, \delta}$ is the normalization constant so that $\int_0^{\ell - \delta} f_{\ell,\delta}(x) dx = 1$. (The equivalence of the two formulae follows immediately from differentiating in $\delta$.) 

By continuity, for all $x < \ell - \delta$ we have $\lim_{\wt \ell \to \ell} f_{\wt \ell,\delta}(x) = f_{\ell,\delta}(x)$, so we can couple $\cD$ and $\wt \cD$ so that the left side lengths of $\cD^\delta$ and $\wt \cD^\delta$ agree with probability $1-o(1)$; since the conditional law of $\cD^\delta$ given its left boundary length $x$ is $\cM_2^\disk(W; x)^\#$ and likewise for $\wt\cD^\delta$, there is a coupling so $\cD^\delta = \wt \cD^\delta$ with probability $1-o(1)$. 

The second claim is clear from the  explicit formula for $f_{\ell,\delta}(x)$.  
\end{proof}

Arguing similarly we can define a disintegration $\{ \cM_2^\disk(W; \ell, \ell')\}_{\ell, \ell'}$ for \emph{all} $\ell, \ell'$. Let $\{ \wt\cM_2^\disk(W; \ell, \ell')\}_{\ell, \ell'}$ be any disintegration with respect to left and right boundary arc lengths (i.e. $\cM_2^\disk(W) = \iint_0^\infty \wt\cM_2^\disk(W; \ell, \ell')\,d \ell\, d \ell'$, and $\wt\cM_2^\disk(W; \ell, \ell')$ is supported on the set of quantum surfaces with left and right boundary lengths $\ell$ and $\ell'$ respectively).
We define $\cM_2^\disk(W;\ell, \ell')$ by taking $(x, x', y, y', z, z') \in \R_+^6$ from a suitable measure supported on the set $\{x + y + z = \ell, x' + y' + z' = \ell'\}$, then sampling $\cD_1 \sim \wt\cM_2^\disk(W; x,x')^\#, \cD \sim \wt\cM_2^\disk(\gamma^2 - W; y,y')^\#$ and $\cD_2 \sim \wt\cM_2^\disk(W; z,z')^\#$, and outputting the concatenation of $\cD_1, \cD, \cD_2$. 
The family $\{\cM_2^\disk(W;\ell, \ell')\}_{\ell, \ell'}$ satisfies a similar continuity property: for any $\delta > 0$ the measures $\{\cM_2^\disk(W;\ell,\ell') \: : \: \ell > \delta\}$  are continuous with respect to the total variation distance of the $\delta$-trimming of the quantum disk.

\section{Pinching a curve-decorated bottleneck yields thin quantum disks}\label{sec-intrinsic}

The goal of this section is to prove the following weaker version of Theorem~\ref{thm-disk-cutting}. Recall from Definition~\ref{def-disint-one} that $\{\cM_2^\disk(W;\ell)\}_{\ell}$ is the disintegration of $\cM_2^\disk(W)$ with respect to left boundary arc length.
\begin{proposition}\label{prop-weak-theorem}
Let $W_1, \dots, W_n \in (0, \frac{\gamma^2}2)$ and $W = \sum W_i > \frac{\gamma^2}2$. Then
\alb
&\cM_2^\disk(W;1)\otimes \cP^\disk(W_1, \dots, W_n) \\
&= c_{W_1, \dots, W_n} \iiint_0^\infty \cM_2^\disk(W_1; 1, \ell_1) \times \cM_2^\disk(W_2; \ell_1, \ell_2) \times \cdots \times \cM_2^\disk(W_n; \ell_{n-1}, \ell_n) \,d \ell_1 \dots d \ell_n.
\ale
\end{proposition}

We start with a curve-decorated quantum wedge $(\cS, h, +\infty, -\infty, \eta_1, \dots, \eta_{n-1}) \sim \cM^\wed(W) \otimes \cP^\disk(W_1, \dots, W_n)$ which is canonically embedded (i.e. $\nu_h(\R_+) = \frac12$).
Recall the field bottleneck event $\mathrm{BN}_{r, \zeta,K, \eps}$ from Proposition~\ref{prop-extrinsic}, and that the field bottleneck is located near $\tau_{-r} = \inf\{ t \: : \: h_t = -r\}$. Proposition~\ref{prop-extrinsic} and the SLE independence statement Lemma~\ref{lem-curve-limit} say that $(h, \eta_1, \dots, \eta_{n-1})$ conditioned on $\mathrm{BN}_{r, \zeta, K, \eps}$ in neighborhoods of $+\infty$ and near $\tau_{-r}$ are almost independent, and that conditioned on $\mathrm{BN}$, the  law of $(h,\eta_1, \dots, \eta_{n-1})$ in neighborhoods of $+\infty$ converges to that of $\cM_2^\disk(W, 1) \otimes \cP^\disk(W_1, \dots, W_n)$ conditioned on curve and boundary lengths being at least $\eps$. Based on this, our proof of Proposition~\ref{prop-weak-theorem} is carried out in four steps. 
\begin{enumerate}[Step 1.]
\item Introduce a ``curve bottleneck'' event $F_{r,K,\eps}$.  Conditioning on $F_{r,K,\eps}$, the pinched region is a conformal welding of thin quantum disks with small offsets near the bottleneck. This is done in Section~\ref{subsec-intrinsic-setup}, building on Section~\ref{subsec-thin-decomp}.

\item Conditioned on $F_{r, K, \eps}$, the conformal welding of thin quantum disks  with small offsets converges in neighborhoods of $+\infty$ as $r \to \infty$ to an exact welding of thin quantum disks. This is done in Section~\ref{subsec-intrinsic-thin-decomp}.

\item ${\P[ \mathrm{BN}_{r, \zeta, K, \eps} \mid F_{r, K, \eps}]} \to 1$ as $r\to \infty, \zeta \to \infty, K \to \infty$. Consequently, the law of $(h, \eta_1, \dots, \eta_{n-1})$  conditioned on $F_{r, K, \eps}$ is close to its law conditioned on $\mathrm{BN}_{r, \zeta, K, \eps} \cap F_{r, K, \eps}$. This is carried out in Section~\ref{subsec-intrinsic-compatibility}.

\item Conditioned on $\mathrm{BN}_{r, \zeta, K, \eps}$, the event $F_{r, K, \eps}$ occurs with uniformly positive probability for large $r,\zeta$, and is almost measurable with respect to the field and curves near $\tau_{-r}$. Combining with Proposition~\ref{prop-extrinsic}, we conclude that $(h, \eta_1, \dots, \eta_{n-1})$ conditioned on $\mathrm{BN}_{r, \zeta, K, \eps} \cap F_{r,K,\eps}$ converges in neighborhoods of $+\infty$ to a curve-decorated thick quantum disk. Comparing with Steps 2 and 3 yield the theorem. This is done in Section~\ref{subsec-intrinsic-thick}.
\end{enumerate}

\subsection{Conditioning on the curve bottleneck event $F_{r,K,\eps}$}\label{subsec-intrinsic-setup}
In this section we define $F_{r,K,\eps}$ and discuss the laws of the quantum surfaces arising upon conditioning on $F_{r, K, \eps}$. 

\begin{figure}[ht!]
\begin{center}
  \includegraphics[scale=0.75]{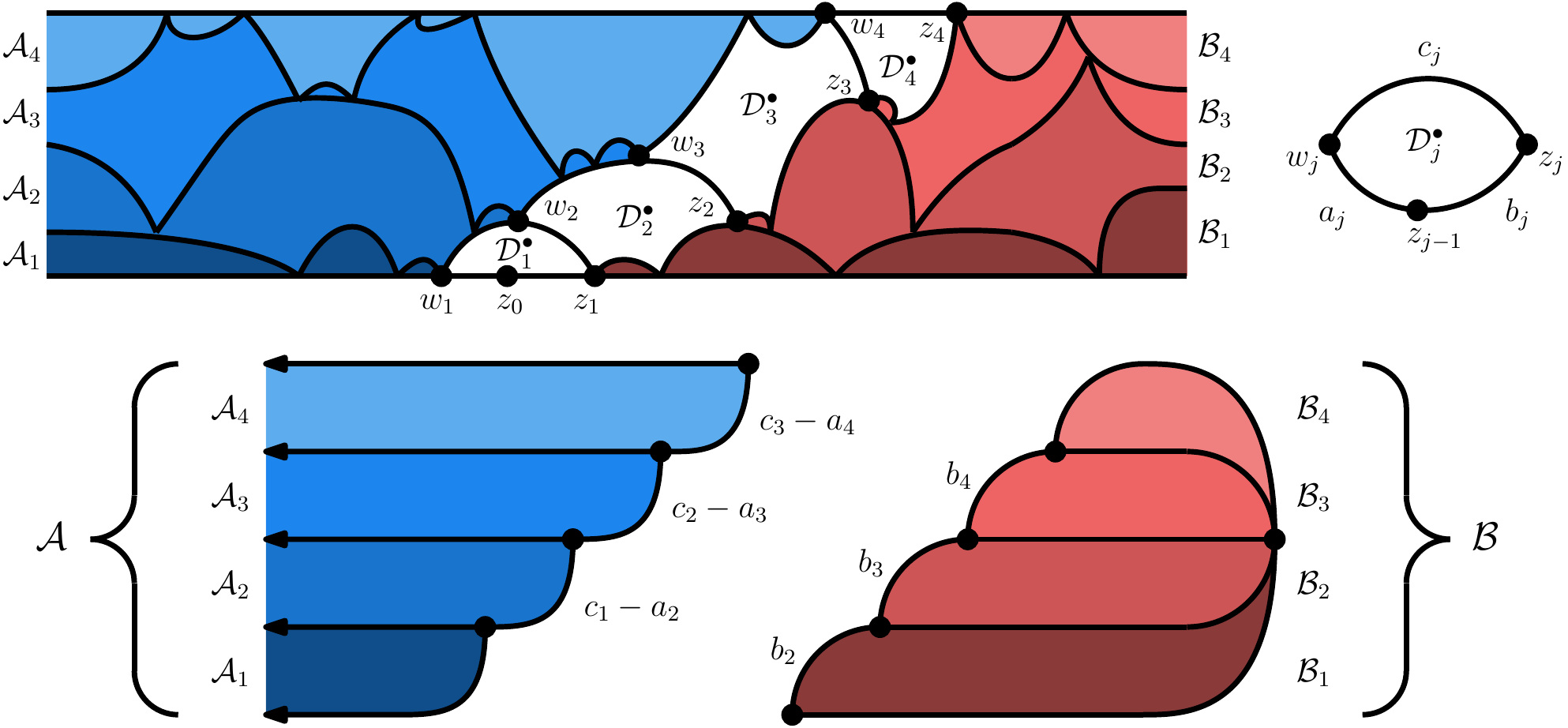}%
 \end{center}
\caption{\label{fig-bottleneck} \textbf{Top left.} Starting with the point $z_0 \in \R$ satisfying $\nu_h((z_0, \infty)) = 1$, for $j = 1, \dots, n$ we iteratively define the thin wedge bead $\cD^\bullet_j \subset \cW_j$ containing $z_{j-1}$, and call its right marked point $z_j$. Removing $\cD_j$ from $\cW_j$ yields an infinite beaded quantum surface $\cA_j$ and a finite beaded quantum surface $\cB_j$. \textbf{Top right.} Each $\cD^\bullet_j$ has three marked boundary points $z_{j-1}, z_j, w_j$ and arc lengths $a_j, b_j, c_j$ (although $w_{j+1}$ lies on $\partial \cD^\bullet_j$ we don't treat it as a marked point). \textbf{Bottom.} Conditioned on $F_{r,C,\eps}$ and on the lengths $(a_j,b_j,c_j, \ell_j)_{j=1}^n$, the multiply-marked quantum surfaces $\cA := \bigcup \cA_j$ and $\cB := \bigcup \cB_j$ are independent. These surfaces are welded with offsets shown in the diagram (note that the nontrivial topology of the quantum disks and wedges are not depicted).}
\end{figure}

We start with the definition of $F_{r, K, \eps}$; see Figure~\ref{fig-bottleneck}.
Let $(\cS, h, +\infty, -\infty)$ be a quantum wedge $\cW$ with weight $W$, and sample $(\eta_1, \dots, \eta_{n-1})\sim \cP^\disk(W_1, \dots, W_n)$ on $\cS$ with $\eta_1$ on the bottom and $\eta_{n-1}$ on top, cutting $\cW$ into independent quantum wedges $\cW_1, \dots, \cW_n$ of weights $W_1, \dots, W_n$ (Theorem~\ref{thm-wedges}). Let $z_0 \in \R$ be the point such that $\nu_h([z_0, \infty)) = 1$, and let $\cD^\bullet_1$ be the thick quantum disk of $\cW_1$ containing $z_0$; call its left and right marked points $w_1, z_1$ respectively. Iteratively for $j = 1, \dots, n$, let $\cD^\bullet_j$ be the thick quantum disk of $\cW_j$ containing $z_{j-1}$ on its boundary, and let $w_j, z_j$ be its left and right marked points. For $j = 1, \dots, n$, from $w_j$ and tracing $\partial \cD^\bullet_j$ in counterclockwise order, let the three boundary arc lengths be $a_j, b_j, c_j$, and let $\ell_j$ be the length of $\eta_j$ from $z_j$ to $+\infty$ (here we write $\eta_n := \R + i\pi$).

Finally, define
\eqb\label{eq-F-abc}
F_{r, K, \eps} := \bigcap_{j=1}^n \left\{
e^{\frac\gamma2(r+K)}b_j \in [1, 2], \quad e^{\frac\gamma2(r+K)}(a_j+b_j) \in [3, 4], \quad  e^{\frac\gamma2(r+K)}c_j \in [5, 6], \quad \ell_j > \eps \right\}.
\eqe
The intervals $[1,2], [3,4], [5,6]$ are chosen to ensure the welding offsets $c_{j-1} - a_{j}, b_j$ for $j=2,\dots, n$ described in Proposition~\ref{prop-F-decomp} are roughly $e^{\frac\gamma2(-r-K)}$ in magnitude; see Figure~\ref{fig-bottleneck}.

\begin{lemma}\label{lem-length-distribution}
Set $(\wt a_j, \wt b_j, \wt c_j) := (e^{\frac\gamma2(r+K)}a_j, e^{\frac\gamma2(r+K)}b_j, e^{\frac\gamma2(r+K)}c_j)$ for $j = 1, \dots, n$. Then in the $r \to \infty$ limit, conditioned on $F_{r,K, \eps}$ the $n+1$ tuples $(\ell_1, \dots, \ell_{n-1},\ell_n), (\wt a_1, \wt b_1, \wt c_1), \dots, (\wt a_n, \wt b_n, \wt c_n)$ jointly converge in distribution to a collection of $n+1$ independent tuples. The limit law of $(\ell_1, \dots, \ell_n)$ has density with respect to Lebesgue measure on $(\eps, \infty)^n$ given by 
\eqb\label{eq-l-limit}
\frac1Z \prod_{j=1}^{n} \left|\cM_2^\disk(W_j; \ell_{j-1}, \ell_j)\right|  \quad\text{ with } \ell_0 := 1,
\eqe
and for each $j=1, \dots, n$, the limit law of $(\wt a_j, \wt b_j, \wt c_j)$ is supported on the set $S_j = \{\wt b_j \in [1,2], (\wt a_j+\wt b_j) \in [3,4], \wt c_j \in [5,6]\}$ and has density with respect to Lebesgue measure
\[\frac1{Z_j} 1_{(\wt a_j, \wt b_j, \wt c_j) \in S_j} \left|\cM_2^\disk(\gamma^2 - W_j; \wt a_j + \wt b_j, \wt c_j) \right|.\]
Here the $Z, Z_1, \dots, Z_n$ are nonexplicit normalization constants. 
\end{lemma}
Implicit in the above lemma is the fact that the integral $\iiint_\eps^\infty \prod_{j=1}^n\left| \cM_2^\disk(W_j, \ell_{j-1}, \ell_j) \right| \,d \ell_1 \dots d\ell_n$ is finite. We show this for $n=2$, and the general case follows similarly. Using Lemma~\ref{lem-scaling-disint},
\alb
\iint_\eps^\infty &\left|\cM_2^\disk(W_1; 1, \ell_1)\right| \left| \cM_2^\disk(W_2; \ell_1, \ell_2)\right|\, d\ell_2 \,d \ell_1 = \int_\eps^\infty \left|\cM_2^\disk(W_1; 1, \ell_1)\right| \ell_1^{-\frac2{\gamma^2}W_2} \int_{\eps/\ell_1}^\infty \left| \cM_2^\disk(W_2; 1, x)\right|\, dx \,d \ell_1 \\
<& \ \eps^{-\frac2{\gamma^2}W_2}\left|\cM_2^\disk(W_2; 1)\right| \int_\eps^\infty \left|\cM_2^\disk(W_1; 1, \ell_1)\right|\, d \ell_1 < \eps^{-\frac2{\gamma^2}W_2}\left|\cM_2^\disk(W_2; 1)\right| \left|\cM_2^\disk(W_1; 1)\right| < \infty.
\ale
\begin{proof}
Although $\ell_0 = 1$ is a constant, we will make statements in terms of $\ell_0$ that generalize to $\ell_j$. We will also slightly abuse notation and use the same symbol for random variables and the dummy variables describing their densities.

We first explain the law of $(\wt a_1, \wt b_1, \wt c_1, \ell_1)$ when we condition on~\eqref{eq-F-abc} for $j=1$.  
 Start with the unconditioned setup, and define $x := \ell_0 - b_1$ and $y = a_1 + b_1$. By Corollary~\ref{cor-thin-wedge-decomp} and Lemmas~\ref{lem-length-exponent} and~\ref{lem-thin-disk-perim}, the law of $(x,y)$ is given by
\[\wh Z_1^{-1} 1_{0<x<\ell_0<x+y} x^{-p_1} y^{p_1-2}\,dx\,dy, \quad \text{ with } p_1 := \frac2{\gamma^2} W_1 \in (0,1). \]
Indeed, $x$ and $y$ are the left side lengths of $\cB_1 \sim \cM_2^\disk(W_1)$ and $\cD_1 \sim \cM_2^\disk(\gamma^2-W_1)$ conditioned on $x < \ell_0 < x+y$. Moreover, by doing a change of variables $(x,y) \mapsto (\alpha x, \alpha y)$, we see that the normalization constant 
$\wh Z_1 =\iint 1_{0<x<\ell_0<x+y} x^{-p_1} y^{p_1-2}\,dx\,dy$ 
has no dependence on $\ell_0$; this is important for subsequent steps where $\ell_j$ is random. 

A change of variables yields that when we condition on the event that  $\wt b_1 \in [1,2]$ and $(\wt a_1+ \wt b_1) \in [3,4]$, the conditional law of $(\wt a_1, \wt b_1) = e^{\frac\gamma2(r+K)}(x+y - \ell_0, \ell_0 -x)$ has density given by a $W_1$-dependent constant times $(1+o_r(1))\ell_0^{-p_1}(\wt a_1 + \wt b_1)^{p_1 - 2}\,d\wt a_1\, d \wt b_1$ on its support. The conditional law of $\wt c_1$ given $\wt a_1 + \wt b_1$ is that of the right boundary length of $\cM_2^\disk(\gamma^2 - W_1; \wt a_1 + \wt b_1)^\#$, and by Lemma~\ref{lem-length-exponent} (with the weight $W = \gamma^2 - W_1$) can be written as
\[\frac{\left|\cM_2^\disk(\gamma^2 - W_1; \wt a_1 + \wt b_1, \wt c_1)\right|}{\left|\cM_2^\disk(\gamma^2 - W_1; \wt a_1 + \wt b_1)\right|}\, d \wt c_1 = \frac{\left|\cM_2^\disk(\gamma^2 - W_1; \wt a_1 + \wt b_1, \wt c_1)\right|}{(\wt a_1 + \wt b_1)^{p_1-2}\left|\cM_2^\disk(\gamma^2 - W_1; 1))\right|} \,d \wt c_1.\]
Similarly, since $x = (1-o_r(1)) \ell_0$ and using Lemma~\ref{lem-thin-disk-perim}, the conditional law of $\ell_1$ given $(\wt a_1, \wt b_1)$ is some $W_1$-dependent constant times $(1+o_r(1))\ell_0^{p_1} \left|\cM_2^\disk(W_1;\ell_0, \ell_1 )\right| \,d \ell_1$. 

By the conditional independence of $\cD_1$ and $\cB_1$ given $\wt a_1, \wt b_1, x$, when we further condition on $\wt c_1 \in [5,6]$ the density of $(\wt a_1, \wt b_1, \wt c_1, \ell_1)$ is a $W_1$-dependent constant times 
\[(1+o_r(1))\left|\cM_2^\disk(\gamma^2 - W_1; \wt a_1 + \wt b_1, \wt c_1)\right| \left|\cM_2^\disk(W_1; \ell_0, \ell_1 )\right|\,d\wt a_1\, d\wt b_1\, d\wt c_1\, d \ell_1.\]

We now understand the law of $(\wt a_1, \wt b_1, \wt c_1, \ell_1)$ when we condition on~\eqref{eq-F-abc} for $j=1$. Iterating for $j = 2, \dots, n$ and using the independence of $\cW_1, \dots, \cW_n$ yields the lemma.
\end{proof}

\begin{proposition}\label{prop-F-decomp}
On the event $F_{r,K,\eps}$, condition on the lengths $(a_j, b_j, c_j, \ell_j)_{j=1}^n$. Then the surfaces $\cB$ and $\cA$ to the right and left, respectively, of $\bigcup_j \cD^\bullet_j$ are a.s.\ conditionally independent. The conditional law of $\cB$ is given by the welding of independent thin quantum disks $\cB_j \sim \cM_2^\disk(W_j; \ell_{j-1} - b_j, \ell_j)^\#$ for $j = 1,\dots, n$, where $\ell_0 := 1$. The conditional law of $\cA$ is given by the welding of independent thin quantum wedges $\cA_j \sim \cM_2^\wed(W_j)$, where the root of $A_j$ is welded to the point on the left boundary of $A_{j+1}$ at distance $c_j - a_{j+1}$ from the root for $j = 1, \dots, n-1$.  See Figure~\ref{fig-bottleneck}.
\end{proposition}
\begin{proof}
This is immediate from Corollary~\ref{cor-thin-wedge-decomp}.
\end{proof}

\subsection{Convergence to welding of thin quantum disks}\label{subsec-intrinsic-thin-decomp}
In this section we prove Proposition~\ref{prop-bulk-law}, which roughly says that when we condition on $F_{r,K,\eps}$ and send $r \to \infty$, the surface $(\cS_+ + \tau_{-r}, h, +\infty, \tau_{-r})$ converges in distribution to a welding of thin quantum disks, with respect to a suitable topology on quantum surfaces. Although $\cB = \bigcup \cB_j$ is a welding of quantum disks whose side lengths do not exactly agree, using Lemma~\ref{lem-trim} it can be coupled to agree (with high probability, near $+\infty$) with a welding of quantum disks whose side lengths \emph{do} agree, yielding Proposition~\ref{prop-bulk-law}.

For a quantum surface $(\cS, h, +\infty,-\infty, \eta_1, \dots, \eta_n)$  embedded in the strip and satisfying $\nu_\psi(\R) > \frac12$, recall from Section~\ref{sec-extrinsic} its  \emph{canonical embedding} satisfies $\nu_{h}(\R_+) = \frac12$. 

\begin{proposition}\label{prop-bulk-law}
Condition $(\cS, h, +\infty,-\infty, \eta_1, \dots, \eta_{n-1})$ on $F_{r,K,\eps}$ and consider its canonical embedding. As $r \to \infty$, in any neighborhood of $+\infty$ excluding $-\infty$, the field and curves converge in distribution to those of the canonical embedding of a sample from
\eqb\label{eq-weak-thm-eps}
Z^{-1}
\iiint_\eps^\infty \cM_2^\disk(W_1; 1, \ell_1) \times \cM_2^\disk(W_2; \ell_1, \ell_2)\times \cdots\times \cM_2^\disk(W_n; \ell_{n-1}, \ell_n)\, d \ell_1 \dots d \ell_n,
\eqe
where $Z$ is a normalization constant, and we understand~\eqref{eq-weak-thm-eps} as a probability measure on field-curves tuples in $\cS$ obtained by conformally welding $n$ quantum surfaces then canonically embedding the resulting curve-decorated surface in $(\cS,+\infty,-\infty)$. 
The topology of convergence is, for each neighborhood of $+\infty$ not containing $-\infty$, the product topology of the weak-* topology for fields and Hausdorff topology for curves. 

Moreover, we have $\P[E'_{r, \zeta,\eps} \mid F_{r, K, \eps}] \to 1$ for fixed $K$ as first $r \to \infty$ then $\zeta \to \infty$;  the event $E'_{r, \zeta, \eps}$ is defined in~\eqref{eq-E'}. 
\end{proposition}
\begin{proof}
We first elaborate on the definition and well-definedness of~\eqref{eq-weak-thm-eps} as a measure on field-curve tuples; a sample from~\eqref{eq-weak-thm-eps} is obtained as follows.
Fix $\wt \ell_0 = 1$ and sample $\wt \ell_1, \dots, \wt \ell_{n-1}$ from the distribution~\eqref{eq-l-limit}. Sample independent quantum disks $\wt \cB_1, \dots, \wt \cB_n$ with $\wt \cB_j\sim \cM_2^\disk(W_j; \wt \ell_{j-1}, \wt \ell_j)^\#$. Conformally weld them by quantum length to obtain a quantum surface $\wt B$ with two marked points and $n-1$ curves, and embed $\wt B = (\cS,\wt h, +\infty,-\infty, \wt \eta_1, \dots, \wt \eta_{n-1})$ via the canonical embedding. 
The a.s.\ existence and uniqueness of this conformal welding follows from that of thin quantum wedges and the local absolute continuity of thin quantum disks with respect to thin quantum wedges. 
\medskip

\noindent \textbf{Proving convergence to~\eqref{eq-weak-thm-eps}.}
Consider a parameter $\delta > 0$; we will send $r \to \infty$ then $\delta \to 0$ in that order, and write $o_r(1)$ (resp. $o_\delta(1)$) for a quantity that tends to zero in probability as $r \to \infty$ (resp. $r\to \infty$ then $\delta \to 0$). 
Sample $\cB_1, \dots, \cB_n$ conditioned on $F_{r, K, \eps}$ and let $\cB_1^\delta, \dots, \cB_n^\delta$ be the $\delta$-trimmings of $\cB_1, \dots, \cB_n$ (so each $\cB_j^\delta$ contains the marked point $+\infty$). Similarly let $\wt \cB_1^\delta, \dots, \wt \cB_n^\delta$ be the $\delta$-trimmings of $\wt \cB_1, \dots, \wt \cB_n$. By Lemma~\ref{lem-trim}, we can couple $(\cB_1^\delta, \dots, \cB_n^\delta) = (\wt \cB_1^\delta, \dots, \wt \cB_n^\delta)$ with probability $1-o_r(1)$. Restrict to this event.

 Let $\wt D$ (resp. $D$) be the region parametrizing $\bigcup_1^n \wt \cB_j^\delta$ (resp. $\bigcup_1^n \cB_j^\delta$). Since $\bigcup_1^n \wt \cB_j^\delta = \bigcup_1^n \cB_j^\delta$ as quantum surfaces, there is a.s.\  a (random) conformal map $\varphi: \wt D \to D$ fixing $+\infty$ so that $h|_D = (\wt h \circ \varphi^{-1} + Q \log |(\varphi^{-1})'|)|_D$. Since $h|_D \in H^{-1}_\mathrm{loc}(D)$ we conclude that $\wt h|_{\wt D} \in H^{-1}_\mathrm{loc}(\wt D)$ also. When we send $\delta \to 0$, the trimming interface in $(\cS, \wt h, +\infty,-\infty)$ goes to $-\infty$ in probability. Therefore \cite[Lemma 2.24]{ag-disk} says that for any neighborhood $U$ of $+\infty$ bounded away from $-\infty$, we have $\sup_{z \in U} e^{\Re z}|\varphi'(z) - 1| = o_\delta(1)$ (the cited lemma only states boundedness of $\sup_{z \in U} |\varphi'(z)-1|$, but the argument gives exponential decay); consequently there is a random constant $c$ for which $\sup_{z \in U} |\varphi (z) - z + c| = o_\delta(1)$. Since both $h$ and $\wt h$ are canonically embedded we have $c = o_\delta(1)$, hence 
\eqb\label{eq-phi-almost-identity}
\sup_{z \in U} |\varphi(z) - z| = o_\delta(1).
\eqe 
This allows us to show that as $r \to \infty$ then $\delta \to 0$, the tuple $(h,\eta_1, \dots,  \eta_{n-1})$ converges to $(\wt h, \wt \eta_1, \dots, \wt \eta_{n-1})$ in distribution: Convergence of the curves in the Hausdorff topology is immediate from~\eqref{eq-phi-almost-identity}. For the field, notice that for $f$ a smooth function compactly supported in $U$ we have
\[(h,f)_\nabla = (\wt h \circ \varphi^{-1} + Q \log |(\varphi^{-1})'|, f)_\nabla = (\wt h, f \circ \varphi)_\nabla + o_\delta(1) = (\wt h, f)_\nabla + o_\delta(1), \]
since $\sup_U \left| \log |(\varphi^{-1})'|\right| = o_\delta(1)$ and $f \circ \varphi \to f$ in probability in the $C^1$ topology. Since this holds for all $f$ we obtain convergence in distribution of the field (in the weak-* topology).
\medskip

\noindent \textbf{Showing $\P[E'_{r, \zeta} \mid F_{r, K, \eps}]\to 1$ as $ r\to\infty$ then $\zeta \to \infty$.} Choose some (random) $x \ll 0$ such that $\nu_{\wt h} (\wt \eta_j \cap (\cS_+ + x)) > \eps$ for $j=1,\dots,n$. 
Let $\Phi_{b}^1$ be the set of smooth functions supported in the rectangle $[x-4,x-1]\times [0,\pi]$ with $\phi \geq 0, \int \phi(x)\,d^2x =1$ and $\|\phi'\|_\infty \leq b$ and define
\[
m(\wt h) := \inf_{\phi \in \Phi_{b}^1}(\wt h, \phi).
\]
Since $\wt h$ is a distribution and $\Phi^1_{b}$ is compact in the space of test functions, $m(\wt h)$ is finite almost surely (see the discussion after Proposition 9.19 in \cite{wedges} for details). Fix a nonnegative function $f$ which is constant on vertical segments, supported on $[x-3, x-2] \times [0,\pi]$, and satisfies $\int f(x)d^2x = 1$.  Then since $\sup_{z \in U}|\varphi'(z) - 1| = o_\delta(1)$ and $\sup_{z \in U} |\varphi(z) - z| = o_\delta(1)$, we conclude that for some $b$ depending only on $f$, we have $|\varphi'|^2f \circ \varphi \in \Phi^1_{b}$ in probability as $r \to \infty$ then $\delta \to 0$. Thus, if we condition on the event $\{ m(\wt h) > -\zeta + 1\}$, then with probability $1-o_r(1)$ we have 
\[(h, f) = (\wt h\circ \varphi^{-1} + Q \log |(\varphi^{-1})'|, f) = (\wt h, |\varphi'|^2f\circ \varphi) + (Q \log |(\varphi^{-1})'|, f) > -\zeta.\] 
Since $f$ is constant on vertical segments, there exists some $t < x$ for which $h_t > -\zeta$; moreover, the quantum lengths of $\eta_j \cap (\cS_+ + t)$ are at least $\eps$ so $E'_{r, \zeta}$ holds. Since $\lim_{\zeta \to \infty} \P[ m(\wt h) > -\zeta + 1] = 1$ we obtain the desired result.
\end{proof}

\subsection{Compatibility of bottlenecks} \label{subsec-intrinsic-compatibility}
In this section, we prove Proposition~\ref{prop-E-given-F}, which roughly speaking says that $\P[E_{r,K}\mid F_{r,K,\eps}] \approx 1$. This is tricky because we are conditioning on the very rare event $F_{r, K, \eps}$. On the other hand, the surface $\cA$ conditioned on $F_{r, K, \eps}$ is simply a welding of independent quantum wedges with welding offsets $\asymp e^{\frac\gamma2(-r-K)}$ (Lemma~\ref{lem-length-distribution} and Proposition~\ref{prop-F-decomp}). Let $\cA + r+K$ denote the quantum surface obtained by adding $r+K$ to the field of $\cA$. 
We define a proxy surface $\wh \cA$ so that the law of the quantum surface $\cA + r + K$ conditioned on $F_{r,K,\eps}$ is absolutely continuous with respect to the law of $\wh \cA$. We obtain estimates on $\wh\cA$ in  Lemma~\ref{lem-comparison-wh-h}, and use these to analyze $\cA$ and hence prove Proposition~\ref{prop-E-given-F}.

First we construct the proxy surface $\wh {\mcl A}$. Consider $\wh \cW = (\cS, \wh h, +\infty, -\infty)$ decorated by curves $(\wh \eta_1, \dots, \wh \eta_{n-1})\sim \cP^\disk(W_1, \dots, W_n)$, conditioned on the following event $F$: defining the point $\wh z_0, \dots, \wh z_n$, quantum surfaces $\wh \cD_j^\bullet, \wh \cA_j$, and lengths $\wh a_j, \wh b_j, \wh c_j$ in the same way as in $F_{r,K,\eps}$, set
\[ 
F = \left\{ 
\wh b_j \in [1, 2], \quad (\wh a_j+\wh b_j) \in [3, 4], \quad  \wh c_j \in [5, 6]
\quad \text{ for }j=1,\dots,n \right\}.
\]
Let $\wh \cA := \bigcup_{i=1}^n \wh \cA_i$ and let $\wh U \subset \cS$ be the unbounded connected component of the set parametrizing $\wh \cA$. 

\begin{lemma}\label{lem-proxy}
The law of the quantum surface $\cA + r + K$ conditioned on the event $F_{r, K, \eps}$ is absolutely continuous with respect to the law of $\wh \cA$, with Radon-Nikodym derivative uniformly bounded for $r,K>0$. 
\end{lemma}
\begin{proof}
By Lemma~\ref{lem-length-distribution} we know that given $F_{r,K,\eps}$, the law of $(\wt a_j, \wt b_j, \wt c_j)_{j=1}^n = (e^{\frac\gamma2(r+K)} a_j, e^{\frac\gamma2(r+K)} b_j, e^{\frac\gamma2(r+K)} c_j)_{j=1}^n$ has a density with respect to Lebesgue measure on $\prod_{j=1}^n \{\wt b_j \in [1,2], (\wt a_j+\wt b_j) \in [3,4], \wt c_j \in [5,6]\}$, and this density is bounded between $\delta$ and $\delta^{-1}$ uniformly for all $r>0$, for some $\delta \in (0,1)$. By the same reasoning, the same is true for $(\wh a_j, \wh b_j, \wh c_j)_{j=1}^n$, so the law of $(\wt a_j, \wt b_j, \wt c_j)_{j=1}^n$ conditioned on $F_{r,K,\eps}$ is absolutely continuous with respect to the law of $(\wh a_j, \wh b_j, \wh c_j)_{j=1}^n$, and the Radon-Nikodym derivative is uniformly bounded in $r$. 

Given $(\wt a_j, \wt b_j, \wt c_j)_{j=1}^n$ the conditional law of the quantum surface $\cA + r + K$ is simply a welding of independent thin quantum wedges with offsets given by $(\wt c_j - \wt a_{j+1})_{j=1}^{n-1}$ (Proposition~\ref{prop-F-decomp}), and the same is true for $\wh \cA$. This yields the lemma. 
\end{proof}

Now we establish some estimates on $\wh A = (\wh U, \wh h, -\infty)$ that hold for any choice of embedding $(U, \wh h^U,-\infty)$ of $\wh A$ in $\cS$ fixing $-\infty$. Recall that for any field $\psi$ on $\cS$ we write $\psi_t$ for the average of $\psi$ on $[t,t+i\pi]$. 
\begin{lemma}\label{lem-comparison-wh-h}
As $K \to \infty$ the following holds in probability:

For any simply connected neighborhood $U\subset \cS$ of $-\infty$ and any conformal map $\varphi: \wh U \to U$ fixing $-\infty$, writing $\wh h^U := \wh h \circ \varphi^{-1} + Q \log |(\varphi^{-1})'|$, there exists $x \in \R$ so that the segment $[x, x+i\pi]\subset U$, $\wh h^U_x \in (-K,K)$, and $\nu_{\wh h^U}((x-K^2, \infty) \cap U) < \frac12 e^{\frac\gamma2 K^3}$. 
\end{lemma}
\begin{figure}[ht!]
\begin{center}
  \includegraphics[scale=0.8]{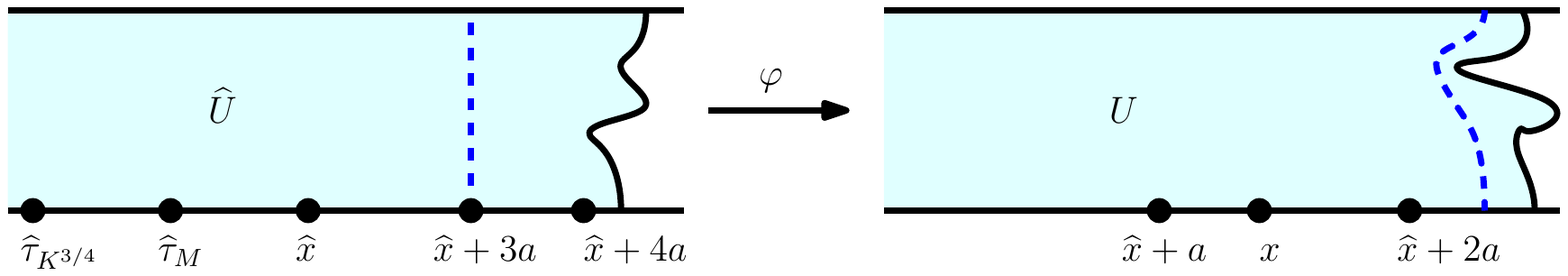}%
 \end{center}
\caption{\label{fig-proxy-surface} Diagram for argument of Lemma~\ref{lem-comparison-wh-h}.}
\end{figure}
\begin{proof}
Let $a,b>0$ be absolute constants we choose later. 
Write $\wh x= \inf\{\Re z  :  z \in \cS \backslash \wh U \} - 4a$. Let $\Phi^1_{a,b}$ be the set of smooth functions supported in the rectangle $[\wh x, \wh x + 3a]\times [0,\pi]$ with $\phi \geq 0, \int \phi(x)\,d^2x =1$ and $\|\phi'\|_\infty \leq b$ and define
\[
m(\wh h) := \inf_{\phi \in \Phi_{a,b}^1}(\wh h, \phi), \quad M(\wh h) := \sup_{\phi \in \Phi_{a,b}^1}(\wh h, \phi).
\]
The random variables 
$m(\wh h)$ and $M(\wh h)$ are a.s. finite since $\wh h$ is a distribution and $\Phi^1_{a,b}$ is compact in the space of test functions.  
Let $f$ be some function supported on $[\wh x +a, \wh x +2a] \times [0,\pi]$ which is nonnegative, constant on vertical segments, and satisfies $\int f(z)\, d^2z = 1$. 

Since the statement of the lemma is translation invariant, we may translate $U$ so that $\lim_{z \to -\infty} |\varphi(z) - z| = 0$.
By \cite[Lemma 2.24]{ag-disk} we see that for some absolute constant $a$ we have $|\varphi(z) - z|, |\varphi'(z)|, |(\varphi^{-1})'| < a$ for all $z$ with $\Re z < \inf\{\Re z  :  z \in \cS \backslash \wh U \} -a$. Consequently, we can choose $b$ large in terms of $a, f$ so that $|\varphi'|^2f \circ \varphi \in \Phi^1_{a,b}$ for \emph{any} $U, \varphi$. Then
\begin{equation*}
\begin{gathered}
(\wh h^U, f) = (\wh h\circ \varphi^{-1} + Q \log |(\varphi^{-1})'|, f) = (\wh h, |\varphi'|^2f\circ \varphi) + (Q \log |(\varphi^{-1})'|, f),\\
m(\wh h) - Qa \leq (\wh h^U, f) \leq M(\wh h) + Qa. 
\end{gathered}
\end{equation*}
The event $\{ m(\wh h) - Qa > -K\} \cap \{ M(\wh h) + Qa < K\}$ holds with probability approaching 1 as $K \to \infty$. On this event, since $(\wh h^U, f) \in (-K,K)$ we can find some $x \in [\wh x + a, \wh x + 2a]$ such that the average of $h$ on $[x, x + i\pi]$ lies in $(-K,K)$. Moreover, notice that the curve $\varphi ([\wh x + 3a, \wh x + 3a + i\pi])$ is contained in $U$ and lies to the right of $[\wh x + 2a, \wh x + 2a + i\pi]$ (since $|\varphi(z) - z|< a$ for $z \in [\wh x + 3a, \wh x + 3a + i\pi]$). Therefore $x$ satisfies the first claim of the lemma.

Finally, we claim that with probability approaching 1 as $K \to \infty$ we have $\nu_{\wh h}((\wh x - K^2 - a, \infty) \cap \wh U) < \frac12 e^{\frac\gamma2 K^3}$; since $|\varphi(z) - z|< a$ for $z$ to the left of $\wh x$, this implies the last assertion of the lemma. Since $F$ has positive probability, it suffices to prove this claim in the setting where $(\cS, \wh h, +\infty, -\infty)$ is \emph{not} conditioned on $F$, so we will assume this. The left-to-right field average process $(\wh h_t)_{t \in \R}$ is Brownian motion started from $+\infty$ to $-\infty$ with variance 2 and downward drift. Define the stopping time $\wh \tau_{s} = \inf \{ t \in \R \: : \: \wh h_t = s\}$. Fix some large $M$, then $\wh h_{\wh x} < M$ with probability $1- o_M(1)$, hence $\wh \tau_M < \wh x$ with probability $1-o_M(1)$. By Brownian motion estimates we have $\wh \tau_{K^{5/2}}  < \wh \tau_M -  K^2 - a$ with probability tending to 1 as $K \to \infty$ for fixed $M$. Finally, since the law of $e^{-\frac\gamma2 K^{5/2}}\nu_{\wh h} (\tau_{K^{5/2}}, \infty)$ does not depend on $K$, we have $\nu_{\wh h}((\tau_{K^{5/2}}, \infty)) < \frac12	e^{\frac\gamma2 K^3}$ with probability $1-o_K(1)$. Sending $K \to \infty$ then $M \to \infty$ and combining these three estimates, we conclude that with probability approaching 1 we have $\nu_{\wh h}((\wh x - K^2 - a, \infty)) < \frac12 e^{\frac\gamma2 K^3}$, as desired. 
\end{proof}

Recall the event $\mathrm{BN}_{r, \zeta, K, \eps}= E_{r,K} \cap E'_{r, \zeta, \eps}$ of Proposition~\ref{prop-extrinsic}. 
Abusing notation, define  
\eqb\label{eq-G}
G_r := \left\{ \bigcup_{j=1}^n \cD_j^\bullet \subset (\cS_+ + \tau_{-r}) \right\};
\eqe
i.e. $G_r$ is the event that the white regions in Figure~\ref{fig-bottleneck} (top left) lie to the right of $\tau_{-r}$. More precisely, for $j = 1, \dots, n$ the curve segments of $\eta_j$ between $w_j$ and $+\infty$ lie in $\cS_+ + \tau_{-r}$ (with $\eta_n = \R + i\pi$).
\begin{proposition}\label{prop-E-given-F}
$\P[\mathrm{BN}_{r,\zeta,K,\eps} \cap G_{r} \mid F_{r, K, \eps}]\to 1$ as $r \to \infty, \zeta \to \infty, K \to \infty$ in that order. 
\end{proposition}
\begin{proof}
In this argument, we use ``with high probability'' as shorthand for  ``with probability approaching 1 as first $r \to \infty$ then $K \to \infty$''.

Let $U$ be the unbounded connected component of the set parametrizing $\cA$, let 
 $y = \inf \{ \Re z : z \in \cS \backslash U\}$, and define $(Y_t)_{t \geq 0}$ to be the field average of $h$ on $[y - t, y-t + i\pi]$. We claim that, since $F_{r, K, \eps}$ is measurable with respect to $(h|_{\cS_+ + y}, \eta_1, \dots, \eta_n)$, when we condition on $F_{r, K, \eps}$, the law of $(Y_t)_{t \geq 0}$ is Brownian motion started at (the random) $Y_0$ with variance 2 and upward drift. This claim follows from a minor modification to the proof of \cite[Lemma 2.10]{ag-disk}, and is essentially a Markov property of the field when we explore it from right to left, analogous to the domain Markov property of GFFs. We leave the details to the reader.

Transferring the high probability estimate Lemma~\ref{lem-comparison-wh-h} from $\wh \cA$ to $\cA + r +K$ using Lemma~\ref{lem-proxy}, we conclude that with high probability we can find a point $x \in \R$ so that $[x, x + i\pi]\subset U$, the average of $h$ on $[x, x+i\pi]$ is between $-r-2K$ and $-r$, and $\nu_h((x - K^2, \infty) \cap \partial U) < \frac12 e^{\frac\gamma2 (-r+K^3)}$. Restrict to the event that these occur and choose $x \leq y$ to be the rightmost point satisfying these constraints. Since the average of $h$ on $[x, x+i\pi]$ is less than $-r$, we see that $\tau_{-r} < x\leq y$; this proves $\P[G_r \mid F_{r, K, \eps}] \to 1$.

We claim that with high probability we have $\tau_{-r} > x - K^2$. Indeed, if $x < y$, then by the Markov property of Brownian motion $(Y_{t + (y-x)} - Y_{y-x})_{t\geq 0}$ is Brownian motion started at 0 with variance 2 and linear upward drift. Thus with high probability $Y_{t + (y-x)}- Y_{y-x} >  2K$ for all $t > K^2$; in particular the field average on any vertical segment to the left of $x - K^2$ is at least $-r$ with high probability. Consequently, $\tau_{-r} > x - K^2$. The case $x=y$ similarly yields $\tau_{-r} > x - K^2$. 

Finally, we conclude that with high probability 
\[ 
\nu_h([\tau_{-r}, \infty)) = \nu_h([\tau_{-r}, \infty)\cap \partial U ) + a_1 + 1 \leq \frac12 e^{\frac\gamma2 (-r + K^3)} + 4e^{\frac\gamma2 (-r-K)} + 1 < 1 + e^{\frac\gamma2(-r+K^3)},
\]
and similarly with high probability $\nu_h([\tau_{-r}, \infty)) > 1$. This shows that $\P[E_{r, K} \mid F_{r, K, \eps}] \to 1$. Combining with the last claim of Proposition~\ref{prop-bulk-law}, we conclude that $\P[\mathrm{BN}_{r, \zeta, K, \eps} \mid F_{r, K,\eps}] \to 1$ as $r, \zeta, K \to \infty$. 
\end{proof}

\subsection{Convergence to thick quantum disk}\label{subsec-intrinsic-thick}
In this section we prove Proposition~\ref{prop-weak-theorem}. Proposition~\ref{prop-bulk-law} shows that $(\cS, h, +\infty,-\infty, \eta_1, \dots, \eta_{n-1})$ conditioned on $F_{r,K,\eps}$ approximates a welding of thin quantum disks. 
From Sections~\ref{subsec-intrinsic-thin-decomp} and~\ref{subsec-intrinsic-compatibility}, this is close to the law of $(\cS, h, +\infty,-\infty, \eta_1, \dots, \eta_{n-1})$ conditioned on $F_{r, K, \eps} \cap G_r \cap \mathrm{BN}_{r, \zeta, K, \eps}$. 

When we condition only on $\mathrm{BN}_{r, \zeta, K, \eps}$, the field and curves in neighborhoods of $+\infty$ resemble those of a quantum disk decorated by macroscopic curves (namely with quantum lengths at least $\eps$), and are almost independent from the field and curves near the bottleneck (Proposition~\ref{prop-extrinsic}). On $\mathrm{BN}_{r, \zeta, K, \eps}$, the event $F_{r, K, \eps} \cap G_r$ is almost determined by the field and curves near the bottleneck, hence the field and curves in neighborhoods of $+\infty$ conditioned on $F_{r, K, \eps} \cap G_r \cap \mathrm{BN}_{r, \zeta, K, \eps}$ look like a quantum disk decorated by macroscopic curves. This concludes the proof of Proposition~\ref{prop-weak-theorem}. 

To that end, we make a general statement about the near-independence of SLE in spatially separated regions in Lemma~\ref{lem-curve-limit} (whose proof we defer to Appendix~\ref{appendix-SLE}), then carry out the argument sketched above. 

\begin{lemma}[Near independence of SLE]\label{lem-curve-limit}
Suppose $(\eta_1, \dots, \eta_{n-1})\sim \cP^\disk(W_1, \dots, W_n)$, and condition on $(\eta_1 \cap \cS_+, \dots, \eta_{n-1} \cap \cS_+)$. Then $(\eta_1 \cap \cS_+, \dots, \eta_{n-1} \cap \cS_+)$-almost surely, as $N \to \infty$ the total variation distance between the conditional law of $(\eta_1 \cap (\cS_- - N), \dots, \eta_{n-1} \cap (\cS_- - N))$ and the unconditioned law of $(\eta_1 \cap (\cS_- - N), \dots, \eta_{n-1} \cap (\cS_- - N))$ goes to zero. 
\end{lemma}

Recall from Section~\ref{subsec-thick-disint} that a quantum surface $(\cS, h, +\infty,-\infty, \eta_1, \dots, \eta_n)$ satisfying $\nu_h(\R) > \frac12$ is \emph{canonically embedded} if $\nu_h(\R_+) = \frac12$. 

\begin{proposition}\label{prop-F-given-E}
Consider the canonically embedded curve-decorated surface $(\cS, h, +\infty, -\infty, \eta_1, \dots, \eta_{n-1})$ conditioned on $F_{r, K, \eps} \cap G_{r} \cap\mathrm{BN}_{r, \zeta, K, \eps}$, where $G_r$ is defined as in~\eqref{eq-G}. Send $r \to \infty, \zeta \to \infty, K \to \infty$ in that order. Then in any neighborhood $U$ of $+\infty$ with $-\infty \not \in U$, the law of $(h|_U, \eta_1 \cap U, \dots, \eta_{n-1} \cap U)$  converges in total variation to the law of $(\wh\psi|_U, \wh \eta_1 \cap U, \dots, \wh \eta_{n-1}\cap U)$, where 
$(\cS, +\infty,-\infty, \wh \psi, \wh \eta_1, \dots, \wh \eta_{n-1})$ is taken from $\cM_2^\disk(W, 1) \otimes \cP^\disk(W_1, \dots, W_n)$ (with canonical embedding) and conditioned on $\{ \nu_{\wh\psi}(\wh\eta_j) > \eps \text{ for } j = 1, \dots, n\}$ (with $\wh \eta_n := \R + i\pi$).  
\end{proposition}
\begin{proof}
In this proof we will send parameters $r \to \infty, \zeta \to \infty, S\to\infty, K \to \infty$ in that order. We will write $o_\zeta(1)$ (resp. $o_S(1)$) for a quantity that goes to zero as $r,\zeta\to\infty$ (resp. $r, \zeta, S \to \infty$) in that order.

First sample $(\cS, h, +\infty,-\infty, \eta_1, \dots, \eta_{n-1})$ conditioned only on $\mathrm{BN}_{r, \zeta, K, \eps}$. By Proposition~\ref{prop-extrinsic} and Lemma~\ref{lem-curve-limit} we know that the joint law of $h(\cdot +\tau_{-r})|_R + r, ((\eta_1 -\tau_{-r})\cap R, \dots, (\eta_{n-1} - \tau_{-r}) \cap R), e^{\frac\gamma2r}(\nu_h([\tau_{-r}, \infty))-1)$ and $(h|_U, \eta_1\cap U, \dots, \eta_{n-1}\cap U)$ converges in total variation as $r, \zeta \to \infty$ to the independent objects $\phi, (\eta_1'\cap R, \dots, \eta_{n-1}' \cap R), V$ and $(\wh\psi|_U, \wh \eta_1 \cap U, \dots, \wh \eta_{n-1}\cap U)$. Here, $\phi$ is the field described in~\eqref{eq-phi}, $(\eta_1', \dots, \eta_{n-1}')\sim\cP^\disk(W_1, \dots, W_n)$, and $V \sim \mathrm{Unif}([0, e^{\frac\gamma2K}])$.

We claim that outside a bad event of probability $o_S(1)$, the event $F_{r, K, \eps} \cap G_r$ is measurable with respect to the tuple $h(\cdot +\tau_{-r})|_R + r, ((\eta_1 -\tau_{-r})\cap R, \dots, (\eta_{n-1} - \tau_{-r}) \cap R), e^{\frac\gamma2r}(\nu_h([\tau_{-r}, \infty))-1)$, and moreover $F_{r, K, \eps} \cap G_r$ occurs with uniformly positive probability as $r, \zeta, S \to \infty$. Assuming this claim, by the previous discussion, when we further condition on $F_{r, K, \eps} \cap G_r$ the law of $(h|_U, \eta_1\cap U, \dots, \eta_{n-1}\cap U)$ only changes by $o_S(1)$ in total variation, and hence is $o_S(1)$-close to $(\wh\psi|_U, \wh \eta_1 \cap U, \dots, \wh \eta_{n-1}\cap U)$. Taking $r,\zeta, S, K\to \infty$ then yields the proposition. 

We turn to the proof of the claim. 
For notational simplicity we work with the embedding of $(\cS, h, +\infty,-\infty, \eta_1, \dots, \eta_{n-1})$ where $\tau_{-r} = 0$ (rather than the canonical embedding). In this setting, the claim is that as we send $r,\zeta,S \to \infty$, outside of a bad event of probability $o_S(1)$, the event $F_{r, K, \eps}\cap G_r$ is determined by $h|_R + r, (\eta_1 \cap R, \dots, \eta_{n-1} \cap R)$ and $e^{\frac\gamma2r}(\nu_h(\R_+)-1)$ and occurs with uniformly positive probability. This follows from the following observations: 
\begin{itemize}
\item The condition $\{ \ell_j > \eps \text{ for all } j = 1,\dots, n\}$ of $F_{r, C, \eps}$ is always fulfilled because we are conditioning on $\mathrm{BN}_{r, \zeta, K, \eps}$. 

\item Recall $z_0 \in \R$ satisfies $\nu_h([z_0, \infty)) = 1$, i.e. $\nu_{h + r}([0, z_0]) = e^{\frac\gamma2 r}(\nu_h(\R_+) - 1)< e^{\frac\gamma2 K}$. 

Since the field average process of $\phi$ has upward drift, and the law of $h|_R + r$ is $o_r(1)$-close in total variation to that of $\phi$, with probability $1-o_S(1)$ we have $\nu_{h+r}([0,S]) > e^{\frac\gamma2 K}$, and on this event, $z_0$ is measurable with respect to $\sigma(h|_R + r, e^{\frac\gamma2r}(\nu_h(\R_+)-1))$. Likewise, with probability $1-o_S(1)$ the curves $(\eta_1\cap R, \dots, \eta_{n-1}\cap R)$ and point $z_0\in R$ are such that they determine $1_{G_r}$, and restricted to $G_r$ they determine $(a_j, b_j, c_j)_{j=1}^n$ (indeed, the event $G_r$ guarantees that the relevant curve segments are to the right of $[0,i\pi]$, and as $S \to \infty$ the right border of $R$ goes to $\infty$ so the curve segments lie in $R$ with high probability).  On this event the lengths $(e^{\frac\gamma2(r+K)}a_j, e^{\frac\gamma2(r+K)}b_j, e^{\frac\gamma2(r+K)}c_j)_{j=1}^n$ (and hence the event $F_{r, K, \eps} \cap G_r$) are measurable with respect to $\sigma(h|_R + r, \eta_1\cap R, \dots, \eta_{n-1} \cap R, e^{\frac\gamma2r}(\nu_h(\R_+)-1))$.

\item $F_{r, K, \eps} \cap G_r$ occurs with uniformly positive probability as $r, \zeta \to \infty$ then $S \to \infty$ because the tuple $(h|_R + r, \eta_1 \cap R, \dots, \eta_{n-1} \cap R)$ and $e^{\frac\gamma2r}(\nu_h(\R_+)-1)$ converge to some limit law as $r, \zeta \to \infty$. 
\end{itemize}
This proves the claim and hence the proposition.
\end{proof}

Now we are ready to prove Proposition~\ref{prop-weak-theorem}.
\begin{proof}[Proof of Proposition~\ref{prop-weak-theorem}]
In this proof, when we say ``with probability approaching 1'' or ``close in distribution'', we mean in the $r \to \infty, \zeta\to \infty, K \to \infty$ limit, and when we say ``in the bulk'', we mean in neighborhoods of $+\infty$ in the canonical embedding. 

Proposition~\ref{prop-E-given-F} tells us that $\P[G_r \cap \mathrm{BN}_{r,\zeta,K,\eps} \mid F_{r, C, \eps}] \to 1$. Therefore Proposition~\ref{prop-bulk-law} tells us that the law of $(h, \eta_1, \dots, \eta_n)$ in the bulk conditioned on $F_{r, C, \eps}\cap G_r \cap \mathrm{BN}_{r,\zeta,K,\eps}$ is close in distribution to that of $\frac1Z\iiint_\eps^\infty \cM_2^\disk(W_1; 1, \ell_1)\cdots \cM_2^\disk(W_n; \ell_{n-1}, \ell_n) \, d\ell_1 \dots d\ell_n$. But Proposition~\ref{prop-F-given-E} tells us that the law of $(h, \eta_1, \dots, \eta_n)$ in the bulk conditioned on $F_{r, C, \eps}\cap G_r\cap\mathrm{BN}_{r,\zeta,K,\eps}$ is close in distribution to that of $\cM_2^\disk(W;1)\otimes \cP^\disk(W_1, \dots, W_n)$ conditioned on the event $A_\eps$ that boundary arc and interface lengths are greater than $\eps$. Thus for fixed $\eps > 0$ and for some constant $c_\eps>0$, we have
\[(\cM_2^\disk(W;1)\otimes \cP^\disk(W_1, \dots, W_n))|_{A_\eps} = c_\eps \iiint_\eps^\infty \prod_{i=1}^n \cM_2^\disk(W_i; \ell_{i-1}, \ell_i) \,d \ell_1 \dots d \ell_n. \]
For any $\eps' > \eps$, by restricting first to $A_\eps$ then to $A_{\eps'}$, we see that $c_\eps = c_{\eps'}$ so the constant does not depend on $\eps$. Sending $\eps \to 0$ yields Proposition~\ref{prop-weak-theorem}.
\end{proof}

\section{Conclusion of the proofs of main results}\label{sec-general}
In this section we extend Proposition~\ref{prop-weak-theorem} to Theorem~\ref{thm-disk-cutting}, and explain the argument modifications needed to obtain Theorem~\ref{thm-sphere-cutting}.

\begin{proof}[Proof of Theorem~\ref{thm-disk-cutting} when $W_1, \dots, W_n, W \neq \frac{\gamma^2}2$]
We discuss the proof of Theorem~\ref{thm-disk-cutting} in three different regimes. 
\medskip

\noindent \textbf{Case 1: $W_1, \dots, W_n \in (0, \frac{\gamma^2}2)$ and $W > \frac{\gamma^2}2$.}

Proposition~\ref{prop-weak-theorem} tells us that 
\alb
&\cM_2^\disk(W;1)\otimes \cP^\disk(W_1, \dots, W_n) \\
&= c_{W_1, \dots, W_n} \iiint_0^\infty \cM_2^\disk(W_1; 1, \wt\ell_1) \times \cM_2^\disk(W_2; \wt\ell_1, \wt\ell_2) \times \cdots \times \cM_2^\disk(W_n; \wt\ell_{n-1}, \wt\ell_n) \,d \wt\ell_1 \dots d \wt\ell_n.
\ale
Add $\frac2\gamma \log \ell$ to the field and apply Lemma~\ref{lem-scaling-disint} $(n+1)$ times. Writing $\ell_j = \ell \wt \ell_j$ for $j = 1, \dots, n-1$ and $\ell' = \ell \wt \ell_n$, then disintegrating with respect to $\ell'$, we have for a.e. $\ell'>0$ that  
\alb
&\cM_2^\disk(W;\ell,\ell')\otimes \cP^\disk(W_1, \dots, W_n) \\
&= c_{W_1, \dots, W_n} \iiint_0^\infty \cM_2^\disk(W_1; \ell, \ell_1) \times \cM_2^\disk(W_2; \ell_1, \ell_2) \times \cdots \times \cM_2^\disk(W_n; \ell_{n-1}, \ell') \,d \ell_1 \dots d \ell_{n-1}.
\ale
Continuity of $\cM_2^\disk(W; \ell, \ell')$ and $\cM_2^\disk(W_n; \ell_{n-1}, \ell')$ in $\ell'$ (see the proof of Proposition~\ref{prop-bulk-law} for the argument by continuity) then yields the result for all $\ell'>0$, establishing Case 1.
\medskip

\noindent \textbf{Case 2: $W_1, \dots, W_n \in (0, \frac{\gamma^2}2)$ and $W < \frac{\gamma^2}2$.}

Choose $W_{n+1} \in (0, \frac{\gamma^2}2)$ so that $\sum_1^{n+1} W_i = W + W_{n+1} > \frac{\gamma^2}2$. By the definition of $\cP^\disk$, one can sample $(\eta_1, \dots, \eta_{n-1}, \eta_n)\sim \cP^\disk(W_1, \dots, W_n, W_{n+1})$ by first sampling $\eta \sim \cP^\disk(W_1 + \dots + W_n, W_{n+1})$, then independently sampling $n-1$ curves in each bounded connected component $D \subset \cS \backslash \eta$ from $\cP_D(W_1, \dots, W_n)$, and concatenating to get $(\eta_1, \dots, \eta_{n-1})$. Therefore, applying Case 1 to the $(n+1)$-tuple $(W_1, \dots, W_{n+1})$ and to the pair $(W, W_{n+1})$ yield
\alb
&\cM_2^\disk(W;\ell,\ell'')\otimes \cP^\disk(W_1, \dots, W_{n+1}) \\
&= c_{W_1, \dots, W_{n+1}} \iiint_0^\infty \cM_2^\disk(W_1; \ell, \ell_1) \times \cdots \times \cM_2^\disk(W_n; \ell_{n-1}, \ell') \times \cM_2^\disk (W_{n+1}; \ell', \ell'')\, d \ell_1 \dots d \ell_{n-1}\, d \ell' \\
&= c_{W_1 + \dots + W_n, W_{n+1}} \int_0^\infty \left(\cM_2^\disk(W_1+\dots + W_n; \ell, \ell_n)\otimes \cP^\disk(W_1, \dots, W_n) \right) \times \cM_2^\disk (W_{n+1}; \ell', \ell'') \, d \ell'.
\ale
Disintegrating with respect to $\ell'$ yields the desired identity for a.e.\ $\ell' > 0$, and continuity extends this to all $\ell'> 0$. Thus we have shown Theorem~\ref{thm-disk-cutting} for Case 2. 
\medskip

\noindent \textbf{Case 3: $W_1, \dots, W_n \in (0, +\infty) \backslash \{ \frac{\gamma^2}2\}$ and $W > \frac{\gamma^2}2$.}

Choose some sufficiently large $N$ and decompose $W_i = W_i^1 + \cdots + W_i^N$ with $W_i^j \in (0, \frac{\gamma^2}2)$ for $1\leq i \leq n$ and $1 \leq j \leq N$, then by Case 1 we have for constants $c_1, \dots, c_n \in (0, \infty)$
\[\cM_2^\disk(W_i; \ell, \ell') = c_i\iiint_0^\infty \cM_2^\disk(W_i^1; \ell, \ell_1) \times\cdots \times \cM_2^\disk(W_i^N; \ell_{N}, \ell')\, d \ell_1 \dots d \ell_{N-1}. \]
Here, the right hand side is a measure on curve-decorated surfaces; forgetting the curves yields the left hand side. 
Applying Case 1 to the $nN$ weights $((W_1^j)_{j=1}^N, \dots, (W_n^j)_{j=1}^N)$ and comparing to the above yields
\[\cM_2^\disk(W;\ell,\ell')\otimes \cP 
= c_{W_1, \dots, W_n} \iiint_0^\infty \cM_2^\disk(W_1; \ell, \ell_1)\times \cdots \times \cM_2^\disk(W_n; \ell_{n-1}, \ell')\, d \ell_1 \dots d \ell_{n-1}\]
where $\cP$ is a probability measure on $(n-1)$-tuples of curves obtaining by forgetting some curves of $\cP^\disk((W_1^j)_{j=1}^N, \dots, (W_n^j)_{j=1}^N)$. The same argument using Theorem~\ref{thm-wedges} yields $\cM^\wed(W) \otimes \cP = \prod_{i=1}^n \cM^\wed(W_i)$, and comparing this with Theorem~\ref{thm-wedges} for weights $(W_1, \dots, W_n)$ yields $\cP = \cP^\disk(W_1, \dots, W_n)$. Thus we have shown Theorem~\ref{thm-disk-cutting} for Case 3. 
\medskip

\noindent \textbf{Case 4: $n = 2$, $W_1 = \frac{\gamma^2}2$ and $W_2=2$.}
This follows from applying Case 3 and sending $\eps \to 0$ when $n = 2$, $W_1 = \frac{\gamma^2}2+\eps$ and $W_2 = 2$. See Appendix~\ref{appendix-A} for details. 
\medskip

\noindent \textbf{Case 5: Either $W = \frac{\gamma^2}2$ or some $W_j = \frac{\gamma^2}2$.}

The case $W = \frac{\gamma^2}2$ follows from Case 4 and the argument of Case 2. The case where $W_j = \frac{\gamma^2}2$ for some $j$ then follows from the argument of Case 3.

\end{proof}

The proof of Theorem~\ref{thm-sphere-cutting} is nearly identical to that of Theorem~\ref{thm-disk-cutting}, so we explain it briefly. First we will show the analog of Proposition~\ref{prop-weak-theorem}, and then extend it to the full result using scaling arguments and Theorem~\ref{thm-disk-cutting}. 

Let $\cC = \R \times [0,2\pi] /\mathord\sim$ be the cylinder (here $\R$ and $\R + 2 \pi i$ are identified under the relation $x \sim x + 2 \pi i$).  

\begin{figure}[ht!]
\begin{center}
  \includegraphics[scale=0.75]{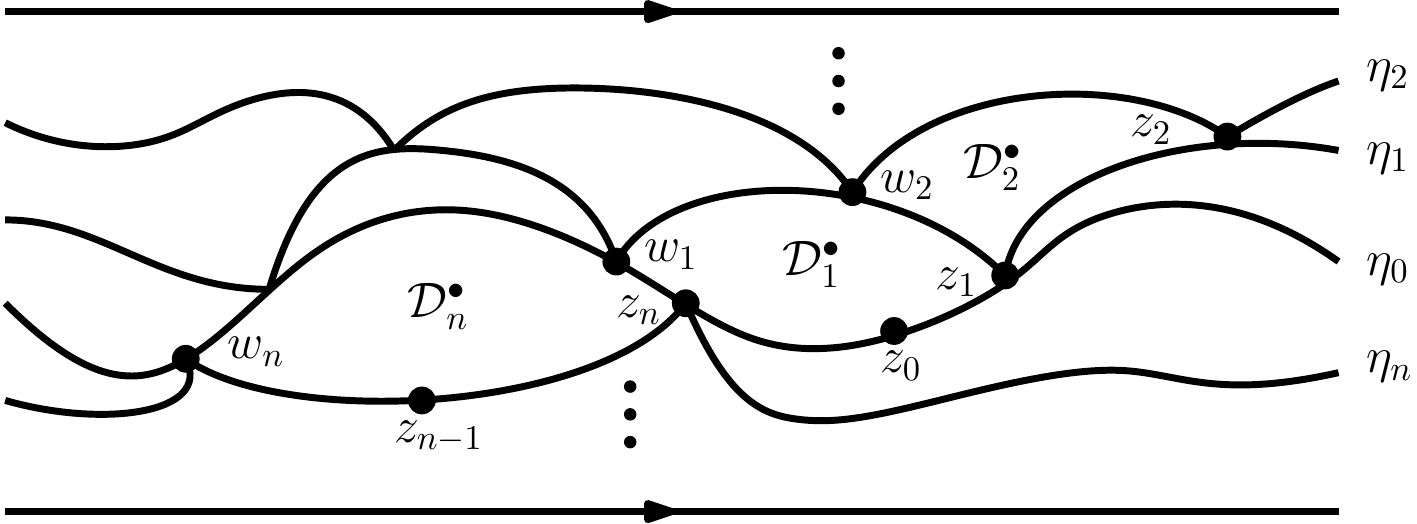}%
 \end{center}
\caption{\label{bottleneck-sphere} Consider a quantum cone $(\cC, h, +\infty, -\infty)$ decorated by the curves $(\eta_0, \dots, \eta_{n-1}) \sim$ $\cP^\sph(W_1, \dots, W_n)$. Let $z_0$ be the  point on $\eta_0$ with $\nu_h$-length from $+\infty$ equal to 1. Iteratively for $j = 1, \dots, n$, let $\cD_j^\bullet$ be the component of $\cW_j$ with $z_{j-1}$ on its boundary, and let $w_j, z_j$ be its left and right marked points respectively. Let $F_{r, C, \eps}$ be the event that $z_n$ lies on $\partial \cD_1^\bullet$ between $w_1$ and $z_0$, and that the length bounds~\eqref{eq-F-sph} hold.}
\end{figure}

\begin{lemma}\label{lem-weak-thm-sph}
Fix $n \geq 2$ and fix $W_1, \dots, W_n \in (0, \frac{\gamma^2}2)$.
Consider a field-curves pair $(\psi, \eta_0, \dots, \eta_{n-1})$ taken from $\cM_2^\sph(W; 1) \otimes \cP^\sph(W_1, \dots, W_n)$. When we condition on $\nu_\psi(\eta_0) = 1$ and cut along $\eta_0, \dots, \eta_{n-1}$, we obtain an $n$-tuple of quantum surfaces with law 
\[
\wh c_{W_1, \dots, W_n} \iiint_0^\infty \cM_2^\disk(W_1; 1, \ell_1) \times \cdots \times \cM_2^\disk(W_n; \ell_{n-1}, 1) \,d \ell_1 \dots d\ell_{n-1}
\]
for some $\wh c_{W_1, \dots, W_n} \in (0, \infty)$. 
\end{lemma}
\begin{proof}
Consider a quantum cone $(\cC, h, +\infty, -\infty)\sim \cM^\mathrm{cone}(W)$ decorated by independent curves $(\eta_0, \eta_1, \dots, \eta_{n-1}) \sim \cP^\sph(W_1, \dots, W_n)$; the curves cut the cone into $n$ independent quantum wedges $\cW_1, \dots, \cW_n$ with weights $W_1, \dots, W_n$ \cite[Theorems 1.2, 1.5]{wedges}. We define the event $F_{r, K, \eps}$ as follows: Let $z_0$ be the point on $\eta_0$ at distance $1$ from $+\infty$, and define the quantum surface $\cD^\bullet_1$ to be the bead of $\cW_1$ with $z_0$ on its boundary; call its left and right marked points $w_1$ and $z_1$ respectively. Iteratively define $\cD_2^\bullet, z_2, w_2, \dots, \cD_n^\bullet, z_n, w_n$ similarly, and define $(a_j, b_j, c_j, \ell_j)_{j=1}^n$ in the same way as in the disk case. Let $F_{r, K, \eps}$ be the event that the following holds, see Figure~\ref{bottleneck-sphere}:
\begin{itemize}
\item For each $j = 1, \dots, n$ we have the inequalities
\eqb \label{eq-F-sph}
e^{\frac\gamma2(r+K)} b_j \in (1, 2), \quad e^{\frac\gamma2(r+K)}(a_j+b_j) \in (3, 4), \quad  e^{\frac\gamma2(r+K)}c_j \in (5, 6), \quad \ell_j > \eps ;
\eqe

\item The point $z_{n}$ lies on the bead $\cD_1^\bullet$, between the points $w_1$ and $z_0$. 
\end{itemize}
The exact choice of the second condition is not too important; any suitable variant of ``the cycle of beads $\cD_1^\bullet, \dots, \cD_n^\bullet$ closes up'' suffices.  This condition is used to prove the analog of Lemma~\ref{lem-length-distribution}.

Conditioning on $F_{r,C,\eps}$ gives a decomposition of the quantum cone into the quantum surfaces $\cA, \cB$, and $(\cD_j^\bullet)_{j=1}^n$, where $\cA$ is infinite and $\cB$ is finite.  
As before, when we condition on the lengths $(a_j, b_j, c_j, \ell_j)_{j=1}^n$, these quantum surfaces become mutually independent.  By following the steps in the proof of Proposition~\ref{prop-weak-theorem}, we obtain Lemma~\ref{lem-weak-thm-sph}. \end{proof}

\begin{proof}[Proof of Theorem~\ref{thm-sphere-cutting}]
The theorem in the case $n \geq 2$ and $W_1, \dots, W_n \in (0, \frac{\gamma^2}2)$ follows from Lemma~\ref{lem-weak-thm-sph} and a scaling argument using Lemma~\ref{lem-scaling-disint} (see Case 1 in the proof of Theorem~\ref{thm-disk-cutting}). For the case $n = 1$ and $W_1 \in (0, \frac{\gamma^2}2)$, we choose any $W_1', W_2' \in (0, \frac{\gamma^2}2)$ with $W_1' + W_2' = W_1$ and apply the $n \geq 2$ case and Theorem~\ref{thm-disk-cutting}. Finally, for the case where $n \geq 1$ and $W_1, \dots, W_n>0$ are arbitrary, we split each thick quantum disk into thin quantum disks as in the proof of Case 3 of Theorem~\ref{thm-disk-cutting}. 
\end{proof}

\section{Application to finite-area mating of trees}\label{sec-MOT}

We now present two applications of our main results. In Section~\ref{subsec-MOT-rederive} we explain a unified derivation of the mating-of-trees theorems for the quantum sphere and disk, building on the mating-of-trees theorem for the $2-\frac{\gamma^2}{2}$ quantum wedge from \cite{wedges}. Since these results are not new, we only provide proof sketches but the proofs can be made complete without substantial difficulty by filling in more details. In Section~\ref{subsec-MOT}, we show a new mating-of-trees theorem for $\cM_2^\disk(\frac{\gamma^2}2)$, which is crucial for several subsequent works \cite{AHS-SLE-integrability, ARS-FZZ}.

\subsection{Mating of trees for weight 2 quantum disk and weight $4-\gamma^2$ quantum sphere }\label{subsec-MOT-rederive}

In this section we explain how our arguments and main theorems yield a systematic treatment of the quantum sphere and disk mating-of-trees theorems. The quantum sphere theorem was originally proved in \cite[Theorem 1.1]{sphere-constructions}, and the quantum disk theorem was shown in \cite{wedges}\footnote{See \cite[Theorem 2.1]{sphere-constructions} and the paragraph before it for discussion on the proof of the $\gamma \in [\sqrt2,2)$ quantum disk mating-of-trees theorem in \cite{wedges}.} for $\gamma \in (\sqrt2, 2)$ and \cite[Theorem 1.1]{ag-disk} for $\gamma \in (0,\sqrt2]$. As these results are already present in the literature, we only sketch the proofs --- for example, we will rely on several facts concerning Brownian motion without justification. 
These proofs demonstrate the robustness of our arguments and the effectiveness of our main results.

We start by recalling the setup for the mating-of-trees framework; see \cite{wedges,ghs-mating-survey} for more details.
Fix $\kappa' = \frac{16}{\gamma^2}$. We can define  \emph{space-filling $\SLE_{\kappa'}$} curves between two marked boundary points of a simply connected domain; see \cite[Section 1.2.3]{ig4} and~\cite[Section 3.6.3]{ghs-mating-survey}. Suppose $\eta'$ is a space-filling $\SLE_{\kappa'}$
drawn on an independent $\gamma$-LQG quantum surface. We parametrize it so it covers a unit of \emph{quantum area} per unit of time. 
Moreover, at each instant, the boundary of the region $\eta'$ has explored is locally absolutely continuous with respect to an $\SLE_\kappa$ curve, and one can use this to argue that the boundary a.s.\ has a well defined quantum length  \cite{shef-zipper}.

\begin{figure}[ht!]
	\begin{center}		\includegraphics[scale=0.7]{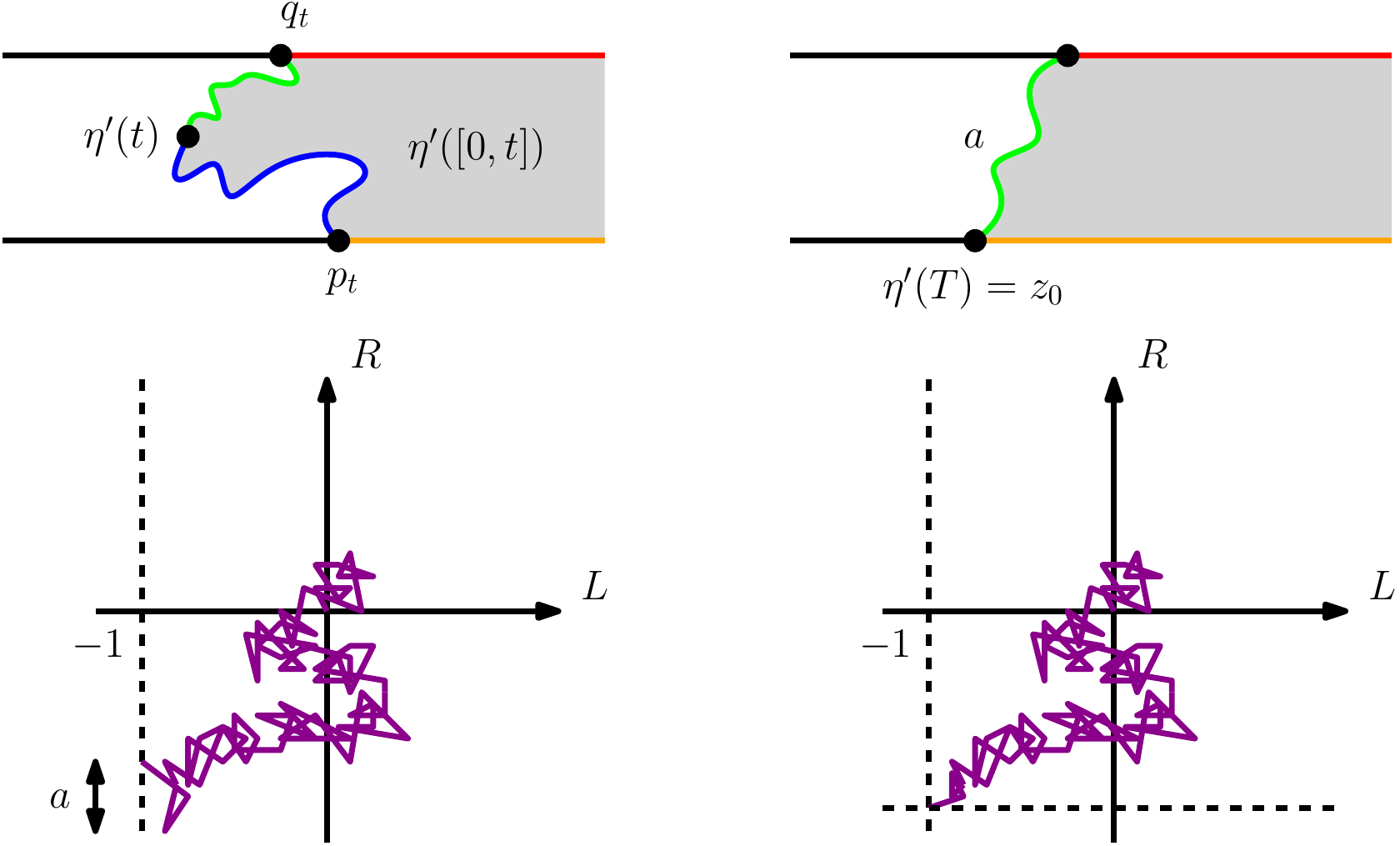}%
	\end{center}
	\caption{ \label{fig-MOT} \textbf{Top left.} Definition of $(L_t,R_t)$  in Theorem~\ref{thm-MOT}: $L_t$ is the quantum length of the boundary arc of $\eta'([0,t])$ from $\eta'(t)$ to $p_t$ minus $\nu_h([p_t, \infty))$ (blue minus orange), and $R_t$ is is the quantum length of the boundary arc of $\eta'([0,t])$ from $\eta'(t)$ to $q_t$ minus $\nu_h([q_t, \infty))$ (green minus red). \textbf{Top right.} Definition of $F_{r,K}$ in the proof of Proposition \ref{prop-pinch-MOT}. \textbf{Bottom left.} The Brownian motion $(L_t, R_t)_{[0,A]}$ conditioned on $F_{r,K}$. \textbf{Bottom right.} Illustration of $(L_t, R_t)_{[0,A]}$, where $( L_t - L_A, R_t - R_A)_{[0,A]} \sim C \int_0^\infty \mu_{\R_+^2}^\gamma ((1,\ell),0) d \ell$.  
}
\end{figure}

We can now state the mating-of-trees theorem for the weight $2-\frac{\gamma^2}2$ quantum wedge. 
For $\gamma \in (0, \sqrt2]$, 
consider a weight $2-\frac{\gamma^2}2$ quantum wedge $(\cS, h, +\infty, -\infty)$, 
decorated with an independent space-filling curve $\eta' $  from $+\infty$ to $-\infty$  parametrized by the quantum area. For each $t>0$ let $p_t\in \R$ and $q_t\in \R + i\pi$ be the leftmost points such that $[p_t, \infty) , [q_t, \infty) \subset \eta'([0,t])$. The boundary $\partial (\eta'([0,t]))$ has a well defined quantum length.  Let $L_t$ be the quantum length of the boundary arc of $\eta'([0,t])$ from $\eta'(t)$ to $p_t$ minus $\nu_h([p_t, \infty))$, and let $R_t$ be  the quantum length of the boundary arc of $\eta'([0,t])$ from $\eta'(t)$ to $q_t$ minus $\nu_h([q_t, \infty))$. See  Figure~\ref{fig-MOT} for an illustration. For $\gamma \in (\sqrt2, 2)$, the weight $2-\frac{\gamma^2}2$ quantum wedge is thin, but the same definition applies with $p_t, q_t$ being the furthest points on the left and right boundaries from the root for which the space-filling curve has filled the boundary arcs from $p_t$ and $q_t$, respectively, to the root. 
\begin{theorem}[{\cite[Theorem 1.9]{wedges}}]\label{thm-MOT} For some $\gamma$-dependent constant $\mathbbm{a}>0$, the process
$(L_t, R_t)_{t \geq 0}$ evolves as Brownian motion with covariances given by 
\eqb\label{eq-cov}
\Var(L_t) = \Var(R_t) =  \mathbbm{a}^2 t, \quad \Cov(L_t, R_t) = -\cos(\pi \gamma^2/4)\mathbbm{a}^2t \quad \text{ for } t \geq 0.
\eqe
\end{theorem}

\begin{remark}
	The variance $\BB a^2$  was not known until the work of the first and the third author with Remy~\cite{ARS-FZZ}, 
	which proves  that $\BB a^2= 2/\sin (\frac{\pi \gamma^2}4)$. This formula is not needed for our paper.
\end{remark}

We will use Theorem~\ref{thm-MOT} and our main results to rederive the mating-of-tree theorems for the quantum sphere and disk.

To that end we need some finite duration variants of the  Brownian motion $(L_t, R_t)_{t \geq 0}$  in Theorem~\ref{thm-MOT}. 
They are constructed through limiting procedures in the same spirit as \cite[Section 3]{lawler-werner-soup}, except that we consider   Brownian path measures  in a cone $\R_+^2 := (0, \infty)^2$ with an endpoint at the vertex, and we use correlated Brownian motion. 
We omit the detailed justification for  the existence of the various limits  because the arguments in \cite{lawler-werner-soup} still apply.

Let $\cK$ be the collection of all continuous planar curves of the form $\eta:[0,t_\eta] \to \R^2$. 
For $\eta\in \cK$ we call $t_\eta$ the duration of $\eta$.
Endow $\cK$ with the metric $d_\cK(\eta_1,\eta_2)= \inf_\theta \{ \sup_{s \in [0, t_{\eta_1}]} |s-\theta(s)| + |\eta_1(s) - \eta_2(\theta(s))|\}$ with the infimum taken over increasing homeomorphisms $\theta:[0, t_{\eta_1}] \to [0,t_{\eta_2}]$.
For $z \in \R_+^2$ let $\mu_{\R_+^2}^\gamma(z)$ denote the probability measure on $\cK$ corresponding to a Brownian motion $Z_t=L_t + R_t i$ started at $z = L_0 + R_0 i$ with covariance~\eqref{eq-cov}, 
and stopped upon hitting $\partial \R_+^2$. For $q \in i \R_+$, let $E_{q, \eps} := \{Z\text{\,\,exits }\R_+^2\text{ in } (q, q + \eps i) \}$, 
where $(q, q + \eps i)$ means the linear segment on $\R^2$ between $q$ and $q + \eps i$.
 Define the weak limit $\mu^\gamma_{\R_+^2}(z,q) = \lim_{\eps \to 0} \frac1\eps \mu^\gamma_{\R_+^2}(z)|_{E_{q,\eps}}$. This is a finite measure on $\cK$, supported on the set of paths from $z$ to $q$.  Similarly, for $p \in \R_+$ we can set $\mu^\gamma_{\R_+^2}(p,q) = \lim_{\eps \to 0} \frac1\eps \mu^\gamma_{\R_+^2}(p + i\eps, q)$. When one of the endpoints of the Brownian excursion is the origin, we need to normalize differently. Let $E_\eps:= \{Z\text{ exits }\R_+^2\text{ in } (0, \eps i)\}$ and define for $z \in \R_+^2$ the finite measure $\mu^\gamma_{\R_+^2}(z,0) = \lim_{\eps \to 0} \eps^{-\frac4{\gamma^2}} \mu^\gamma_{\R_+^2}(z)|_{E_\eps}$; Lemma~\ref{lem-corr-BM-rate} shows that $-\frac4{\gamma^2}$ is the correct exponent. We also define for $p \in \R_+$ the finite measure $\mu^\gamma_{\R_+^2}(p,0) = \lim_{\eps \to 0} \frac1\eps \mu^\gamma_{\R_+^2}(p+i\eps, 0)$.  Finally define the Brownian bubble measure $\mu^\gamma_{\R_+^2}(0,0) = \lim_{\eps \to 0} \eps^{-\frac4{\gamma^2} -1} \mu^\gamma_{\R_+^2}(\eps, 0)$.

We note that, as an immediate consequence of the scale invariance of Brownian motion and the exponents in the above definitions, for any $\lambda > 0$, $z \in \R_+^2$, and $x > 0$ we have
\eqb\label{eq-scale-BM}
\mu^\gamma_{\R_+^2}(\lambda z, 0) = \lambda^{-\frac4{\gamma^2}} (T_\lambda)_* \mu^\gamma_{\R_+^2}(z, 0)\quad \textrm{and}\quad \mu^\gamma_{\R_+^2}(\lambda x, 0) = \lambda^{-\frac4{\gamma^2}-1} (T_\lambda)_* \mu^\gamma_{\R_+^2}(x, 0),
\eqe
where $T_\lambda: \cK \rta \cK$  is the rescaling operator given by $T_\lambda(\eta):= \lambda \eta(\lambda^{-2} \,\cdot \,)$.

We first prove the variant of Theorem~\ref{thm-MOT} with the weight $2-\frac{\gamma^2}2$ quantum wedge replaced by a   weight $2-\frac{\gamma^2}2$ quantum disk.
As in Proposition~\ref{prop-weak-theorem}, we first restrict to the case when one boundary arc of the disk has quantum length  1.

\begin{proposition}\label{prop-pinch-MOT}
Suppose we are in the setting of Theorem~\ref{thm-MOT} but with the quantum wedge replaced by a quantum surface $(\cS,\wt h,+\infty,-\infty)$ sampled from $\cM_2^\disk(2 - \frac{\gamma^2}2; 1)$.
Let $(L_t, R_t)_{[0,A]}$ be the boundary length process, where $A$ is the random quantum area of the quantum disk.
Then for some constant $C >0$, the law of $( L_t - L_A, R_t - R_A)_{[0,A]}$ is given by $C \int_0^\infty \mu_{\R_+^2}^\gamma (1 + \ell i,0)\, d \ell$. 
\end{proposition}
\begin{proof}[Sketch of proof]
We focus on the $\gamma \in (0, \sqrt2)$ case first, so $2-\frac{\gamma^2}2 > \frac{\gamma^2}2$.
Let $(\cS, h, +\infty, -\infty, \eta')$ be a weight $2-\frac{\gamma^2}2$ quantum wedge decorated by a space-filling SLE$_{\kappa'}$ between the two marked points. We define the field bottleneck event $\mathrm{BN}_{r, \zeta, K} = E_{r,K} \cap E'_{r,\zeta}$ as in Proposition~\ref{prop-extrinsic}, but set instead $E'_{r,\zeta}:= \{ \sigma_{-\zeta} < \infty\}$. That is, remove the curve length condition.  

We now define a curve bottleneck event, see Figure~\ref{fig-MOT} (top right).
Decorate the quantum wedge by an independent space-filling $\SLE_{\kappa'}$ curve $\eta'$ from $+\infty$ to $-\infty$ parametrized by quantum area, and let $(L_t, R_t)_{\R_+}$ be its boundary length process.  
Let $z_0$ be the point on the left boundary arc (i.e., $\R$) of the quantum wedge at distance $1$ from the root, and let $A$ be the time that $\eta'$ hits $z_0$. Let $a$ be the quantum length of $(\partial \eta([0,A])) \backslash \partial \cS$. Define the curve bottleneck event $F_{r,K} = \{ e^{\frac\gamma2(r+K)}a \in [1,2]\}$.

Conditioned on $F_{r, K}$, by Theorem~\ref{thm-MOT} the process $(L_t, R_t)_{[0, A]}$ evolves as Brownian motion with covariances~\eqref{eq-cov} stopped at the random time $A = \inf \{ t \: : \: L_t = -1\}$, and conditioned on $(R_{A} - \inf_{[0, A]} R_t) \in e^{\frac\gamma2(-r-K)}[1,2]$.  Purely using Brownian motion techniques, conditioned on $F_{r,K}$, in the $r \to \infty$ and $ K \to \infty$ limit, the process $(L_t, R_t)_{[0,	 A]}$ converges in distribution to $C \int_0^\infty \mu_{\R_+^2}^\gamma ( 1 + \ell i,0)\, d \ell$ for an appropriate constant $C$. See Figure~\ref{fig-MOT} (bottom).

This convergence in distribution of $(L_t, R_t)_{[0,A]}$ allows us to prove an analog of Proposition~\ref{prop-bulk-law}, where the limiting curve-decorated quantum surface is defined in terms of its peanosphere Brownian motion. More precisely, the process $(L_t, R_t)_{[0,A]}$ defines an equivalence relation $\sim$ on $[0,A]$, and the quotient space $[0,A]/\mathord\sim$ can be viewed as a topological disk decorated with a space-filling curve. Since the Brownian motion locally determines the field and curve in Theorem~\ref{thm-MOT} \cite[Theorem 1.11]{wedges}, this topological curve-decorated disk can be endowed with a conformal structure, so it can be viewed as a certain curve-decorated $\gamma$-LQG surface. At this step, we do not identify this limit as a space-filling $\SLE_{\kappa'}$-decorated quantum disk. 

Finally, we have the counterpart of Lemma~\ref{lem-curve-limit}: space-filling $\SLE_{\kappa'}$ is almost independent in spatially separated domains. As in Appendix~\ref{appendix-SLE}, this can be done by an imaginary geometry argument, using \cite[Lemma 2.4]{gms-harmonic} and \cite[Proposition 2.5(a)]{ag-disk}.

With these ingredients, the same argument in Section~\ref{sec-intrinsic} showing Proposition~\ref{prop-weak-theorem} can be applied: Conditioning on $\mathrm{BN}_{r, \zeta, K}$,  the bulk law of the field and curves is almost independent from their law near the field bottleneck, and sending $r,\zeta, K\to\infty$, the curve-decorated surface converges to a sample from $\cM_2^\disk(2 - \frac{\gamma^2}2; 1)$ decorated by an independent space-filling $\SLE_{\kappa'}$ curve. Conditioned on $\mathrm{BN}_{r, \zeta, K}$, the event $F_{r,K}$ is almost measurable with respect to the field and curve near the field bottleneck, and hence further conditioning on $F_{r,K}$ does not change the limit law. Finally, since $\P[\mathrm{BN}_{r,\zeta,K} \mid F_{r,K}] \approx 1$, the limit law of the  boundary length process conditioned on $\mathrm{BN}_{r,\zeta,K} \cap F_{r,K}$ is $C \int_0^\infty \mu_{\R_+^2}^\gamma ( 1+ \ell i,0) \,d \ell$. This yields Proposition~\ref{prop-pinch-MOT} for $\gamma \in (0,\sqrt2)$.

We adapt this argument for $\gamma = \sqrt2$ (it does not immediately apply since  $2-\frac{\gamma^2}2 = \frac{\gamma^2}2$ is the critical weight for thick quantum disks). Pick $W \in (0, \frac{\gamma^2}2)$ and consider a weight $2-\frac{\gamma^2}2 + W$ thick quantum wedge decorated by an $\SLE_\kappa(-\frac{\gamma^2}2; W -2)$ curve $\eta$. By Theorem~\ref{thm-wedges} the quantum surface $\cW_1$ (resp.\ $\cW_2$) to the left (resp.\ right) of $\eta$ is a weight $2-\frac{\gamma^2}2$ (resp.\ weight $W$) quantum wedge. 

Draw a space-filling $\SLE_{\kappa'}$ curve $\eta'$ on $\cW_1$. Then Theorem~\ref{thm-MOT} yields the boundary length process of $\eta'$ on $\cW_1$. Let $z_0$ be the point on $\R$ so $[z_0, \infty)$ has quantum length 1, let $D$ be the region explored by $\eta'$ up until it hits $z_0$, and let $U$ be the union of $D$ with the bounded components of $\cS \backslash D$. Define $F_{r,K} = \{ e^{\frac\gamma2(r+K)} \nu_h(\partial U \backslash \partial \cS) \in [1,2] \}$. Conditioning on $F_{r, K}$ and sending $r\to\infty$ and $K \to \infty$ (in that order),
we understand the limiting law of $(L_t, R_t)_{[0,A]}$, and the limiting decorated quantum surface is a weight $2-\frac{\gamma^2}2 + W$ thick quantum disk decorated by an $\SLE_{\kappa}(-\frac{\gamma^2}2; W-2)$ curve $\eta$ and a space-filling $\SLE_{\kappa'}$ curve $\eta'$ in the region to the left of $\eta$. By Theorem~\ref{thm-disk-cutting-2} cutting along $\eta$ gives a weight $2-\frac{\gamma^2}2$ quantum disk decorated by space-filling $\SLE_{\kappa'}$. 

For $\gamma \in (\sqrt2, 2)$, we have $2-\frac{\gamma^2}2 < \frac{\gamma^2}2$ so the weight $2-\frac{\gamma^2}2$ quantum wedge is thin. On the left boundary arc of a curve-decorated weight $2-\frac{\gamma^2}2$ quantum wedge, let $z_0$ be the point at distance 1 from the root, and let $\cD$ be the chain of quantum disks between $z_0$ and the root (not including the quantum disk containing $z_0$). Condition on the left side length of $\cD$ being at least $1-\eps$. Then Corollary~\ref{cor-thin-wedge-decomp} tells us that the pinched quantum surface is a quantum disk, and Theorem~\ref{thm-MOT} gives the  boundary length process.  Finally, sending $\eps \to 0$, the boundary length process converges to $C \int_0^\infty \mu_{\R_+^2}^\gamma (1 + \ell i,0) \,d \ell$.
\end{proof}
Observe that in the above proof, just as in the original proofs of the quantum disk and sphere mating-of-trees theorems, the proof  of Proposition~\ref{prop-pinch-MOT} is easier when $\gamma \in (\sqrt2,2)$. In this regime, the nontrivial topology that arises in mating of trees creates natural bottlenecks.

\begin{corollary}\label{cor-MOT}
In the setting of Proposition~\ref{prop-pinch-MOT} we replace  $\cM_2^\disk(2 - \frac{\gamma^2}2; 1)$ by $\cM_2^\disk(2-\frac{\gamma^2}2)$.
Then for some constant $C >0$, the law of $( L_t - L_A, R_t - R_A)_{[0,A]}$ is given by $C \iint_0^\infty \mu_{\R_+^2}^\gamma (\ell + \ell' i, 0)\, d \ell\, d\ell'$. 
\end{corollary}
\begin{proof}
This is a consequence of Proposition~\ref{prop-pinch-MOT} by scaling, using Lemma~\ref{lem-scaling-disint} and \eqref{eq-scale-BM}.
\end{proof}

We now state and give an alternative proof of the sphere variant of the mating-of-trees theorem \cite[Theorem 1.1]{sphere-constructions}. Let $\cP_\mathrm{SF}$ denote the law of space-filling $\SLE_{\kappa'}$ between two boundary points in a simply connected domain, and extend its definition to domains which are chains of disks by concatenation. The probability measure $\cP^\sph_\mathrm{SF}$ for space-filling $\SLE_{\kappa'}$ loops on a one-pointed sphere can be defined by arbitrarily picking a second point on the sphere, drawing a pair of curves from $\cP^\sph(2-\frac{\gamma^2}2, 2-\frac{\gamma^2}2)$, dividing the sphere into two (possibly beaded) parts, independently  sampling space-filling curves in each part  from $\cP_\mathrm{SF}$, and concatenating them; see \cite[Footnote 4]{wedges}.

\begin{theorem}[Quantum sphere mating-of-trees] 
Consider a sample $(\wh \C, h, 0, \infty, \eta')\sim \cM_2^\sph(4-\gamma^2) \otimes \cP^\sph_\mathrm{SF}$ with $\eta'$ parametrized by quantum area. 
Let $(L_t, R_t)$ be, respectively, the left and right quantum boundary lengths of $\eta([0,t])$. 
For some $C>0$, the law of the process $(L_t, R_t)$ is $C \mu^\gamma_{\R_+^2}(0,0)$ weighted by Brownian excursion duration. 
\end{theorem}
\begin{proof}[Sketch of proof]
Write $W = 2- \frac{\gamma^2}2$. Then by Theorem~\ref{thm-sphere-cutting} and the definition of $\cP^\sph_\mathrm{SF}$, we have 
\[\cM_2^\sph(4-\gamma^2) \otimes \cP^\sph_\mathrm{SF} = \wh c_{W,W} \iint_0^\infty \left(\cM_2^\disk(W; \ell, \ell') \otimes \cP_\mathrm{SF} \right) \times \left(\cM_2^\disk(W; \ell', \ell) \otimes \cP_\mathrm{SF} \right)\, d \ell\, d \ell'.\]
Thus, using Corollary~\ref{cor-MOT}, the  boundary length process has law given by 
\[C\iint_0^\infty \mu^\gamma_{\R_+^2}(0, \ell + \ell' i) \times \mu^\gamma_{\R_+^2}(\ell + \ell' i, 0)\, d \ell \,d \ell', \]
where $\mu^\gamma_{\R_+^2}(0, \ell + \ell' i)$ is defined via time-reversal from $\mu^\gamma_{\R_+^2}(\ell + \ell' i,0)$,
and a sample $(\gamma_1, \gamma_2)$ 
from $\mu^\gamma_{\R_+^2}(0, \ell + \ell' i) \times \mu^\gamma_{\R_+^2}(\ell + \ell' i,0)$
is interpreted as a path in $\R_+^2$ from 0 to 0 by concatenating $\gamma_1$ and $\gamma_2$.

We now show that the duration-weighted cone excursion measure agrees with the above law. 
Since the Brownian excursion measure can be written as a disintegration over the excursion duration $\mu_{\R_+^2}^\gamma(0,0) = \int_0^\infty \mu_{\R_+^2}^\gamma(0,0;t)\, dt$, 
the duration-weighted cone excursion measure can be written as $\iint_0^\infty \mu_{\R_+^2}^\gamma(0,0; t_1 + t_2)\, dt_1\,  dt_2$. By the Markov property of Brownian motion we have the path decomposition
\alb
\iint_0^\infty \mu_{\R_+^2}^\gamma(0,0; t_1 + t_2) \,dt_1\, dt_2 &=  \iiiint_0^\infty \mu_{\R_+^2}^\gamma(0,\ell+\ell'i; t_1) \times  \mu_{\R_+^2}^\gamma(\ell+\ell'i,0; t_2)\, d\ell\, d\ell' \, dt_1\, dt_2 \\
&=  \iint_0^\infty \mu_{\R_+^2}^\gamma(0,\ell+\ell'i)\times \mu_{\R_+^2}^\gamma(\ell+\ell'i,0)\, d\ell \,d \ell',
\ale
as desired.
\end{proof}

\begin{figure}[ht!]
	\begin{center}
		\includegraphics[scale=0.65]{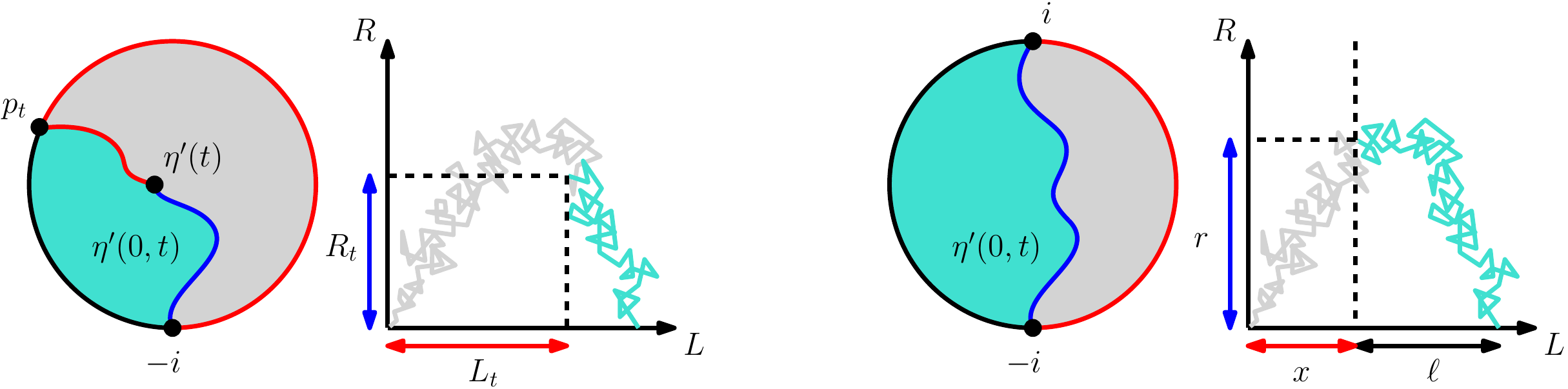}%
	\end{center}
	\caption{\label{fig-mot-gamma2-2} 
		 \textbf{Left:} Illustration of $(L_t, R_t)$ in Theorem \ref{thm-MOT-disk}, namely $(L_t, R_t)$ is the boundary length process of $(\D, h, \eta', -i)$ sampled from $\cM_1^\disk(2; a)^\# \otimes \cP^\circlearrowright_\mathrm{SF}$. With the space-filling curve parametrized by quantum area, the red quantum length is $L_t$ and the blue quantum length $R_t$. The process $(L_t, R_t)|_{[0,\mu_h(\D)]}$ is then a Brownian excursion in the cone $\R_+^2$. \textbf{Right: } Illustration of the proof of Proposition \ref{prop-left-right-law-2}. Sample $(\D, {h}, -i, i, \eta')$ from  $\cM_2^\disk(2;\ell, x) \otimes \cP^\circlearrowright_\mathrm{SF}$. The time $\eta'$ hits $i$ equals the stopping time $\tau = \inf \{ t : L_t \leq x\}$ of $(L_t, R_t)_{[0,\mu_h(\D)]}$. Since the blue interface is a certain $\SLE_\kappa(\rho_-;\rho_+)$ curve by SLE duality, the conditional law of the quantum surface $(\eta'([0,\tau]), h, -i, i)$ given $R_\tau$ is $\cM_2^\disk(\frac{\gamma^2}2; \ell, R_\tau)^\#$. This gives a Brownian motion description of the quantum area and lengths of a sample of $\cM_2^\disk(\frac{\gamma^2}2)$, which implies Proposition \ref{prop-left-right-law-2}.
}
\end{figure}

Finally, we can give an alternative proof of the disk variant of the mating-of-trees theorem \cite{wedges, ag-disk}. %
For $a>0$, let $(\D, h, -i, i)$ be a quantum surface sampled from $\cM_2^\disk(2;\frac a2, \frac a2)^\#$, and let $\cM_1^\disk(2; a)^\#$ be the law of the quantum surface $(\D, h, -i)$.  We can define the probability measure $\cP^\circlearrowright_\mathrm{SF}$ on space-filling $\SLE_{\kappa'}$ loops in a disk with a marked boundary point. %
This measure is obtained by sampling a curve from $\cP_\mathrm{SF}$ in a two-pointed domain $(\D, -i, x)$ and sending $x \to -i$ in the clockwise direction; see \cite[Appendix A.3]{bg-CRT} for details. For a sample $(\D, h, \eta', -i)\sim\cM_1^\disk(2; a)^\# \otimes \cP^\circlearrowright_\mathrm{SF}$, writing $A = \mu_h(\D)$, we define $(L_t, R_t)_{t \in [0,A]}$ as follows. Parametrize $\eta'$ so that at time $t$ we have $\mu_h(\eta'([0,t])) = t$. Let $p_t \in \partial \D$ be the furthest point on the clockwise arc from $-i$ so that the arc from $-i$ to $p_t \in\eta'([0,t])$, and let $L_t$ be the $\nu_h$-length of the boundary arc of $\eta'([0,t])$ from $\eta'(t)$ to $p_t$ plus the $\nu_h$-length of the clockwise arc of $\D$ from $p_t$ to $-i$. Let $R_t$ be the $\nu_h$-length of the boundary arc of $\eta'([0,t])$ from $\eta'(t)$ to $-i$. 
We call $(L_t, R_t)_{[0,A]}$ the \emph{boundary length process} of $(\D, h, \eta', -i)$. See Figure~\ref{fig-mot-gamma2-2}, left. 

\begin{theorem}[Quantum disk mating-of-trees]\label{thm-MOT-disk}
	For any $a>0$, the law of the  boundary length process $(L_t, R_t)$ of $\cM_1^\disk(2; a)^\# \otimes \cP^\circlearrowright_\mathrm{SF}$ is $\mu_{\R_+^2}^\gamma (a,0)^\#$. 
\end{theorem}
\begin{proof}[Sketch of proof]
Embed a quantum disk from $\cM_1^\disk(2; a)^\#$ as $(\D, h, -i)$ so that the boundary points $-i, -1, 1$ divide $\partial \D$ into arcs with quantum lengths $a/3$, and for $\eps > 0$ let $w_\eps \in \partial \D$ be the point so the clockwise arc from $w_\eps$ to $-i$ has quantum length $\eps$. Since marked points on $\cM_2^\disk(2)$ are independently and uniformly distributed according to the quantum length measure \cite[Proposition A.8]{wedges}, the quantum surface $(\D, h, -i, p_\eps)$ has law $\cM_2^\disk(2; a-\eps, \eps)^\#$. 

Let $\eta_\eps$ be an independent $\SLE_\kappa(-\frac{\gamma^2}2; \frac{\gamma^2}2-2)$ curve in $\D$ from $-i$ to $w_\eps$, and let its quantum length be $X$. Conditioned on $X$, by Theorem~\ref{thm-disk-cutting-2}  the quantum surface $\cD_\eps$ to the left of $\eta_\eps$ has conditional law $\cM_2^\disk(2 - \frac{\gamma^2}2; a-\eps, X)^\#$.  Draw an independent space-filling $\SLE_{\kappa'}$ curve $\eta'_\eps \sim \cP_\mathrm{SF}$ from $-i$ to $w_\eps$ in $\cD_\eps$. As $\eps \to 0$, we have $X \to 0$ in probability and hence, by Proposition \ref{prop-pinch-MOT}, the law of the boundary length process of $\eta'_\eps$ converges to $\mu_{\R_+^2}^\gamma(a,0)^\#$. Also as $\eps \to 0$, the curve $\eta_\eps$ degenerates to the point $-i\in \partial \D$ in probability, and hence the law of $\eta'_\eps$ converges to $\cP^\circlearrowright_\mathrm{SF}$. This yields the theorem.
\end{proof}

\subsection{Area and length distribution of  the weight $\gamma^2/2$ quantum disk}
\label{subsec-MOT}
The main result of this subsection is the following proposition, which expresses the area and length distribution of the weight $\gamma^2/2$ quantum disk in terms of Brownian measures, and gives an explicit formula for the joint law of the two boundary lengths.
\begin{proposition}\label{prop-MOT}
	There are constants $C,C_0 \in (0, \infty)$ such that for all $\ell, r > 0$, 
	\eqb
	|\cM^\disk_2(\frac{\gamma^2}2; \ell, r)| 
	= C_0|\mu^\gamma_{\R_+^2}(\ell, ri)|	= C \frac{(\ell r)^{4/\gamma^2 - 1}}{(\ell^{4/\gamma^2}+r^{4/\gamma^2})^2}.
	\label{eq:joint-g2-2}
	\eqe
	Moreover, the quantum area of a sample from $\cM^\disk_2(\frac{\gamma^2}2; \ell, r)^\#$ agrees in law with the duration of a sample from $\mu^\gamma_{\R_+^2}(\ell, ri)^\#$.
\end{proposition}
Before giving the proof of the proposition we state a counterpart of \eqref{eq:joint-g2-2} for the weight 2 quantum disk. The constants $C$ appearing in the propositions is identified in our companion work \cite{AHS-SLE-integrability}.
\begin{proposition}\label{prop-left-right-law-2}
	There is a constant $C>0$ such that 
	\[|\cM_{2}^\disk(2;\ell, r)| = C (\ell + r)^{-\frac4{\gamma^2}-1} \quad \text{for }\ell,r>0.
	\]
\end{proposition}
\begin{proof}
	By Lemma~\ref{lem-length-exponent}, the law of the total quantum boundary length of $\cM_2^\disk(2)$ is $C \ell^{-\frac4{\gamma^2}}\, d\ell$. By \cite[Proposition A.8]{wedges}, if we sample $(\D, h, -1, 1)$ from $\cM_2^\disk(2)$ and independently sample $x, y \in \partial \D$ from $\frac{\nu_h}{\nu_h(\partial \D)}$, then the quantum surface $(\D, h, x, y)$ still has law $\cM_2^\disk(2)$. Thus the joint law of the left and right boundary lengths of $\cM_2^\disk(2)$ is given by $C (\ell + r)^{-\frac4{\gamma^2} - 1} \, d\ell \, dr$.
\end{proof}

We will now give the proof of Proposition~\ref{prop-MOT} 
using Theorem~\ref{thm-MOT-disk}
and exact formulas for the associated planar Brownian motion. We first state a variant of Theorem~\ref{thm-MOT-disk}.

\begin{lemma}
	\label{lemma-MOT-weight2}
	There is a constant $C$ for which the following holds. 
	Let $\ell, x>0$, and sample $(\D, h, \eta', -i, i)$ from $\cM_2^\disk(2; \ell, x) \otimes \cP^\circlearrowright_\mathrm{SF}$. Then the law of $(L_t, R_t)_{[0,\mu_h(\D)]}$ is $C \mu^\gamma_{\R_+^2} ( \ell + x, 0)$, where the process $(L_t, R_t)_{[0,\mu_h(\D)]}$ is defined as in Theorem~\ref{thm-MOT-disk}.
\end{lemma}
\begin{proof}
	Theorem~\ref{thm-MOT-disk} states that if $(\D, h, \eta', -i, i)$ is sampled from $\cM_1^\disk(2; \ell+x)^\# \otimes \cP^\circlearrowright_\mathrm{SF}$ then the law of $(L_t, R_t)_{[0,\mu_h(\D)]}$ is $\mu^\gamma_{\R_+^2} (\ell + x, 0)^\#$. By \cite[Proposition A.8]{wedges} a sample from $\cMtwo(2;\ell, r)^\#$ can be obtained by sampling from $\QD_1(\ell + x)^\#$ and marking another boundary point to get length $\ell$ and $x$ boundary arcs. Since $|\cMtwo(2;\ell, x)| \propto (\ell + x)^{-\frac4{\gamma^2}-1}$ by  Proposition~\ref{prop-left-right-law-2} and $|\mu^\gamma_{\R_+^2}(\ell + x, 0)| \propto (\ell + x)^{-\frac4{\gamma^2}-1}$ by~\eqref{eq-scale-BM}, we get the result. 
\end{proof}

\begin{proof}[Proof of Proposition \ref{prop-MOT}]
	See Figure~\ref{fig-mot-gamma2-2}, right. 
	For $(\D, {h}, -i, i, \eta')$ sampled from  $\cM_2^\disk(2;\ell, x) \otimes \cP^\circlearrowright_\mathrm{SF}$, let $A = \mu_{h}(\D)$,  let $\tau$ be the time that $\eta'$ hits $i$, and let $D_1 = \eta'([0,\tau])$.
	In \cite[Section 3.6.3]{ghs-mating-survey} it is given a construction of $\eta'$ such that, when $\eta$ is the right boundary of a certain $\SLE_{\kappa'}(\kappa'-6)$ curve in $(\D, -i, i)$, then $D_1$ is the region to the left of $\eta$ in $\D$. SLE duality (see \cite[Theorem 5.1]{zhan-duality1} and \cite[Theorem 1.4]{ig1}) says that the law of the right boundary of an $\SLE_{\kappa'}(\rho_-'; \rho_+')$ curve is $\SLE_\kappa(\rho_-; \rho_+)$ where $\rho_- = \frac\kappa2 - 2 + \frac\kappa4 \rho_-'$ and $\rho_+= \kappa +\frac\kappa4 \rho_+' -4$; since $\kappa'  = \frac{16}\kappa, \rho_-' = 0$ and $\rho_+' = \kappa'-6$, the law of $\eta$ is $\SLE_\kappa(\frac\kappa2-2; -\frac\kappa2)$.
	Hence, by Theorem~\ref{thm-disk-cutting-2}, 
	the marginal law of the quantum surface $(D_1,h|_{D_1},-i,\eta'(\tau))$ is given by 
	\[\int_0^\infty C_1  \left|\cMtwo(2-\frac{\gamma^2}2; r, x)\right| \cdot \cMtwo(\frac{\gamma^2}2; \ell, r) \, dr\quad \textrm{for some constant }C_1\in(0,\infty).\]

	Lemma~\ref{lemma-MOT-weight2} tells us that for some $C_2\in (0,\infty)$ the  boundary length process
	$(L_t, R_t)_{[0,\mu_{h}(\D)]}$ of $(\D, {h}, -i, i, \eta')$ has law $C_2 \mu^\gamma_{\R_+^2} (\ell + x, 0)$. The strong Markov property of Brownian motion applied to the stopping time $\tau$ yields a path decomposition \cite[Proposition 1.11]{lawler-notes}: the law of $((L_t-x , R_t)_{[0,\tau]}, (L_{s + \tau}, R_{s+\tau})_{[0, A - \tau]})$ is $C_2\int_0^\infty \mu^\gamma_{ \R_+^2} (\ell, ri) \times \mu^\gamma_{\R_+^2}(x+ri,0) \, dr$. Hence the marginal law of $(L_t-x, R_t)_{[0,\tau]}$ is 
	\eqb \label{eq-MOT-peanosphere-split}
	\int_0^\infty C_2 \left|\mu^\gamma_{\R_+^2}(x+ri,0)\right| \cdot  \mu^\gamma_{\R_+^2} (\ell, ri) \, dr.
	\eqe
	
	Disintegrating over $r$, we see that for every $\ell, r, x > 0$ we have
	\[C_1  \left|\cMtwo(2-\frac{\gamma^2}2; r, x)\right| \cdot \left|\cMtwo(\frac{\gamma^2}2; \ell, r)  \right| = C_2 \left|\mu^\gamma_{\R_+^2}(x+ri,0)\right| \cdot \left|  \mu^\gamma_{\R_+^2} (\ell, ri)\right|.\]
	Integrating over $x > 0$, we have $\int_0^\infty |\cMtwo(2-\frac{\gamma^2}2; r, x)| \, dx = |\cMtwo(2-\frac{\gamma^2}2; r)| \propto r^{-\frac4{\gamma^2}+1}$ by Lemma~\ref{lem-length-exponent} (if $2-\frac{\gamma^2}2 \geq \frac{\gamma^2}2$) or Lemma~\ref{lem-thin-disk-perim} (if $2-\frac{\gamma^2}2 < \frac{\gamma^2}2$), 
	Moreover,  by~\eqref{eq-scale-BM},
	\[
	\int_0^\infty \left|\mu^\gamma_{\R_+^2}(x+ri,0)\right| \, dx   = r^{-\frac{4}{\gamma^2}}  \int_0^\infty \left|\mu^\gamma_{\R_+^2}(\frac{x}r+i,0)\right| \, dx  \propto r^{-\frac4{\gamma^2}+1}.
	\] 
	Thus we get $|\cMtwo(\frac{\gamma^2}2; \ell, r)| = C|\mu^\gamma_{\R_+^2}(\ell, ri)|$ for $\ell, r > 0$ for some $C>0$. The second claim follows from the fact that quantum area corresponds to the Brownian excursion duration. 
\end{proof}

\begin{remark}[Mating-of-trees for weight $\frac{\gamma^2}2$ quantum disk]
	The argument of Proposition~\ref{prop-MOT} shows that for some constant $C>0$, the following holds for all $\ell, r>0$. Sample a quantum disk from $\cMtwo(\frac{\gamma^2}2; \ell, r)$ and decorate it by an independent \emph{space-filling $\SLE_{\kappa'}(0; \frac{\kappa'}2-4)$ curve} between its marked boundary points \cite{ig4}. Then a suitably defined process $(L_t, R_t)$, which can be viewed as a boundary length process for the SLE-decorated LQG surface, has law $C \mu_{\R_+^2}(\ell, ri)$. 
\end{remark}

\appendix

\section{Extension of welding results to $W=\frac{\gamma^2}{2}$}
\label{appendix-A}
In the proof of Theorem~\ref{thm-disk-cutting}  in Section \ref{sec-general} we break the argument into five cases. The first three cases, which we have proved using the input from Sections~\ref{sec-extrinsic}---\ref{sec-intrinsic}, can be summarized as follows.
\begin{proposition}
	Theorem \ref{thm-disk-cutting} holds in the case when $W_1,\dots,W_n,W\in(0,\infty)\setminus\{\frac{\gamma^2}{2} \}$. 
	\label{prop-gn2o2}
\end{proposition}
In this appendix we will complete Case 4 of the proof of  Theorem \ref{thm-disk-cutting} based on Proposition~\ref{prop-gn2o2}. 
Namely, we  extend Theorem \ref{thm-disk-cutting} to include the situation where $n=2$ and one  weight equals $\frac{\gamma^2}{2}$.
\begin{proposition}
	Theorem \ref{thm-disk-cutting} holds in the case when $n=2$, $W_1=2$,
	and $W_2=\frac{\gamma^2}{2}$. 
	\label{prop-g2o2}
\end{proposition}
We prove  Proposition~\ref{prop-g2o2} by considering the $W_1=2$ and $W_2=\frac{\gamma^2}{2}+\eps$ case and then sending $\eps\downarrow 0$.
We start by proving  a continuity result for thick disks in the weight parameter.
\begin{lem}
	Let $W\geq \frac{\gamma^2}{2}$, $\eps\geq 0$, and $a\in(0,1)$. Let 
	$\cD^\eps\sim (\cM_2^\disk(W+\eps)|_{A(a)})^\#$, where $A(a)$ is the event that the left and right LQG boundary lengths of $\cD^\eps$ are in $(a,a^{-1})$. Let $h^\eps$ be such that $(\cS,h^\eps,+\infty,-\infty)$ is an 
	embedding of $\cD^\eps$ for which $\R_+$ has LQG length $a$. %
	Then $h^\eps$ restricted to any bounded set converges to $h^0$ in total variation distance.
	
	The same result holds if we let $A_{\op{L}}(a)$ (resp.\ $A_{\op{R}}(a)$) be the event that the left (resp.\ right) LQG boundary length of the surface is in $(a,a^{-1})$ and we replace $A(a)$ by $A_{\op{L}}(a)$ (resp.\ $A_{\op{R}}(a)$) in the above statement, where, for the case of $A_{\op{R}}(a)$, we embed the surfaces such that $\R_++i\pi$ (instead of $\R_+$) has LQG length $a$.
	\label{lem-weight-cont} 
\end{lem}
\begin{proof}
	Consider Definition \ref{def-thick-disk} with weight parameter $W+\eps$ instead of $W$ and write $\cc_\eps$ instead of $\cc$. For $\delta>0$, let $\wt A(\delta)$ denote the event that $\cc_\eps\geq \delta$. On the event $\wt A(\delta)$ and normalizing the measure from which $\cc_\eps$ is sampled to be a probability measures, $\cc_\eps$ converges %
	in total variation distance to $\cc_0$ as $\eps\rta 0$. Let $h_1^\eps$ and $h_2^\eps$ be the field $\wh h$ in Definition \ref{def-thick-disk} projected onto $\cH_1(\cS)$ and $\cH_2(\cS)$, respectively. Notice that the law of $h_2^\eps$ does not depend on $\eps$, while $h_1^\eps$ restricted to any bounded set converges to $h_1^0$ for the total variation distance. %
	Combining the convergence results for $\cc_\eps$ and $h_1^0$, we get that the lemma holds with $(\cM_2^\disk(W+\eps)|_{\wt A(\delta)})^\#$ instead of $(\cM_2^\disk(W+\eps)|_{A(a)})^\#$. Furthermore, this convergence is joint with convergence of the event that the left and right LQG boundary lengths defined by $h_1^\eps+h_2^\eps+\cc_\eps$ are in $(a,a^{-1})$. We obtain the lemma from this by using that $(\cM_2^\disk(W+\eps)|_{A(a)\cap \wt A(\delta)})^\#$ converges in total variation distance to $(\cM_2^\disk(W+\eps)|_{A(a)})^\#$ as $ \delta\rta 0$, uniformly in $\eps\in[0,1]$. The proof for the events $A_{\op{L}}(a)$ and $A_{\op{R}}(a)$ is identical.
\end{proof}
We will also need a continuity result for SLE.
\begin{lemma}
	For $\eps\geq 0$ let $\eta_\eps$ be an $\op{SLE}_\kappa(0;\gamma^2/2+\eps-2)$ on $(\cS,+\infty,-\infty)$, let $D^1_\eps\subset\cS$ (resp.\ $D^2_\eps\subset\cS$) be the domain below (resp.\ above) $\eta_\eps$, 
	let $\phi_1^\eps:\cS\to D^\eps_1$ be the conformal map which is fixing $\pm\infty$ and $0$, and
	let $\phi_2^\eps:\cS\to D^\eps_2$ be the conformal map which is fixing $\pm\infty$ and $i\pi$. Then $\phi_1^\eps$ (resp.\ $\phi_2^\eps$) is converging uniformly in law to $\phi_1$ (resp.\ $\phi_2$) on compact subsets of $\cS\cup\R$ (resp.\ $\cS\cup( \R+i\pi )$), and the convergence is joint for $\phi_1$ and $\phi_2$.
	\label{lem:sle-conv}
\end{lemma}

 Our proof of the lemma relies on the following result, which is a variant of \cite[Lemma 6.1]{kemppainen-book}. %
The main extension as compared to \cite{kemppainen-book} is that the set $A$ is not required to be bounded away from $\R$. 
\begin{lemma} 
	Let $\eta$ and $\wt\eta$ be curves in $\BB H$ from 0 to $\infty$ with Loewner driving function $(W_t)_{t\geq 0}$ and $(\wt W_t)_{t\geq 0}$, respectively, and let $(g_t)_{t\geq 0}$ and $(\wt g_t)_{t\geq 0}$ denote the Loewner maps. For any $\eps\in(0,1)$ there is a $\delta\in(0,1)$ such that if
	$$ 
	A = \{ (t,z) \in [0,T]\times\ol{\BB H}\,:\,  \inf_{s\in[0,t]}|g_s(z)-W_s|>\eps \}
	\quad\text{and}\quad
	\sup_{t\in[0,T]}|W_t-\wt W_t|\leq\delta.
	$$ 
	then 
	$$
	\sup_{(t,z)\in A}|g_t(z)- \wt g_t(z)| < \eps.
	$$ 
	\label{lem:kempp}
\end{lemma}
\begin{proof}
	In the proof of \cite[Lemma 6.2]{kemppainen-book} it was argued that 
	\eqb
	|g_t(z)-\wt g_t(z)|
	\leq \sup_{s\in[0,t]}|W_s- \wt W_s|( \exp(\sqrt{ I(t) \wt I(t)}) - 1 ),
	\label{eq:loewner-delta}
	\eqe
	where
	\eqbn
		I(t) = 2\int_0^t |g_s(z)-W_s|^{-2}\,ds,\qquad
		\wt I(t) = 2\int_0^t |\wt g_s(z)-\wt W_s|^{-2}\,ds. 
	\eqen
	Set
	\eqbn
	\delta = \frac{\eps}{3}( \exp(4t\eps^{-2}) - 1 )^{-1} \wedge \frac{\eps}{10}.
	\eqen
	Suppose there is a time $s_0\in[0,t]$ such that 
	\eqb
	|\wt g_{s_0}(z)-\wt W_{s_0}|\leq\eps/2.
	\label{eq:W-delta}
	\eqe
	Let $s_0$ be the smallest time satisfying this requirement. Then we get from \eqref{eq:loewner-delta} that
	$|g_{s_0}(z)-\wt g_{s_0}(z)|<\eps/3$, so by the triangle inequality
	\eqbn
	|\wt g_{s_0}(z)-\wt W_{s_0}|
	\geq |g_{s_0}(z)-W_{s_0}|
	-|\wt g_{s_0}(z)-g_{s_0}(z)|
	-|\wt W_{s_0}-W_{s_0}|
	>\eps-\frac{\eps}{3}-\frac{\eps}{10}>\frac{\eps}{2}.
	\eqen
	This is a contradiction to the definition of $s_0$, and we conclude that there is no time $s_0\in[0,t]$ satisfying \eqref{eq:W-delta}. Since there is no time $s_0\in[0,t]$ satisfying \eqref{eq:W-delta}, we have $I(t)<2t\eps^{-2}$ and $\wt I(t)<8t\eps^{-2}$, so the right side of \eqref{eq:loewner-delta} is smaller than $\delta(\exp(4t\eps^{-2})-1)<\eps$. 
\end{proof}

\begin{proof}[Proof of Lemma \ref{lem:sle-conv}]
	We will prove the lemma in the setting of the upper half-plane $\BB H$ instead of $\cS$. Precisely, if $\wh\eta_\eps$ is an $\op{SLE}_\kappa(0;\gamma^2/2+\eps-2)$ on $(\BB H,0,\infty)$, $D^{\op{L}}_\eps\subset\BB H$ (resp.\ $D^R_\eps\subset\BB H$) is the domain to the left (resp.\ right) of $\wh\eta_\eps$, and
	$\phi_{\op{L}}^\eps:\BB H\to D^\eps_{\op{L}}$ 
	(resp.\ $\phi_{\op{R}}^\eps:\BB H\to D^\eps_{\op{R}}$) is the conformal map which is fixing $\infty$, $0$, and $1$ (resp.\ $-1$), then we will argue that $\phi_{\op{L}}^\eps$ (resp.\ $\phi_{\op{R}}^\eps$) is converging uniformly in law to $\phi_{\op{L}}$ 
 (resp.\ $\phi_{\op{R}}$) on compact subsets of $\BB H\cup\R_-$ (resp.\ $\BB H\cup \R_+$), and that the 
 convergence is joint for $\phi_{\op{L}}$ and $\phi_{\op{R}}$. The lemma immediately follows from this upon considering the conformal map $z\mapsto-\log z$ from $(\BB H,0,\infty)$ to $(\cS,+\infty,-\infty)$.
 
 Consider a coupling of the Loewner driving functions $(W^\eps(t))_{t\geq 0}$ of $\wh\eta_\eps$ such that $W^\eps$ converges uniformly to $W^0$ on compact sets a.s.; we leave the proof of existence of such a coupling as an exercise to the reader. Let $T>0$ be a large constant to be chosen later, and let $(g^\eps_t)_{t\geq 0}$ be the centered forward Loewner maps of $\wh\eta_\eps$. We can write %
 $(\phi_{\op{L}}^\eps)^{-1}=\psi^\eps_T\circ\theta^\eps_T \circ g^\eps_T$, where $\theta^\eps_T:\BB H\to\BB H$ is given by $\theta^\eps_T(z)=\frac{z-g^\eps_T(0^-)}{g^\eps_T(0^-)-g^\eps_T(-1)}$ and $\psi_T^\eps$ is the conformal map fixing $0$, $-1$, and $\infty$ which sends the domain to the left of $\theta^\eps_T \circ g^\eps_T(\wh\eta_\eps|_{[T,\infty)})$ to $\BB H$.
 
 By \ref{lem:kempp}, $g_T^\eps$ converges uniformly to $g_T^0$ on any subset of $\ol{\BB H}$ which is bounded away from $\wh\eta_0([0,T])$.
 In particular, $g_T^\eps(-\delta)\rta g_T^0(-\delta)$ a.s.\ as $\eps\rta 0$ for any fixed $\delta>0$, and by harmonic measure considerations we get further that $g_T^\eps(0^-)\rta g_T^0(0^-)$ a.s.\ as $\eps\rta 0$. This implies further that  $\theta^\eps_T$ converges uniformly to $\theta^0_T$ since $g^\eps_T(-1)\rta g^0_T(-1)$ and $g^\eps_T(0^-)\rta g^0_T(0^-)$ a.s.
 
 Fix a compact set $K\subset\BB H\cup\R_-$. The harmonic measure of $\phi_{\op{L}}^{-1}(\wh\eta_\eps|_{[T,\infty)})$ as seen from any point $z\in K$ goes to zero as $T$ goes to $\infty$, uniformly in $\eps$ and $z$; this follows e.g.\ from Brownian motion considerations and by using the fact that $\theta^\eps_T \circ g^\eps_T(\wh\eta_\eps|_{[T,\infty)})$ converges in Carath\'eodory topology as $\eps\rta 0$. Therefore, for any $\delta>0$ we can find $T$ sufficiently large such that $|\psi_T^\eps(z)-z|<\delta$ for all $z\in(\psi_T^\eps)^{-1}(K)$. Combining this with the previous paragraph we get a.s.\ convergence of $\phi_{\op{L}}^\eps$ to $\phi_{\op{L}}$ as desired. Convergence of $\phi_{\op{R}}^\eps$ to $\phi_{\op{R}}$ follows by a similar argument.
\end{proof}

\begin{proof}[Proof of Proposition \ref{prop-g2o2}]
For $\eps\geq 0$, $a>0$, and with $A(a)$ as in Lemma \ref{lem-weight-cont} let $\cD^\eps\sim(\cM^{\disk}_2(W_1+\frac{\gamma^2}{2}+\eps)|_{A(a)})^\#$, let $\eta^\eps\sim \cP^\disk(W_1, \frac{\gamma^2}{2}+\eps)$, and let $\cD^\eps_1$ and $\cD^\eps_2$ denote the surfaces to the left and right, respectively, of $\eta^\eps$. Note that all the considered surfaces and $\eta^\eps$ are sampled from probability measures. By Proposition \ref{prop-gn2o2}, if $\eps>0$ then $\cD^\eps_1$ and $\cD^\eps_2$ have the law of surfaces sampled from  $(\cM^{\disk}_2(W_1)|_{A_{\op{L}}(a)})^\#$ and $(\cM^{\disk}_2(\frac{\gamma^2}{2}+\eps)|_{A_{\op{R}}(a)})^\#$, respectively, and the surfaces are independent conditioned on the event that the right LQG boundary length of the former surface is equal to the left LQG boundary length of the latter surface. To conclude it is sufficient to argue that $\cD^0_1$ and $\cD^0_2$ have the law of surfaces sampled from $(\cM^{\disk}_2(W_1)|_{A_{\op{L}}(a)})^\#$ and $(\cM^{\disk}_2(\frac{\gamma^2}{2})|_{A_{\op{R}}(a)})^\#$, respectively, again such that the surfaces are independent given the same condition on the LQG boundary lengths as before. This is sufficient since it gives \eqref{eq-disk-cutting-2} restricted to the events $A(a),A_{\op{L}}(a),A_{\op{R}}(a)$ and with all the measures normalized to be probability measures. We get the case of non-probability measures by choosing the constant $c_{W_1,\frac{\gamma^2}{2}}$ appropriately so that the measures on the left and right sides of \eqref{eq-disk-cutting-2} have the same total mass, and sending $a\rta 0$ we can remove the constraint on the boundary lengths.
	
	Let $h^\eps$ be the field on $\cS$ such that $(\cS,h^\eps,+\infty,-\infty)$ is the embedding of $\cD^\eps$ for which $\R_+$ has LQG length $a$. Define $h^\eps_1,h^\eps_2$ in the same way for $\cD^\eps_1,\cD^\eps_2$, respectively, except we require that $h^\eps_2$ induces the same length on $\R_++i\pi$ as $h$. Let $D^\eps_1,D^\eps_2\subset\cS$ denote the domains to the left and right, respectively, of $\eta^\eps$. For $j=1,2$ let $\phi^\eps_j:\cS\to D^\eps_j$ be the conformal map such that $h^\eps|_{D^\eps_j}$ and $h^\eps_j$ are related by doing a coordinate change as in \eqref{eq-quantum-surface}, so in particular $\phi^\eps_1$ fixes $\pm\infty$ and 0 while
	$\phi^\eps_2$ fixes $\pm\infty$ and $i\pi$.
	
	By Lemma \ref{lem-weight-cont}, $h^\eps$ restricted to any compact set converges in total variation distance to $h^0$ as $\eps\rta 0$. By Lemma \ref{lem:sle-conv}, $\phi^\eps_1$ and $\phi^\eps_2$ converge jointly in law to $\phi^0_1$ and $\phi^0_2$, respectively, uniformly on compact subsets of $\cS\cup\R$ and $\cS\cup(\R+i\pi)$. Since $h^\eps$ and $(\phi^\eps_1,\phi^\eps_2)$ are independent there is a coupling so that they converge jointly a.s.\ to $h^0$ and $(\phi^0_1,\phi^0_2)$, respectively. Since  $\phi^\eps_1$ and $\phi^\eps_2$ are conformal we also get a.s.\ uniform convergence of their derivatives on compact sets. Therefore, for any smooth compactly supported test function $f$ on $\cS$,
	\eqbn
	\begin{split}
		(h^\eps_1,f) = 
		(h^\eps\circ\phi^\eps_1+Q\log|(\phi^\eps_1)'|,f)\rta 
		(h^0   \circ\phi^0_1   +Q\log|(\phi^0_1)'|,f)
		=(h^0_1,f)\,\,\,\text{a.s.},
	\end{split}
	\eqen
	so $h^\eps_1$ converges a.s.\ to $h^0_1$ for the weak-* topology. The same holds for $h^\eps_2$. By Lemma \ref{lem-weight-cont}, any limit of $h^\eps_1,h^\eps_2$ describe quantum surfaces sampled from $(\cM^{\disk}_2(W_1)|_{A_{\op{L}}(a)})^\#$ and $(\cM^{\disk}_2(\frac{\gamma^2}{2}+\eps)|_{A_{\op{R}}(a)})^\#$, respectively, conditioned on the right boundary arc of the former surface having an equal length as the left boundary arc of the latter surface. In particular, the limiting fields $h^0_1,h^0_2$ have the desired laws, which concludes the proof. 
\end{proof}

\section{Proof of SLE local independence lemma}\label{appendix-SLE}
The goal of this section is to prove Lemma~\ref{lem-curve-limit}. To clarify the picture we work in a bounded domain. Let $D$ be the square $[-1,1]^2$ and let $x = -i, y = i$. Let $U$ = $[-1,1] \times [-1,0]$ denote the lower half of $D$ and let $U_\eps = B_\eps(i) \cap D$ for $\eps > 0$. 

\begin{proposition}\label{prop-curve-limit-variant}
Suppose $n \geq 2$ and $W_1, W_2, \dots, W_n > 0$. Sample curves $(\eta_1, \dots, \eta_{n-1})\sim \cP^\disk(W_1, \dots, W_n)$ from $x$ to $y$. Let $\eta_j^\mathrm{start}$ be the initial segment of $\eta_j$ run until it exits $U$, and let $\eta_j^\eps$ be the initial segment of the time-reversal of $\eta_j$ run until it  exits $U_\eps$. Then the total variation distance between the following two laws is $1-o_\eps(1)$:
\begin{itemize}
\item The joint law of $(\eta_1^\mathrm{start}, \dots, \eta_{n-1}^\mathrm{start})$ and $(\eta_1^\eps,\dots, \eta_{n-1}^\eps)$. 

\item The joint law of $(\wt \eta_1^\mathrm{start}, \dots, \wt \eta_{n-1}^\mathrm{start})$ and $(\eta_1^\eps,\dots,  \eta_{n-1}^\eps)$, where $(\wt \eta_1, \dots, \wt \eta_{n-1})$ is independently sampled from $\cP_D(W_1, \dots, W_n)$ and $\wt \eta_j^\mathrm{start}$ is defined analogously as $\eta_j^\mathrm{start}$ for each $j$.  
\end{itemize}
\end{proposition}
Before giving the proof of this proposition,
we explain how it yields Lemma~\ref{lem-curve-limit}.
\begin{proof}[Proof of Lemma~\ref{lem-curve-limit}]
Proposition~\ref{prop-curve-limit-variant} yields a variant of Lemma~\ref{lem-curve-limit} where, instead of taking the intersections of the curves with $\cS_+$ and $\cS_- - N$, we instead take the curve tips run until they exit these two domains. Because the curve tips never revisit their starting points, there is some random $T > 0$ for which the restrictions of the curve tips to $\cS_+ + T$ and $\cS_- - N - T$ agree with the restrictions of the curves to these regions.
Therefore Lemma~\ref{lem-curve-limit} follows by looking at the curve tips intersected with $\cS_+ + M$ and $\cS_- - N - M$, and sending $M \to \infty$ and then $N \to \infty$. 
\end{proof}

We will understand the single curve ($n=2$) case of Proposition~\ref{prop-curve-limit-variant} using the framework of \emph{imaginary geometry} \cite{ig1, ig2}. Then we explain the minor modifications needed for the general $n$ regime. 

Consider the $n=2$ case and drop the subscript on the curve, i.e., $\eta := \eta_1$. Let $\rho_j = W_j - 2$ for $j=1,2$, so the curve $\eta$ is an $\SLE_\kappa(\rho_1;\rho_2)$ curve in $(D,x,y)$. One can couple $\eta$ with an appropriate Dirichlet boundary GFF $h^\IG$ in $D$, such that $\eta$ is an \emph{angle $\frac\pi2$ flow line} of $h^\IG$. Precisely, when we parametrize by $(\bbH, 0, \infty)$ the imaginary geometry GFF has boundary values $\frac{\pi}{\sqrt \kappa}(\frac\kappa4+ \rho_2)$ on $\R_+$ and $-\frac\pi{\sqrt\kappa}(2-\frac\kappa4+\rho_1)$ on $\R_-$, and $h^\IG$ has boundary values derived from this by an imaginary geometry coordinate change as defined in \cite{ig1}.
By \cite[Theorem 1.1]{ig1}, $\eta$ is a deterministic function of $h^\IG$ and $\eta^\mathrm{start}$ is determined by $h^\IG|_U$. 

Although the reversibility of $\SLE_\kappa(W_1-2; W_2-2)$ could suggest that $\eta^\eps$ is a deterministic function of $h^\IG|_{U_\eps}$, this turns out not to be the case. Instead we need to use the machinery of \emph{counterflow lines}. Let $\kappa' = \frac{16}\kappa$. One can couple with $h^\IG$ a certain $\SLE_{\kappa'}(\frac{\kappa'}2 - 2 + \frac{\kappa'}4 \rho_1; \kappa' - 4 + \frac{\kappa'}4 \rho_2)$ curve $\eta'$ from $y$ to $x$ such that $\eta'$ is a deterministic function of $h^\IG$, and, writing $\eta'^\delta$ for the initial segment of $\eta'$ run until it exists $U_\delta$, the segment $\eta'^\delta$ is determined by $h|_{U_\delta}$ \cite[Theorem 1.1]{ig1}. 

\begin{figure}[ht!]
\begin{center}
  \includegraphics[scale=0.75]{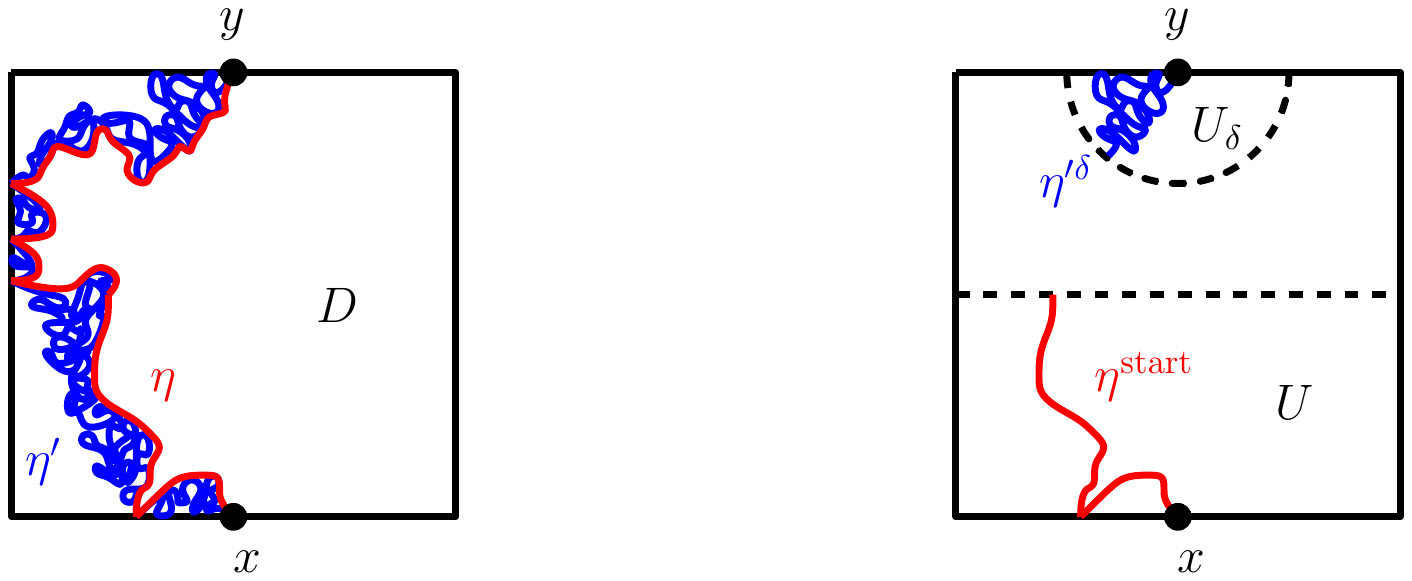}%
 \end{center}
\caption{\label{fig-indep-curves} \textbf{Left.} For the imaginary geometry $h^\IG$, the angle $\frac\pi2$ flow line $\eta$ is a.s.\ the right boundary of the counterflow line $\eta'$ (Lemma~\ref{lem-B-right-bdy}). \textbf{Right.} The initial segments $\eta^\mathrm{start}$ and $\eta'^\delta$ are a.s.\ determined by $h|_U$ and $h|_{U_\delta}$, respectively.}
\end{figure}

\begin{lemma}[{\cite[Theorem 1.4]{wedges}}]\label{lem-B-right-bdy}
Almost surely $\eta$ is the right boundary of $\eta'$. 
\end{lemma}
\begin{lemma}\label{lem-B-avoid-y}
Fix $\delta > 0$. On an event of probability $1-o_\eps(1)$, the curve segment $\eta^\eps$ is determined by $h|_{U_\delta}$. 
\end{lemma}
\begin{proof}
Write $E_{\delta, \eps}$ for the event that the curve $\eta'$ does not revisit $U_\eps$ after leaving $U_\delta$. Since $\eta'$ a.s.\ does not hit $y$ after leaving $U_\delta$, and $\eta'$ is a continuous curve, we conclude that $\P[E_{\delta, \eps}] = 1-o_\eps(1)$. The assertion then follows since on $E_{\delta, \eps}$, by Lemma~\ref{lem-B-right-bdy} $\eta^\eps$ is determined by $\eta'^\delta$, which is determined by $h|_{U_\delta}$. 
\end{proof}
\begin{lemma}\label{lem-B-local-indep}
Let $\wt h^\IG$ be independently sampled, with the same law as $h^\IG$. Then the total variation distance between the laws of $(h^\IG|_U, h^\IG|_{U_\delta})$ and $(\wt h^\IG|_U, h^\IG|_{U_\delta})$ is $1-o_\delta(1)$. 
\end{lemma}
\begin{proof}
We work in the strip $(\cS, +\infty, -\infty)$ instead. The corresponding imaginary geometry field $\wh h^\IG$ in $\cS$ has constant boundary conditions on $\R$ and $\R + i\pi$ (with different values on each line). Let $\wh U$ be a neighborhood of $+\infty$ excluding $-\infty$, and let $\wh U_N = (-\infty, -N) \times [0,\pi]$. Let $\wh V = \cS \backslash \wh U$, and let $I_1 = \partial \wh V \cap \R$, $I_2 = \partial \wh V \cap (\R+i\pi)$, and $I = \partial \wh V \backslash (I_1 \cup I_2)$. 
The Markov property of the GFF tells us that $\wh h^\IG|_{\wh V}$ conditioned on $\wh h^\IG|_{\wh U}$ is a mixed boundary GFF with constant boundary conditions on $I_1$ and $I_2$, and Dirichlet boundary conditions on $I$ determined by $\wh h^\IG|_{\wh U}$. By 
\cite[Proposition 2.5 (a)]{ag-disk}, as $N \to \infty$, the law of $\wh h^\IG|_{\wh U_N}$ given $\wh h^\IG|_{\wh U}$ is within $o_N(1)$ in total variation from its unconditioned law. Mapping back to the square domain $D$, this yields the lemma.
\end{proof}

Now we can prove the proposition. 
\begin{proof}[Proof of Proposition~\ref{prop-curve-limit-variant}]
For the single curve case $n=2$, as we send first $\eps \to 0$ then $\delta \to 0$, outside an event of probability $o_\eps(1)$ the segments $\eta^\mathrm{start}$ and $\eta^\eps$ are respectively determined by $h|_U$ and $h|_{U_\delta}$ (Lemma~\ref{lem-B-avoid-y}), and $(h|_U, h|_{U_\delta})$ is $o_\delta(1)$-close in total variation to $(\wt h|_U, h|_{U_\delta})$ where $\wt h$ is an independent copy of $h$ (Lemma~\ref{lem-B-local-indep}). Therefore $(\eta^\mathrm{start}, \eta^\eps)$ is close in total variation to $(\wt \eta^\mathrm{start}, \eta^\eps)$, as desired. 

We now explain the general $n$ regime. We may couple the tuple $(\eta_1, \dots, \eta_{n-1})$ with an appropriate imaginary geometry field $h^\IG$ so that each $\eta_j$ is a flow line of $h^\IG$ with a certain angle. Lemma~\ref{lem-B-avoid-y} applies for each curve $\eta_j$, and Lemma~\ref{lem-B-local-indep} still applies for $h^\IG$, so the same argument applies.
\end{proof}

\section{Brownian motion computations}
In this appendix we carry out some Brownian motion computations which are needed in Section \ref{subsec-MOT}.
\begin{lemma}\label{lem-corr-BM-rate}
	Consider planar Brownian motion $Z=(Z_t)_{t\geq 0}$ with covariance~\eqref{eq-cov} started at $z \in \R_+^2$ and run until it exits $\R_+^2$. Then for some $C >0$ we have $\P[Z\,\,\text{exits in }(0, \eps i)] = (1+o_\eps(1))C \eps^{-\frac{4}{\gamma^2}}$ as $\eps \to 0$. 
\end{lemma}
\begin{proof}
	First perform the shear transformation $z \mapsto \Lambda z$ with $\Lambda = \frac1 {\mathbbm a} \begin{pmatrix}
		\frac1{\sin \theta} & \frac1{\tan \theta} \\
		0 & 1
	\end{pmatrix}$ and $\theta = \frac{\pi \gamma^2}4$, transforming Brownian motion with covariances~\eqref{eq-cov} in $\R_+^2$ to standard Brownian motion in the cone $\{w \: : \: \arg w \in (0, \theta)\}$. Then map $w \mapsto w^{\frac4{\gamma^2}}$ to get Brownian motion in $\bbH$. This maps the interval $(0, i\eps) \subset \partial \R_+^2$ to an interval of length proportional to $\eps^{\frac4{\gamma^2}}$ in $\partial \bbH$ so $\P[Z\,\,\text{exits in }(0, \eps i)] = (C+o_\eps(1))\eps^{\frac4{\gamma^2}}$ for some $C$.
\end{proof}

Recall the measures $\mu^\gamma_{\R_+^2}(z,w)$ defined in Section~\ref{subsec-MOT}.

\begin{lemma}\label{lem-excursion-partition-fn}
	There is a constant $C > 0$ so that for all $\ell, r > 0$ we have 
	\[|\mu^\gamma_{\R_+^2}(\ell, r i)| = C \ell^{\frac4{\gamma^2}-1} r^{\frac4{\gamma^2}-1} (\ell^{\frac4{\gamma^2}} + r^{\frac4{\gamma^2}})^{-2}.\]
\end{lemma}
\begin{proof}
	We will use the boundary Poisson kernel for standard Brownian motion in $\bbH$, given by $H_\bbH(x, y) = \frac1\pi (x-y)^{-2}$ for $x, y \in \R$; this follows from the limit of the bulk Poisson kernel $\lim_{\delta \to 0} \delta^{-1} H_\bbH( x + \delta i, y ) = \lim_{\delta \to 0} \delta^{-1} \cdot \frac{\delta}{\pi ((x-y)^2 + \delta^2)}$. 
	
	First perform the shear transformation $z \mapsto \Lambda z$ with $\Lambda = \frac1 {\mathbbm a} \begin{pmatrix}
		\frac1{\sin \theta} & \frac1{\tan \theta} \\
		0 & 1
	\end{pmatrix}$ and $\theta = \frac{\pi \gamma^2}4$, transforming Brownian motion with covariances~\eqref{eq-cov} in $\R_+^2$ to standard Brownian motion in the cone $\cC_\theta := \{w \: : \: \arg w \in (0, \theta)\}$. Then, writing $p = \frac{1}{\mathbbm a \sin \theta} \ell$ and $q = \frac{1}{\mathbbm a \sin \theta} r$, we have
	\[|\mu^\gamma_{\R_+^2}(\ell, ri)| = C \lim_{\delta \to 0} \lim_{\eps \to 0} \frac1{\delta\eps} \P_{p+i\delta}[\wt E_{q, \eps}],\]
	where $\P_{z}$ corresponds to Brownian motion started at $z$, and $\wt E_{q, \eps}$ is the event that Brownian motion exits $\cC_\theta$ on the boundary interval $[qe^{i\theta}, (q+ \eps) e^{i\theta}]$, and $C>0$ is a constant. 
	
	Now map from $\cC_\theta$ to $\bbH$ by $w \mapsto w^{\frac4{\gamma^2}}$ to see that 
	\[|\mu^\gamma_{\R_+^2}(\ell, ri)| = C \lim_{\delta \to 0} \lim_{\eps \to 0} \frac1{\delta\eps} \P_{p^{\frac4{\gamma^2}}+i\delta p^{\frac4{\gamma^2}-1}}[\wh E_{q^{\frac4{\gamma^2}}, \eps q^{\frac4{\gamma^2}-1}}],\]
	where $\wh E_{q^{\frac4{\gamma^2}}, \eps q^{\frac4{\gamma^2}-1}}$ is the event that Brownian motion exits $\bbH$ on the interval between $-q^{\frac4{\gamma^2}}$ and $-q^{\frac4{\gamma^2}} - \eps q^{\frac4{\gamma^2}-1}$. Taking the limit, we see that 
	$$
	|\mu^\gamma_{\R_+^2}(\ell,ri)| = C p^{\frac4{\gamma^2} - 1} q^{\frac4{\gamma^2} - 1} H_\bbH(p^{\frac4{\gamma^2}}, q^{\frac4{\gamma^2}})  = C \frac{p^{\frac4{\gamma^2} - 1} q^{\frac4{\gamma^2} - 1} }{(p^{\frac4{\gamma^2}}+q^{\frac4{\gamma^2}})^2}.
	$$ 
	Restating this in $\ell$ and $r$ yields the lemma. 
\end{proof}

\bibliographystyle{hmralphaabbrv}
\bibliography{cibib}

\end{document}